\documentclass[12pt]{article}
\usepackage{amsmath,amssymb,amsthm, amscd}

\newtheorem{Lem}{Lemma}
\newtheorem{Prop}[Lem]{Proposition}
\newtheorem{Thm}[Lem]{Theorem}
\newtheorem{Rem}[Lem]{Remark}
\newtheorem{Claim}{Claim}
\newtheorem{Cor}[Lem]{Corollary}
\newtheorem{Example}[Lem]{Example}

\begin{document}
\begin{center}
{\bf Birational geometry for the covering of a nilpotent orbit closure II}
\end{center}
\vspace{0.4cm}

\begin{center}
{\bf Yoshinori Namikawa}
\end{center}
\vspace{0.2cm}

\begin{abstract}
Let $O$ be a nilpotent orbit of a complex semisimple Lie algebra $\mathfrak{g}$ and let 
$\pi: X \to \bar{O}$ be the finite covering associated with the universal covering of $O$. 
In the previous article \cite{Part 1} we have explicitly constructed a {\bf Q}-factorial terminalization 
$\tilde{X}$ of $X$ when $\mathfrak{g}$ is classical. In this article we count how many non-isomorphic {\bf Q}-factorial 
terminalizations $X$ has. We construct the universal Poisson deformation of $\tilde{X}$ over $H^2(\tilde{X}, \mathbf{C})$ 
and look at the action of the Weyl group $W(X)$ on $H^2(\tilde{X}, \mathbf{C})$. The main result is an explicit geometric 
description of $W(X)$. \vspace{0.2cm}

\end{abstract}

\begin{center}
{\bf Introduction}
\end{center} 

Let $O$ be a nilpotent orbit of a complex semisimple Lie algebra $\mathfrak{g}$ and 
let $\pi_0: X^0 \to O$ be an etale covering. Then $\pi_0$ extends to a finite 
covering $\pi: X \to \bar{O}$. A birational projective morphism $f: \tilde{X} \to X$ is called a {\bf Q}-factorial terminalization if $\tilde{X}$ has only {\bf Q}-factorial terminal singularities and $f$ is crepant (i.e. $K_{\tilde X} = f^*K_X$, where $K_X$ and $K_{\tilde X}$ are respectively the 
canonical divisors of $X$ and $\tilde{X}$). Two {\bf Q}-factorial terminalizations $f_1: \tilde{X}_1 \to X$ and $f_2: \tilde{X}_2 \to X$ 
are called {\em isomorphic} if there is an isomorphism $\tau: \tilde{X}_1 \to \tilde{X}_2$ such that $f_1 = f_2 \circ \tau$. 
In the present article we count how many non-isomorphic {\bf Q}-factorial terminalizations $X$ has. The main idea is to construct the universal Poisson deformation $\tilde{\mathcal X} 
\to H^2(\tilde{X}, \mathbf{C})$ of $\tilde{X}$ and look at the action of the Weyl group $W(X)$ 
on the base space $H^2(\tilde{X}, \mathbf{C})$. For the definition and the properties of the Weyl group $W(X)$, see \cite{Na 2}, \S 1 and \cite{Na 4}. In fact, $H^2(\tilde{X}, \mathbf{R})$ is 
divided into a finite number of chambers, and the movable cone $\mathrm{Mov}(\tilde{X})$ is 
a union of some chambers. The Weyl group $W(X)$ acts on $H^2(\tilde{X}, \mathbf{R})$ 
and $\mathrm{Mov}(\tilde{X})$ is a fundamental domain for this $W(X)$-action. 
Each chamber inside $\mathrm{Mov}(\tilde{X})$ corresponds to a {\bf Q}-factorial terminalization of $X$. Hence
$$\sharp \{\mathrm{{\bf Q}-factorial}\; \mathrm{terminalizations}\;  \mathrm{of} \;  X\} =  
\sharp \{\mathrm{chambers}\}/\vert W(X) \vert.$$   
Let us recall how we constructed a {\bf Q}-factorial terminalization $\tilde{X}$ of $X$. 
Let $G$ be a complex semisimple Lie group $G$ with $Lie(G) = \mathfrak{g}$ and $\pi_1(G) = \{1\}$. Let $Q$ be a parabolic subgroup of $G$.  Fix a maximal torus $T$ of $Q$. Then we 
have a Levi decomposition $Q = U\cdot L$. Here $U$ is the unipotent radical and 
$L$ is the Levi subgroup containing $T$. Let $\mathfrak{n}$ be the nilradical of $\mathfrak{q}$. Then $\mathfrak{n}$ coincides with the Lie algebra of $U$. Let $\mathfrak{l}$ be the Lie algebra of $L$. Then we have a direct sum decomposition 
$\mathfrak{q} = \mathfrak{n} \oplus \mathfrak{l}$. Let $O'$ be a nilpotent orbit of $\mathfrak{l}$. Then there is a unique 
nilpotent orbit $O$ of $\mathfrak{g}$ such that $O$ meets $\mathfrak{n} + O'$ in a non-empty Zariski open set of $\mathfrak{n} + O'$. In such a case we say that $O$ is induced from 
$O'$ and write $O = \mathrm{Ind}^{\mathfrak g}_{\mathfrak l}(O')$. 
There is a generically finite map 
$$\mu: G \times^Q(\mathfrak{n} + \bar{O}') \to \bar{O} \:\:\: ([g, z] \mapsto Ad_g(z)),$$ which we call a generalized Springer map. 
Assume that an $L$-equivariant etale cover $(X')^0 \to O'$ is given. Let $X' \to \bar{O}'$ be the associated finite covering.
Then the $Q$-action on $\mathfrak{n} + \bar{O}'$ 
lifts to a Q-action on $\mathfrak{n} \times X'$.
Hence we can make $G \times^Q(\mathfrak{n} 
\times X')$ and get a commutative diagram 
\begin{equation} 
\begin{CD} 
G \times^Q(\mathfrak{n} \times X') @>{\mu'}>> Z \\ 
@V{\pi'}VV @VVV \\ 
G \times^Q (\mathfrak{n} + \bar{O}') @>{\mu}>> \bar{O},       
\end{CD} 
\end{equation}  
where $Z$ is the Stein factorization of $\mu \circ \pi'$. When $Z = X$, we say that the $G$-cover $X^0$ of $O$  
is birationally induced from an $L$-cover $(X')^0$ of $O'$. When $X^0$ is not birationally induced from any other one except itself, 
$X^0$ is called a birationally rigid cover. It is proved in [Ma, Corollary 4.3], [LMM, \S 7] that, for any 
$G$-cover $X^0$ of $O$, there exists a unique $(L, O', (X')^0)$ such that $\mu'$ is birational and 
$(X')^0$ is a birationally rigid cover. Moreover, in such a case, $X'$ has only {\bf Q}-factorial terminal singularities and 
$\mu'$ gives a {\bf Q}-factorial terminalization of $X$. In \cite{Part 1} we have given an explicit algorithm 
for finding such $(L, O', (X')^0)$ when $\mathfrak{g}$ is a classical simple Lie algebra and $X^0$ is the universal covering of 
$O$.   

Therefore, in this article, we start with the situation when a commutative diagram  
\begin{equation} 
\begin{CD} 
G \times^Q(\mathfrak{n} \times X') @>{\mu'}>> X \\ 
@V{\pi'}VV @VVV \\ 
G \times^Q (\mathfrak{n} + \bar{O}') @>{\mu}>> \bar{O},       
\end{CD} 
\end{equation}  
is already given and $\mu'$ is a {\bf Q}-factorial terminalization of $X$.    
Put $$N_G(L, T) := \{g \in G\; \vert \; gTg^{-1} = T, \; gLg^{-1} = L\}.$$
Then $N_G(L, T)$ acts on $\mathfrak{l}$. We are mainly interested in the case 
when $X^0$ is the universal covering of $O$. So 
we assume the following three conditions 

(0)  If we write $O'$ as $L/L^z$ for $z \in O'$, then $(X')^0 = L/(L^z)^0$, where 
$(L^z)^0$ is the identity component of $L^z$. 

(1) $X'$ has only {\bf Q}-factorial terminal singularities. 

(2) $N_G(L, T)$ stabilizes $O'$. 

In fact, (0) is always satisfied when $X^0$ is the universal covering of $O$. 
Moreover, (2) is satisfied for all $(L, O', X')$ constructed in \cite{Part 1}. 

Let $\mathfrak{k}$ be the center of $\mathfrak{l}$. Then the solvable radical 
$\mathfrak{r}$ of $\mathfrak{q}$ is written as $$\mathfrak{r} = \mathfrak{k} \oplus 
\mathfrak{n}.$$ The adjoint $Q$-action stabilizes $\mathfrak{r} + \bar{O}' \subset 
\mathfrak{r} \oplus \mathfrak{l}$ and one can make a space $G \times^Q (\mathfrak{r} + \bar{O}')$. Write $x \in \mathfrak{r} + \bar{O}'$ as $x = x_1 + x_2 + x_3$, where $x_1 \in \mathfrak{k}$, 
$x_2 \in \mathfrak{n}$ and $x_3 \in \bar{O}'$. Define a map 
$$\eta: G \times^Q (\mathfrak{r} + \bar{O}')  \to \mathfrak{k}$$ by $[g, x] \in G \times^Q (\mathfrak{r} + \bar{O}') \mapsto x_1$. Since $G \times^Q(\mathfrak{r} \times X')$ is a finite cover 
of $G \times^Q(\mathfrak{r} + \bar{O}')$, we have a map $G \times^Q(\mathfrak{r} \times X') 
\to \mathfrak{k}$. This map is a flat map with the central fiber $G \times^Q (\mathfrak{n} \times X')$. Moreover, it is actually the universal Poisson deformation of   
$G \times^Q (\mathfrak{n} \times X')$ and $\mathfrak{k}$ is identified with $H^2(G \times^Q (\mathfrak{n} \times X'), \mathbf{C})$. We put $$Z_{\mathfrak{l}, X'} := 
\mathrm{Spec}\; \Gamma (G \times^Q (\mathfrak{r} \times X'), \mathcal{O}).$$
As the subscript indicates, $Z_{\mathfrak{l}, X'}$ only depends on the Levi part $\mathfrak{l}$ 
and $X'$.  
By definition, there is a map $Z_{\mathfrak{l}, X'} \to \mathfrak{k}$ whose central fiber is 
$X$. By the general theory of the Poisson deformation 
of a conical symplectic variety, the Weyl group $W(X)$ for $X$ acts on both $\mathfrak{k}$ and 
$Z_{\mathfrak{l}, X'}$ so that the induced map $Z_{\mathfrak{l}, X'}/W(X) \to \mathfrak{k}/W(X)$ gives the universal Poisson deformation of $X$: 
\begin{equation} 
\begin{CD} 
G \times^Q(\mathfrak{r} \times X') @>>> Z_{\mathfrak{l}, X'}/W(X) \\ 
@VVV @VVV \\ 
\mathfrak{k} @>>> \mathfrak{k}/W(X)     
\end{CD} 
\end{equation}
The main result of this article is an explicit description of the $W(X)$ (Theorem \ref{Theorem 9}). 
In \cite{Na 1} we already observed that a certain finite group $W'$ acts on $\mathfrak{k}$. 
The group $W'$ is defined as $W' := N_G(L,T)/N_L(T)$, where $$N_L(T) := \{g \in L\; \vert \; 
gTg^{-1} = T\}.$$ Here we put $$Y_{\mathfrak{l}, O'} := \mathrm{Spec}\; \Gamma 
(G \times^Q (\mathfrak{r} + O'), \mathcal{O}).$$ The map $Z_{\mathfrak{l}, X'} \to \mathfrak{k}$ factorizes as $$Z_{\mathfrak{l}, X'} \to Y_{\mathfrak{l}, O'} \to \mathfrak{k}.$$
Then $\mathrm{Aut}(Z_{\mathfrak{l}, X'}/Y_{\mathfrak{l}, O'})$ is isomorphic 
to $\mathrm{Aut}((X')^0/{O}')$.   
The group $W'$ acts on $Y_{\mathfrak{l}, O'}$, but does not act on $Z_{\mathfrak{l}, X'}$ 
in general.  Instead of $W'$ we consider a new group $\tilde{W}'$ obtained as an extension 
of $W'$ by $\mathrm{Aut}((X')^0/O')$: 
$$ 1 \to  \mathrm{Aut}((X')^0/O') \to \tilde{W}' \to W' \to 1$$ 
Then $\tilde{W}'$ acts on $Z_{\mathfrak{l}, X'}$ in such a way that the map 
$Z_{\mathfrak{l}, X'} \to \mathfrak{k}$ is $\tilde{W}'$-equivariant. 
Since $W'$ fixes the origin $0 \in \mathfrak{k}$, an element of $\tilde{W}'$ 
acts on the central fiber $X$ of $Z_{\mathfrak{l}, X'} \to \mathfrak{k}$ as an element 
of $\mathrm{Aut}(X/\bar{O})$. In this way we have a homomorphism 
$$\rho_X: \tilde{W}' \to \mathrm{Aut}(X/\bar{O}).$$ We then define 
$$W_X := \mathrm{Ker}(\rho_X).$$ 
This $W_X$ equals the Weyl group $W(X)$ (Theorem \ref{Theorem 9}). We remark that Losev \cite{Lo}, Proposition 4.6 has given a similar description of $W(X)$ when $X'$ is the normalization of $\bar{O'}$. In this case we do not need 
to take $\tilde{W}'$.

In \S 2 we first show that $\rho_X$ is surjective.  This is enough to calculate $\vert W_X \vert$ and the number of 
non-isomorphic {\bf Q}-factorial terminalizations of $X$ (cf. Theorem \ref{Theorem 14}). 
But we need more information on $\rho_X$ to know the structure of $W_X$.      
When $\mathfrak{g}$ is a classical simple Lie algebra, we study $\rho_X(\tilde{w})$ for an element $\tilde{w} \in \tilde{W}'$ 
such that its image $w \in W'$ is a reflection.    
In \S 3, we actually count the numbers of {\bf Q}-factorial terminalizations of $X$ in two examples as an application of the theory developed in the previous sections.
\vspace{0.2cm}

{\em Acknowledgement}. The author thanks Ivan Losev for letting him know Proposition 4.7 
in \cite{Lo} after the first version of this paper (arXiv:1912.01729) was written. 
\vspace{0.2cm}

\begin{center}
\S\:  {\bf 1}. Universal Poisson deformations of {\bf Q}-factorial terminalizations and Weyl groups
\end{center}

Let $G$ be a complex semisimple Lie group and let $Q$ be a parabolic subgroup 
of $G$. We take a maximal torus $T$ of $G$ so that $T \subset Q$. 
Let $Q = U \cdot L$ be the Levi decomposition with $T \subset L$. 
Let $\mathfrak{n}$ be the nilradical of $\mathfrak{q}$. Then $\mathfrak{n}$ coincides with the Lie algebra of $U$. Let $\mathfrak{l}$ be the Lie algebra of $L$. Then we have a direct sum decomposition 
$\mathfrak{q} = \mathfrak{n} \oplus \mathfrak{l}$. 
Let $O' \subset \mathfrak{l}$ be a nilpotent orbit. Then we have a generalized 
Springer map $$\mu: G \times^Q(\mathfrak{n} + \bar{O}') \to \mathfrak{g}, \:\: 
[g, x] \mapsto Ad_g(x).$$ The image of $\mu$ is the closure $\bar{O}$ of a certain nilpotent orbit 
$O$ of $\mathfrak{g}$. We write $O'$ as $L/L^z$ for $z \in O'$. Put $(X')^0 := L/(L^z)^0$ and 
consider an etale covering $(X')^0 \to O'$.  Let $X' \to \bar{O}'$ be the associated finite covering. Then we have a finite covering 
$\mathfrak{n} \times X' \to \mathfrak{n} \times \bar{O}' (= \mathfrak{n} + O')$. 
The $Q$-action on $\mathfrak{n} + \bar{O}'$ lifts to a $Q$-action on $\mathfrak{n} \times  
X'$. 
We assume  

(1) $X'$ has only {\bf Q}-factorial terminal singularities. 

We can make $G \times^Q (\mathfrak{n} \times X')$ and get a commutative 
diagram 
\begin{equation} 
\begin{CD} 
G \times^Q(\mathfrak{n} \times X') @>{\mu'}>> X \\ 
@V{\pi'}VV @V{\pi}VV \\ 
G \times^Q (\mathfrak{n} + \bar{O}') @>{\mu}>> \bar{O},       
\end{CD} 
\end{equation}  
where $X$ is the Stein factorization of $\mu \circ \pi'$. 
$X$ is a conical symplectic variety and $\mu'$ gives a {\bf Q}-factorial terminalization 
of $X$.  

Let $\mathfrak{k}$ be the center of $\mathfrak{l}$. We put $\mathfrak{t} := 
Lie(T)$. Then the root system $\Phi$ is determined from $(\mathfrak{g}, \mathfrak{t})$. 
The quotient ${\mathfrak l}_{\mathrm{s.s}} := \mathfrak{l}/\mathfrak{k}$ is a 
semisimple Lie algebra, and $\mathfrak{t}/\mathfrak{k}$ is a Cartan subalgebra of 
$\mathfrak{l}_{\mathrm{s.s}}$.   
The root system $\Phi_{\mathfrak{l}_{\mathrm{s.s}}}$ for $(\mathfrak{l}_{\mathrm{s.s}},  \mathfrak{t}/\mathfrak{k})$ is a root subsystem of $\Phi$. We have 
$$\mathfrak{k} = \{h \in \mathfrak{t}\; \vert \; \alpha (h) = 0,\; \forall \alpha \in 
\Phi_{\mathfrak l}\}.$$ 
We define $$\mathfrak{k}^{reg} := \{h \in \mathfrak{k}\; \vert \; 
\alpha (h) \ne 0\; \forall \alpha \in \Phi \setminus \Phi_{\mathfrak{l}_{\mathrm{s.s}}}\}.$$   
The solvable radical 
$\mathfrak{r}$ of $\mathfrak{q}$ is written as $$\mathfrak{r} = \mathfrak{k} \oplus 
\mathfrak{n}.$$ The adjoint $Q$-action stabilizes $\mathfrak{r} + \bar{O}' \subset 
\mathfrak{r} \oplus \mathfrak{l}$ and one can make a space $G \times^Q (\mathfrak{r} + \bar{O}')$. Write $x \in \mathfrak{r} + \bar{O}'$ as $x = x_1 + x_2 + x_3$, where $x_1 \in \mathfrak{k}$, 
$x_2 \in \mathfrak{n}$ and $x_3 \in \bar{O}'$. Define a map 
$$\eta: G \times^Q (\mathfrak{r} + \bar{O}')  \to \mathfrak{k}$$ by $[g, x] \in G \times^Q (\mathfrak{r} + \bar{O}') \mapsto x_1$. This is a well-defined map; in fact, 
for $q \in Q$, we have $(Ad_q(x_1))_1 = x_1$, $Ad_q(x_2) \in \mathfrak{n}$ and $Ad_q(x_3) \in 
\mathfrak{n} + \bar{O}'$. By definition $\eta^{-1}(0) = G \times^Q (\mathfrak{n} + \bar{O}')$ 
and $\eta$ gives a flat deformation of $G \times^Q (\mathfrak{n} + \bar{O}')$ over 
$\mathfrak{k}$. Moreover, the $Q$-action on $\mathfrak{r} + \bar{O}'$ lifts to a $Q$-action 
on $\mathfrak{r} \times X'$. Then $G \times^Q (\mathfrak{r} \times X')$ is a finite covering 
of $G \times^Q(\mathfrak{r} + \bar{O}')$. We have a commutative diagram 
\begin{equation} 
\begin{CD}
G \times^Q(\mathfrak{n} \times X') @>>> G \times^Q(\mathfrak{r} \times X') \\
@V{\pi'}VV @V{\Pi'}VV \\   
G \times^Q(\mathfrak{n} + \bar{O}') @>>> G \times^Q(\mathfrak{r} + \bar{O}') \\ 
@VVV @V{\eta}VV \\ 
\{0\} @>>> \mathfrak{k},    
\end{CD} 
\end{equation} 
where $\eta \circ \Pi'$ gives a flat deformation of $G \times^Q(\mathfrak{n} \times X')$ 
over $\mathfrak{k}$.
Note that $\mathfrak{k}$ is naturally identified with $H^2(G \times^Q(\mathfrak{r} \times X'), 
\mathbf{C})$. In fact, $H^2(G \times^Q(\mathfrak{r} \times X'), 
\mathbf{C}) = \mathrm{Pic}(G \times^Q(\mathfrak{r} \times X')) \otimes_{\mathbf Z}\mathbf{C}$ 
and we have isomorphisms 
$$\mathrm{Hom}_{alg.gp}(L, \mathbf{C}^*) \cong \mathrm{Pic}(G/Q) \cong  
\mathrm{Pic}(G \times^Q(\mathfrak{r} \times X'))$$
The isomorphism $\mathfrak{t} \to \mathfrak{t}^*$ determined by the Killing form 
sends $\mathfrak{k} \subset \mathfrak{t}$ isomorphically onto the subspace  
$\mathrm{Hom}_{alg.gp}(L, \mathbf{C}^*)\otimes_{\mathbf Z}\mathbf{C} \subset 
\mathfrak{t}^*$. \vspace{0.2cm}

Let us consider a fiber of $\eta$. 
For $t \in \mathfrak{k}$, 
$\eta^{-1}(t) = G \times^Q (t + \mathfrak{n} + \bar{O}')$. \vspace{0.2cm} 

\begin{Lem}\label{Lemma 1} (\cite{Na 3}, Lemma (3.2)) 
For $t \in \mathfrak{k}^{reg}$, any orbit of 
the $Q$-variety $t + \mathfrak{n} + \bar{O}'$ is of the form $Q(t + y)$ 
with $y \in \bar{O'}$.
\end{Lem}
\vspace{0.2cm}

Let $W$ be the Weyl group of 
$\mathfrak{g}$ with respect to $\mathfrak{t}$. $W$ acts on $\mathfrak{t}$. Let $\mathfrak{g}//G$ be the adjoint quotient of $\mathfrak{g}$. 
The map $\mathfrak{t}/W \to \mathfrak{g}//G$ is an isomorphism by Chevalley's restriction theorem.    
We consider the adjoint quotient map 
$\chi : \mathfrak{g} \to \mathfrak{g}//G \cong \mathfrak{t}/W$. Let $G \cdot (\mathfrak{r} + \bar{O}') \subset \mathfrak{g}$ be the space of all $G$-translates of the elements of $\mathfrak{r} + \bar{O}'$.
We denote by the same $\chi$ the restriction of $\chi$ to $G \cdot (\mathfrak{r} + \bar{O}')$.  
By \cite{Na 3}, p. 21, the following diagram commutes 
\begin{equation} 
\begin{CD} 
G \times^Q(\mathfrak{r} + \bar{O}') @>>> G\cdot (\mathfrak{r} + \bar{O}') \\ 
@V{\eta}VV @V{\chi}VV \\ 
\mathfrak{k} @>{\iota}>> \mathfrak{t}/W,       
\end{CD} 
\end{equation}  
where the top horizontal map is defined by $[g, x] \mapsto Ad_g(x)$ and $\iota$ is the composite 
of the natural inclusion $\mathfrak{k} \to \mathfrak{t}$ and the quotient map 
$\mathfrak{t} \to \mathfrak{t}/W$.  
By Lemma \ref{Lemma 1}, $$G\cdot(\mathfrak{r} + \bar{O}') = 
\overline{G\cdot (\mathfrak{k}^{reg} + \bar{O'})}.$$ 
Note that $\mathfrak{k}$ only depends on the Levi part $\mathfrak{l}$ of 
$\mathfrak{q}$. Hence $G\cdot(\mathfrak{r} + \bar{O}')$ depends only on 
$\mathfrak{l}$ and $O'$. 

The commutative diagram above induces a commutative diagram 
\begin{equation} 
\begin{CD} 
G \times^Q(\mathfrak{r} + \bar{O}') @>>> \mathfrak{k} \times_{\mathfrak{t}/W} G\cdot (\mathfrak{r} + \bar{O}') \\ 
@V{\eta}VV @VVV \\ 
\mathfrak{k} @>{id}>> \mathfrak{k}       
\end{CD} 
\end{equation}

We put $N_G(T) := \{g \in G \: \vert \: gTg^{-1} = T\}$, $N_G(L,T) := \{g \in N_G(T) \: \vert 
\: gLg^{-1} = L\}$ and $N_L(T) := \{g \in L \: \vert \: gTg^{-1} = T\}$. 
Then $T \subset N_L(T) \subset N_G(L,T) \subset N_G(T)$. The Weyl group $W$ for $G$ is defined as $W := N_G(T)/T$ and the Weyl group $W(L)$ for $L$ is defined as $W(L) := N_L(T)/T$. We have inclusions $$W(L) \subset N_G(L,T)/T \subset W.$$ 
Define 
\begin{equation}\label{W'}
W' := N_G(L,T)/N_L(T)
\end{equation}
Obviously $W'$ is the quotient of $N_G(L,T)/T$ by 
$W(L)$. Then $W'$ acts on $\mathfrak{k}$ and 
the map $\iota : \mathfrak{k} \to \mathfrak{t}/W$ factorizes as $\mathfrak{k} 
\to \mathfrak{k}/W' \to \mathfrak{t}/W$. Here $\mathfrak{k}/W'$ is the normalization 
of $\iota (\mathfrak{k})$.  Let $G \cdot (\mathfrak{r} + \bar{O}')^n$ be the normalization 
of $G \cdot (\mathfrak{r} + \bar{O}')$. Then the map $G\cdot (\mathfrak{r} + \bar{O}')^n 
\to \mathfrak{t}/W$ factors through $\mathfrak{k}/W'$. Hence  
$$(\mathfrak{k} \times_{\mathfrak{t}/W} G\cdot (\mathfrak{r} + \bar{O}'))^n 
= (\mathfrak{k} \times_{\mathfrak{k}/W'} G\cdot (\mathfrak{r} + \bar{O}')^n)^n.$$ 
We denote this variety by $Y_{\mathfrak{l}, O'}$. 
As the notation indicates,   
$Y_{\mathfrak{l}, O'}$ depends only on $\mathfrak{l}$ and $O'$. 
Let $W'$ act on $\mathfrak{k}$ in the usual way and act trivially on $G\cdot (\mathfrak{r} + \bar{O}')^n$. Then $W'$ acts on $\mathfrak{k} \times_{\mathfrak{k}/W'} G\cdot (\mathfrak{r} + \bar{O}')^n$; hence acts on its normalization $Y_{\mathfrak{l}, O'}$. Let $\tilde{O}'$ be 
the normalization of $\bar{O}'$. Then the normalization of $G \times^Q(\mathfrak{r} + \bar{O}')$ is $G \times^Q(\mathfrak{r} \times \tilde{O}')$. 
There is a natural map $$G \times^Q(\mathfrak{r} \times \tilde{O}') \to Y_{\mathfrak{l}, O'}.$$
By definition $N_G(L,T)$ acts on $\mathfrak{l}$. We assume that 
\vspace{0.15cm}

(2) $N_G(L,T)$ stabilizes $O'$. 
\vspace{0.15cm}

Then the map is a birational projective morphism by \cite{Na 3}, Proposition 3.6.     
Here let us consider the composite $G \times^Q(\mathfrak{r} \times X') \to 
G \times^Q(\mathfrak{r} \times \tilde{O}') \to Y_{\mathfrak{l}, O'}$ and take its Stein factorization 
$Z_{{\mathfrak l}, X'}$. Then we get a commutative diagram 
\begin{equation} 
\begin{CD} 
G \times^Q(\mathfrak{r} \times X') @>{\mathfrak{m}'_Q}>> Z_{\mathfrak{l}, X'} \\ 
@VVV @VVV \\ 
G \times^Q(\mathfrak{r} \times \tilde{O}') @>{\mathfrak{m}_Q}>> Y_{\mathfrak{l}, O'}
\end{CD} 
\end{equation}
We claim that $Z_{\mathfrak{l}, X'}$ depends only on $\mathfrak{l}$ and $X'$ as the 
subscript indicates. More precisely, 
if we start with another parabolic subgroup $Q'$ of $G$ with Levi part $L$ and with 
the same $X' \to \bar{O}'$, then we can make a similar diagram for 
$G \times^{Q'} (\mathfrak{r}' \times X')$. We claim that the same $Z_{\mathfrak{l}, X'}$ is obtained as the Stein factorization. The claim follows from the next lemma. In fact, by the lemma, 
$G \times^Q (\mathfrak{r} +X')$ and $G \times^{Q'}(\mathfrak{r}' \times X')$ have 
the same function field, and $Z_{\mathfrak{l}, X'}$ is nothing but the normalization of 
$Y_{\mathfrak{l}, O'}$ in this field. \newpage  

\begin{Lem}\label{Lemma 2} Let $Q$ and $Q'$ be parabolic subgroups of $G$ containing 
$L$ as the Levi parts. Then there is a commutative diagram  
$$ G \times^Q(\mathfrak{r} \times X') ---\to G \times^{Q'}(\mathfrak{r}' \times X')$$ 
$$  \downarrow \hspace{3.0cm} \downarrow $$
$$ G \times^Q(\mathfrak{r} + \bar{O}') ---\to G \times^{Q'}(\mathfrak{r}' + \bar{O}') $$
$$  \downarrow \hspace{3.0cm} \downarrow $$
$$ \mathfrak{k} \times_{\mathfrak{t}/W} G\cdot (\mathfrak{r} + \bar{O}') = \mathfrak{k} \times_{\mathfrak{t}/W} G\cdot (\mathfrak{r} + \bar{O}') $$
$$  \downarrow \hspace{3.0cm} \downarrow $$
$$  \mathfrak{k}  \hspace{1.3cm}  = \hspace{1.3cm}  \mathfrak{k}$$
where the 1-st and the 2-nd horizontal maps are both birational maps. 
\end{Lem} \vspace{0.2cm}

{\em Proof}. The fiber of  the map $G \times^Q(\mathfrak{r} + \bar{O}')  \to \mathfrak{k}$ over $t \in \mathfrak{k}^{reg}$ is $G \times^Q(t + \mathfrak{n} + \bar{O}')$.
By Lemma \ref{Lemma 1}, $$(*) \:\:\:\: G \times^Q(t + \mathfrak{n} + \bar{O}') = G \times^Q (\cup_{y \in \bar{O}'}
Q\cdot(t + y)) = \cup_{y \in \bar{O}'}G \times^Q Q\cdot(t + y) $$   
Denote by $Q_{t+y}$ the stabilizer group of $t + y$. Note that $Q_{t+y} \subset L$. In fact, 
assume that $q \in Q_{t+y}$. Put $z = t + x$. Since $t \in \mathfrak{k}$ and $y \in \mathfrak{l}$, we have $[t, y] = 0$; hence, $z = t + y$ is the Jordan-Chevalley decomposition of $z$. Then $Ad_q(z) = Ad_q(t) + Ad_q(y)$ is also the Jordan-Chevalley decomposition of 
$z = Ad_q(z)$. This means that $Ad_q(t) = t$ and $Ad_q(x) = x$. In particular, $q$ is contained in the centralizer $Z(t)$ of $t$ (in $Q$). $Z(t) = L$ because $t \in \mathfrak{k}^{reg}$. Therefore $Q_{t+y} \subset L$. 
We have $$Q\cdot(t+y) = Q \times^{Q_{t+y}}
(t+y) = Q \times^L L \times^{Q_{t+y}} (t+y) = Q \times^L L\cdot(t+y).$$  
The right hand side of $(*)$ can be written as 
$$\cup_{y \in \bar{O}'}G \times^Q Q\cdot(t + y) = \cup_{y \in \bar{O}'} G \times^Q 
(Q \times^L L\cdot(t+y)) = \cup_{y \in \bar{O}'} G \times^L L\cdot(t+y) = 
G \times^L (t + \bar{O}').$$ Changing $Q$ by $Q'$, we similarly get 
$$G \times^{Q'}(t + \mathfrak{n}' + \bar{O}') = G \times^L (t + \bar{O}').$$
As a consequence we have $$G \times^Q(t + \mathfrak{n} + \bar{O}') =  
G \times^{Q'}(t + \mathfrak{n}' + \bar{O}').$$ This means that the birational map 
$G \times^Q(\mathfrak{r} + \bar{O}') ---\to G \times^{Q'}(\mathfrak{r}' + \bar{O}')$ is 
an isomorphism over any point of $\mathfrak{k}^{reg}$. 

We next look at a fiber of the map $G \times^Q(\mathfrak{r} \times X')  \to \mathfrak{k}$
over $t \in \mathfrak{k}^{reg}$. The fibre is  $G \times^Q((t + \mathfrak{n}) \times X')$. 
We can write $$(t + \mathfrak{n}) \times X' = \cup_{\tilde{y} \in X'} Q \cdot (t , \tilde{y})$$  because $t + \mathfrak{n} + \bar{O}' = 
\cup_{y \in \bar{O}'} Q \cdot (t + y)$.       
For $\tilde{y} \in X'$, we denote by $y \in \bar{O}'$ its image by the map $X' \to \bar{O}$. 
Let $Q_{(t , \tilde{y})}$ be the stabilizer group of $(t, \tilde{y})$. Then $Q_{(t, \tilde{y})} 
\subset Q_{t + y}$. In particular, $Q_{(t, \tilde{y})} \subset L$. Then, by the same argument 
as above, we have $$G \times^Q((t + \mathfrak{n}) \times X') = G \times^L (t \times X').$$  
Since the right hand side depends only on $L$ and $X'$, we conclude that 
$$G \times^Q((t + \mathfrak{n}) \times X') =  
G \times^{Q'}((t + \mathfrak{n}') \times X').$$
Therefore we have a birational map $$G \times^Q(\mathfrak{r} \times X') ---\to G \times^{Q'}(\mathfrak{r}' \times X'),$$ which is an isomorphism over any point of $\mathfrak{k}^{reg}$. $\square$

By definition the map $\mathfrak{m}'_Q : G \times^Q (\mathfrak{r} \times X') \to Z_{\mathfrak{l}, X'}$ is a birational map. For $t \in \mathfrak{k}$ it induces a map $\mathfrak{m}'_{Q,t}: G \times^Q ((t + \mathfrak{n}) \times X') \to Z_{\mathfrak{l}, X', t}$ 
between fibers.  \vspace{0.2cm}

\begin{Lem}\label{Lemma 3} The map $\mathfrak{m}'_{Q,t}$ is a birational morphism for every  
$t \in \mathfrak{k}$. For $t \in \mathfrak{k}^{reg}$, the map $\mathfrak{m}'_{Q,t}$ is an 
isomorphism. In particular, $\mathfrak{m}'_Q$ is an isomorphism in codimension one. \end{Lem} 

{\em Proof}. When $t = 0$, the map $\mathfrak{m}_Q: G \times^Q(\mathfrak{r} \times \tilde{O}') 
\to Y_{\mathfrak{l}, \bar{O}'}$ induces a birational morphism $\mathfrak{m}_{Q,0}: G \times^Q (\mathfrak{n} \times \tilde{O}') \to (G \cdot (\mathfrak{n} + \bar{O}'))^s$. Here $(G \cdot (\mathfrak{n} + \bar{O}'))^s$ is the Stein factorization of the generalized Springer map 
$\mu: G \times^Q (\mathfrak{n} \times \tilde{O}') \to G \cdot (\mathfrak{n} + \bar{O}')$. 
We have a commutative diagram 
\begin{equation} 
\begin{CD} 
G \times^Q(\mathfrak{n} \times X') @>{\mathfrak{m}'_{Q,0}}>> Z_{\mathfrak{l}, X',0} \\ 
@V{\pi'}VV @VVV \\ 
G \times^Q (\mathfrak{n} \times \tilde{O}') @>{\mathfrak{m}_{Q,0}}>> (G \cdot (\mathfrak{n} + \bar{O}'))^s     
\end{CD} 
\end{equation}
The composite $\mathfrak{m}_{Q,0} \circ \pi'$ is a generically finite surjective map.  
Since $\mathfrak{m}'_Q$ has only connected fibers, $\mathfrak{m}'_{Q,0}$ also has connected 
fibers. This means that $\mathfrak{m}'_{Q,0}$ is a birational map. The space 
$\mathfrak{r} + \bar{O}'$ has a natural scaling $\mathbf{C}^*$-action.  Let $s: \mathbf{C}^* 
\to \mathrm{Aut}(\mathfrak{r} + \bar{O}')$ be the corresponding homomorphism. We define 
a new $\mathbf{C}^*$-action, say $\sigma$, on $\mathfrak{r} + \bar{O}'$ by the composite 
of the map $\mathbf{C}^* \to \mathbf{C}^* \:\: (t \to t^2)$ and $s$. Then $\sigma$ lifts to 
a $\mathbf{C}^*$-action on $\mathfrak{r} \times X'$. Thus $G \times^Q (\mathfrak{r} \times X')$ has 
a $\mathbf{C}^*$-action. We also introduce a $\mathbf{C}^*$-action on $\mathfrak{k}$ 
by the composite of $\mathbf{C}^* \to \mathbf{C}^* \:\: (t \to t^2)$ and 
the scaling action $\mathbf{C}^* \to \mathrm{Aut}(\mathfrak{k})$. Then the map 
$G \times^Q (\mathfrak{r} \times X') \to \mathfrak{k}$ is $\mathbf{C}^*$-equivariant.
By this $\mathbf{C}^*$-action we see that $\mathfrak{m}'_{Q,t}$ is a birational map for 
every $t \in \mathfrak{k}$.    

We next look at the fibers over $t \in \mathfrak{k}^{reg}$.  
We then have a commutative diagram  
\begin{equation} 
\begin{CD} 
G \times^Q((t + \mathfrak{n}) \times X') @>{\mathfrak{m}'_{Q,t}}>> Z_{\mathfrak{l}, X',t} \\ 
@VVV @VVV \\ 
G \times^Q (t + \mathfrak{n} \times \tilde{O}') @>{\mathfrak{m}_{Q,t}}>> Y_{\mathfrak{l}, O', t}     
\end{CD} 
\end{equation}

By \cite{Na 3}, Proposition 3.6, $\mathfrak{m}_{Q,t}$ is an isomorphism. Moreover, by Lemma \ref{Lemma 2}, 
$G \times^Q (t + \mathfrak{n} + \bar{O}') = G \times^L (t + \bar{O}')$. Taking the normalizations of both sides, we have $G \times^Q ((t + \mathfrak{n}) \times \tilde{O}') = G \times^L (t \times \tilde{O}')$. Moreover we have  
$G \times^Q ((t + \mathfrak{n}) \times X') = G \times^L (t \times X')$. 
Hence $Z_{\mathfrak{l}, O', t}$ is the Stein factorization of the finite surjective map 
$G \times^L (t \times X') \to G \times^L (t \times \tilde{O}')$. Then we see that $\mathfrak{m}'_{Q,t}$ 
is an isomorphism. $\square$  \vspace{0.2cm}

\begin{Prop}\label{Proposition 4} The maps $Z_{\mathfrak{l}, X'} \to \mathfrak{k}$ and $Y_{\mathfrak{l}, O'} \to \mathfrak{k}$ are both flat. Their central fibers $Z_{\mathfrak{l}, X', 0}$ and $Y_{\mathfrak{l}, O',0}$ are both normal varieties. Moreover, $Z_{\mathfrak{l}, X', 0} = X$. \end{Prop}
\vspace{0.2cm}  

{\em Proof}. This is similar to \cite{Na 2}, Lemma 2.1. $\square$ 
\vspace{0.2cm}

Let $\tilde{O}'$ be the normalization of $\bar{O}'$. 
We first introduce a Poisson structure on $G \times^Q (\mathfrak{r} \times \tilde{O}')$ over $\mathfrak{k}$. In order to do this, we consider the $G$-equivariant map 
$$\mu_t: G \times^Q (t + \mathfrak{n} \times \tilde{O}') \to \mathfrak{g}, \;\;\; [g, t + x + y] \mapsto 
Ad_g(t + x + y)$$ for each $t \in \mathfrak{k}$. Then the image of this map coincides with the closure of an adjoint 
orbit on $\mathfrak{g}$, say $O_{\mu_t}$. The pull-back $\mu_t^*{\omega}_t$ of the Kostant-Kirillov form $\omega_t$ on $O_{\mu_t}$ gives a symplectic 2-form on $G \times^Q (t + \mathfrak{n} + {O}')$ and $\{\mu_t^*\omega_t\}_{t \in \mathfrak{k}}$ determines a relative symplectic 2-form 
on $G \times^Q(\mathfrak{r} + O')$ over $\mathfrak{k}$ (cf. \cite{Na 3}, Proposition 4.2). 
This relative symplectic 2-form determines a Poisson structure on $G \times^Q(\mathfrak{r} + O')$ and the Poisson structure uniquely extends to a Poisson structure on $G \times^Q ( \mathfrak{r} \times \tilde{O}')$. We next introduce a Poisson structure on $G \times^Q(\mathfrak{r} 
\times X')$ over $\mathfrak{k}$. There is an etale finite map $$G \times^Q(\mathfrak{r} \times (X')^0) 
\to G \times^Q(\mathfrak{r} + O').$$ If we pull back 
the relative symplectic 2-form on $G \times^Q(\mathfrak{r} + O')$ by this map, we get 
a relative symplectic 2-form on $G \times^Q(\mathfrak{r} \times (X')^0)$. This relative symplectic 
2-form determines a Poisson structure on $G \times^Q(\mathfrak{r} \times (X')^0)$. Moreover, 
such a Poisson structure uniquely extends to a Poisson structure on $G \times^Q(\mathfrak{r} \times X')$. Finally these Poisson structures respectively induce Poisson structures on $Y_{\mathfrak{l}, O'}$ 
and $Z_{\mathfrak{l}, X'}$ because 
$$Y_{\mathfrak{l}, O'} = \mathrm{Spec}\; 
\Gamma (G \times^Q(\mathfrak{r} + O'), \mathcal{O}), \;\;\; 
Z_{\mathfrak{l}, X'} = \mathrm{Spec}\; 
\Gamma (G \times^Q(\mathfrak{r} \times X'), \mathcal{O}).$$ 
Therefore we have: 

\begin{Prop}\label{Proposition 5} \em $G \times^Q (\mathfrak{r} \times X') \to \mathfrak{k}$ is a Poisson deformation of $G \times^Q(\mathfrak{n} \times X')$. $Z_{\mathfrak{l}, X'} \to \mathfrak{k}$ is a 
Poisson deformation of $X$. \end{Prop}
\vspace{0.2cm}

By assumption, $N_G(L, T)$ acts on $O'$. For each $w \in N_G(L, T)$, the automorphism  $Ad_w: O' \to O'$ lifts to an automorphism $\theta_w: (X')^0 \to (X')^0$. In fact, take a point 
$z \in O'$ and write $O' = L/L^z$. By the assumption $Ad_w(z) \in O'$. $Ad_w: L \to L$ sends 
$L^z$ onto $L^{Ad_w(z)}$. Hence $Ad_w$ sends the identity component $(L^z)^0$ of $L^z$ to 
the identity component $(L^{Ad_w(z)})^0$ of $L^{Ad_w(z)}$. In particular, $Ad_w$ induces 
an isomorphism $L/(L^z)^0 \to L/(L^{Ad_w(z)})^0$. Recall that $(X')^0 = L/(L^z)^0$ as 
an $L$-variety. Take points $\tilde{z}, \: \widetilde{Ad_w(z)} \in (X')^0$ lying over $z, \: Ad_w(z) 
\in O'$. Then the isomorphism determines an automorphism $\theta_w: (X')^0 \to (X')^0$ which 
sends $\tilde{z}$ to $\widetilde{Ad_w(z)}$.    
Note that, for a covering transformation $\sigma \in \mathrm{Aut}((X')^0/O')$, the 
composite $\sigma \circ \theta_w$ also gives a lifting of $Ad_w$. For each $w \in N_G(L, T)$, 
we consider all liftings of $\theta_w$ of $Ad_w$. Then these form a group $\widetilde{N_G(L, T)}$ and one has an exact sequence 
$$ 1 \to \mathrm{Aut}((X')^0/O') \to \widetilde{N_G(L, T)} \to N_G(L, T) \to 1.$$
$\widetilde{N_G(L, T)}$ acts on $O'$ by means of the map $\widetilde{N_G(L, T)} \to N_G(L, T)$.  
In general, the $N_G(L,T)$-action on $O'$ does not lift to an action on 
$(X')^0$, but the $\widetilde{N_G(L, T)}$-action on $O'$ lifts to an action on $(X')^0$. 
This $\widetilde{N_G(L, T)}$-action extends to an action on $X'$.  
For $\tilde{w} \in \widetilde{N_G(L, T)}$ we denote by $w \in N_G(L, T)$ its image by the 
map $\widetilde{N_G(L, T)} \to N_G(L, T)$.  We define an isomorphism 
$$\phi_{\tilde w}: G \times^Q (\mathfrak{r} \times X') \to G \times^{Ad_w(Q)}(Ad_w(\mathfrak{r}) \times X'), \:\:\: [g, (y, \tilde{z})] \mapsto [gw^{-1}, (Ad_w(y), \tilde{w}(\tilde{z}))].$$
As is already remarked, $N_G(L, T)$ (or its quotient group $W'$) acts on $Y_{\mathfrak{l}, O'}$.   
We have a commutative diagram 
\begin{equation} 
\begin{CD} 
G \times^Q (\mathfrak{r} \times X') @>{\phi_{\tilde w}}>> G \times^{Ad_w(Q)}(Ad_w(\mathfrak{r}) \times X') \\ 
@VVV @VVV \\ 
Y_{\mathfrak{l}, O'} @>{w}>> Y_{\mathfrak{l}, O'}     
\end{CD} 
\end{equation}     
By Lemma \ref{Lemma 2}, $G \times^Q (\mathfrak{r} \times X')$ and $G \times^{Ad_w(Q)}(Ad_w(\mathfrak{r}) \times X')$ are birational (not by $\phi_{\tilde w}$, but as varieties over $Y_{\mathfrak{k}, O'}$), hence they have the same function field. By $\phi_{\tilde w}$, the group $\widetilde{N_G(L, T)}$ acts on the function field. 
We have defined $Z_{\mathfrak{l}, X'}$ as the normalization of $Y_{\mathfrak{l}, O'}$ in this 
field. Hence $\widetilde{N_G(L, T)}$ acts on $Z_{\mathfrak{l}, X'}$ and the map $Z_{\mathfrak{l}, X'} \to \mathfrak{k}$ is $\widetilde{N_G(L, T)}$-equivariant. 
Note that $\widetilde{N_G(L, T)}$-action on $Z_{\mathfrak{l}, X'}$ is not effective. 
Since $N_L(T)$-action on $O'$ lifts to an action on $(X')^0$, $N_L(T)$ is regarded as 
a subgroup of $\widetilde{N_G(L,T)}$ and it is characterized as the subgroup consisting of the elements 
which act trivially on $Z_{\mathfrak{l}, X'}$. Define 
$$\widetilde{W'} :=  \widetilde{N_G(L, T)}/N_L(T).$$ Then we have an exact sequence 
$$ 1 \to \mathrm{Aut}((X')^0/O') \to \widetilde{W'} \to W' \to 1$$ and 
$\widetilde{W}'$ acts on $Z_{\mathfrak{l}, X'}$. Let us consider the 
central fiber $Z_{\mathfrak{l}, X', 0}$ of the map $Z_{\mathfrak{l}, X'} \to \mathfrak{k}$. 
Then $Z_{\mathfrak{l}, X', 0} = X$ and $\widetilde{W'}$ acts on $X$ as a covering transformation of $X \to \bar{O}$. Therefore we have a homomorphism $$\rho_X: \widetilde{W'} \to \mathrm{Aut}(X/\bar{O}).$$ 
We define $W_X$ as the kernel of $\rho_X$: $$W_X := \mathrm{Ker}(\rho_X).$$ 
Let us consider the composite $\mathrm{Aut}((X')^0/O') \to \widetilde{W'} \stackrel{\rho_X}\to 
\mathrm{Aut}(X/\bar{O})$. This is an injection. In fact, an element $\tau \in \mathrm{Aut}((X')^0/O')$ determines an automorphism $$\phi_{\tau}: G \times^Q (\mathfrak{n} \times X') \to 
G \times^Q (\mathfrak{n} \times X'), \:\: [g, (x, \tilde{y})] \mapsto [g, (x, \tau{\tilde y})],$$ which induces  
an automorphism of $X$ over $\bar{O}$.  If $\rho_X (\tau) = id$, then we get 
a commutative diagram 
\begin{equation} 
\begin{CD} 
G \times^Q (\mathfrak{n} \times X') @>{\phi_{\tau}}>> G \times^Q(\mathfrak{n} \times X') \\ 
@VVV @VVV \\ 
X @>{id_X}>> X     
\end{CD} 
\end{equation}     
This means that $\phi_{\tau} = id$ and, hence $\tau = id$. 
We regard $\mathrm{Aut}((X')^0/O')$ as a subgroup of 
$\mathrm{Aut}(X/\bar{O})$ in this way. 
Then we have a commutative diagram of exact sequences 
\begin{equation} 
\begin{CD}
@. @. 1 @. 1 @. @. \\ 
@. @. @VVV  @VVV @. @. \\
@. 1 @>>> W_X @>>>  \pi (W_X) @>>> 1 \\  
@. @VVV  @VVV  @VVV @. \\ 
1 @>>> \mathrm{Aut}((X')^0/O') @>>> \widetilde{W'} @>{\pi}>> W' @>>> 1 \\ 
@..  @V{id}VV  @V{\rho_X}VV   @VVV  @.. \\ 
1 @>>>  \mathrm{Aut}((X')^0/O') @>>> \mathrm{Im}(\rho_X) @>>> 
\mathrm{Im}(\rho_X)/\mathrm{Aut}((X')^0/O') @>>> 1  \\
\end{CD} 
\end{equation}\vspace{0.2cm}     
Therefore, $W_X$ can be regarded as a subgroup of $W'$.  \vspace{0.2cm}

{\em Twisting} (cf. \cite{Na 1}, pp.386-388) 
Define 
$$\mathcal{S}(L) := \{Q \subset G \:\: \vert \:\: Q \: \mathrm{is \: a \: parabolic \: subgroup \: which  \: contains \: L \: as \: Levi \: part}\}.$$  
Take an element $Q$ from $\mathcal{S}(L)$. Then we obtain a new element $Q' \in 
\mathcal{S}(L)$ from $Q$ by making an operation called a {\em twisting}.  We shall explain this. 
Fix a maximal torus $T$ of $L$. Put $\mathfrak{t} = Lie(T)$ and $\mathfrak{l} = 
Lie (L)$. Let $\Phi \subset \mathfrak{t}^*$ be the root system determined by $(\mathfrak{g}, \mathfrak{t})$. If we take a basis $\Delta$ of $\Phi$, then we can make a Dynkin diagram $D$ attaching an element (= simple root)  of $\Delta$ to each vertex.  Moreover, if we mark 
some vertices of $D$ with back, then it determines a parabolic subalgebra $\mathfrak{q}$ 
of $\mathfrak{g}$. More precisely, let $I \subset \Delta$ be the set of simple roots corresponding to unmarked vertices. Let $\Phi_I$ be the root subsystem of $\Phi$ generated 
by $I$, and put $\Phi_I^- := \Phi_I \cap \Phi^-$, where $\Phi^-$ is the set of negative roots. 
Then  
$$\mathfrak{q} = \mathfrak{t} \oplus \bigoplus_{\alpha \in \Phi_I^-} \mathfrak{g}_{\alpha} 
\oplus \bigoplus_{\alpha \in \Phi^+}\mathfrak{g}_{\alpha} $$ 
is a parabolic subalgebra of $\mathfrak{g}$. 
For $\beta \in \Delta - I$, we make a new marked Dynkin diagram $\bar{D}$ from 
$D$ by erasing the marking of $\beta$. By definition, the set $\bar{I}$ of unmarked 
simple roots for $\bar{D}$ is $I \cup \{\beta\}$. The new marked Dynking diagram $\bar{D}$ 
determines a parabolic subalgebra $\bar{q}$. 
The unmarked vertices of $\bar{D}$ determines 
a Dynkin subdiagram, which is decomposed into the disjoint sum of the connected component 
containing $\beta$ and the union of other components. Correspondingly, we have a decomposition $\bar{I} = I_{\beta} \cup I'_{\beta}$ with $\beta \in I_{\beta}$. 
Let $\bar{\mathfrak{l}}$ be the Levi part of $\bar{q}$ containing $\mathfrak{t}$, and 
let $\mathfrak{z}(\bar{\mathfrak l})$ be is its center.  Then $\bar{\mathfrak{l}}/\mathfrak{z}(\bar{\mathfrak l})$ is decomposed into the direct sum of simple factors. 
Let $\bar{\mathfrak{l}}_{\beta}$ be the simple factor corresponding to $I_{\beta}$, and let 
$\bar{\mathfrak{l}}'_{\beta}$ be the direct sum of other simple factors. 
Then $$\bar{\mathfrak{l}}/\mathfrak{z}(\bar{\mathfrak l}) = \bar{\mathfrak{l}}_{\beta} \oplus 
\bar{\mathfrak{l}}'_{\beta}.$$ 
Put $\mathfrak{q}_{\beta} := \mathfrak{q} \cap \bar{\mathfrak{l}}_{\beta}$. Then 
$\mathfrak{q}_{\beta}$ is a parabolic subalgebra of $\bar{\mathfrak{l}}_{\beta}$.   
Let $\pi: \bar{\mathfrak{l}} \to \bar{\mathfrak{l}}/\mathfrak{z}(\bar{\mathfrak l})$ be the 
quotient map. Note that 
$$\bar{\mathfrak{q}} = \bar{\mathfrak l} \oplus \mathfrak{n}(\bar{\mathfrak q}), $$
$$\mathfrak{q} = \pi^{-1}(\mathfrak{q}_{\beta} \oplus \bar{\mathfrak{l}}'_{\beta}) \oplus 
\mathfrak{n}(\bar{\mathfrak q}). $$
Let $\Phi_{I_{\beta}}$ be the root subsystem of $\Phi$ generated by 
$I_{\beta}$. We put $\mathfrak{t}_{\beta} : = \mathfrak{t}/\mathfrak{z}(\bar{\mathfrak l}) 
\cap \bar{\mathfrak{l}}_{\beta}$. Then $\Phi_{I_{\beta}}$ is the root system for 
$(\bar{\mathfrak{l}}_{\beta}, \mathfrak{t}_{\beta})$. 
The root system $\Phi_{I_{\beta}}$ has a standard involution $-1$. 
There is a unique involution $\phi_{\beta}$ of $\bar{\mathfrak{l}}_{\beta}$ which stabilizes 
$\mathfrak{t}_{\beta}$ and which acts on $\Phi_{I_{\beta}}$ as $-1$. 
We define 
$$\mathfrak{q}' := \pi^{-1}(\phi_{\beta}(\mathfrak{q}_{\beta}) \oplus \bar{\mathfrak{l}}'_{\beta}) \oplus \mathfrak{n}(\bar{\mathfrak q}).$$ 
Then $\mathfrak{q}'$ is a parabolic subalgebra of $\mathfrak{g}$, which contains $\mathfrak{l}$ as the Levi part.  We say that $\mathfrak{q}'$ is a {\em twist} of $\mathfrak{q}$ by $\beta$. 
\vspace{0.2cm}

\begin{Example}\label{Example 6} (i) Assume that $\mathfrak{q}$ is determined by the following marked Dynkin 
diagram $D$ \end{Example} 
  
\begin{picture}(400,20) 
\put(0,-3){- - }\put(15,0){\line(1,0){10}}\put(28,0){\circle*{5}}\put(25,5){$\alpha_0$}\put(33,0){\line(1,0){22}}\put(60,0){\circle{5}}
\put(55,5){$\alpha_1$} 
\put(65,0){\line(1,0){20}}\put(90,0){\circle{5}}\put(85,5){$\alpha_2$}
\put(95,0){\line(1,0){10}}\put(110,-3){- - -}\put(135,0){\line(1,0){10}}
\put(150,0){\circle*{5}}\put(145,5){$\alpha_k$}\put(150,0){\line(1,0){10}}
\put(165,-3){- - -}\put(190,0){\line(1,0){10}}\put(205,0){\circle{5}}
\put(195,5){$\alpha_{n-1}$}\put(210,0){\line(1,0){20}}\put(235,0){\circle{5}}
\put(230,5){$\alpha_n$}\put(240,0){\line(1,0){20}}\put(265,0){\circle*{5}} 
\put(260,5){$\alpha_{n+1}$}\put(268,0){\line(1,0){10}} 
\put(283,-3){- - -}    
\end{picture}  \vspace{1.0cm} 

We make a new marked Dynkin diagram $\bar{D}$ by erasing the marking of $\alpha_k$. 
The unmarked vertices of $\bar{D}$ determine a non-connected Dynkin subdiagram. 
Its connected component containing $\alpha_k$ is of type $A_n$. 
The twist $\mathfrak{q}'$ of $\mathfrak{q}$ by $\alpha_k$ is then determined by 
the marked Dynkin diagram $D'$:
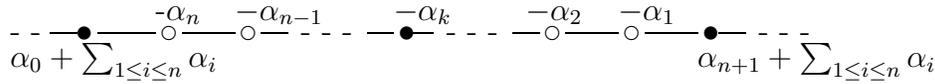
\vspace{0.5cm} 

\begin{picture}(400,20) 
\put(0,-3){- - }\put(15,0){\line(1,0){10}}\put(28,0){\circle*{5}}\put(0,-13){$\alpha_0 + \sum_{1 \le i \le n}\alpha_i$}\put(33,0){\line(1,0){22}}\put(60,0){\circle{5}}
\put(55,5){-$\alpha_n$} 
\put(65,0){\line(1,0){20}}\put(90,0){\circle{5}}\put(85,5){$-\alpha_{n-1}$}
\put(95,0){\line(1,0){10}}\put(110,-3){- - -}\put(135,0){\line(1,0){10}}
\put(150,0){\circle*{5}}\put(145,5){$-\alpha_k$}\put(150,0){\line(1,0){10}}
\put(165,-3){- - -}\put(190,0){\line(1,0){10}}\put(205,0){\circle{5}}
\put(195,5){$-\alpha_2$}\put(210,0){\line(1,0){20}}\put(235,0){\circle{5}}
\put(230,5){$-\alpha_1$}\put(240,0){\line(1,0){20}}\put(265,0){\circle*{5}} 
\put(260,-13){$\alpha_{n+1} + \sum_{1 \le i \le n}\alpha_i$}\put(268,0){\line(1,0){10}} 
\put(283,-3){- - -}    
\end{picture}  \vspace{1.0cm} 

Other simple roots in $D$ different from $\{\alpha_i\}_{0 \le i \le n+1}$ 
do not have any effect by the twist. \vspace{0.2cm} 

(ii) Assume that $\mathfrak{q}$ is determined by the marked Dynkin diagram 
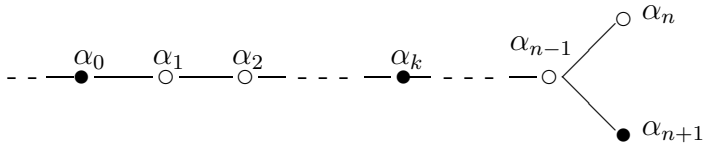
\vspace{0.5cm}
 
\begin{picture}(400,20) 
\put(0,-3){- - }\put(15,0){\line(1,0){10}}\put(28,0){\circle*{5}}\put(25,5){$\alpha_0$}\put(33,0){\line(1,0){22}}\put(60,0){\circle{5}}
\put(55,5){$\alpha_1$} 
\put(65,0){\line(1,0){20}}\put(90,0){\circle{5}}\put(85,5){$\alpha_2$}
\put(95,0){\line(1,0){10}}\put(110,-3){- - -}\put(135,0){\line(1,0){10}}
\put(150,0){\circle*{5}}\put(145,5){$\alpha_k$}\put(150,0){\line(1,0){10}}
\put(165,-3){- - -}\put(190,0){\line(1,0){10}}\put(205,0){\circle{5}}
\put(190,10){$\alpha_{n-1}$}
\put(210,0){\line(1,1){20}}
\put(210,0){\line(1,-1){20}}
\put(233,22){\circle{5}}\put(240,22){$\alpha_n$}
\put(233,-22){\circle*{5}}\put(240,-22){$\alpha_{n+1}$}
\end{picture}  \vspace{1.0cm}

The twist $\mathfrak{q}'$ of $\mathfrak{q}$ by $\alpha_k$ is then determined by 
the marked Dynkin diagram 
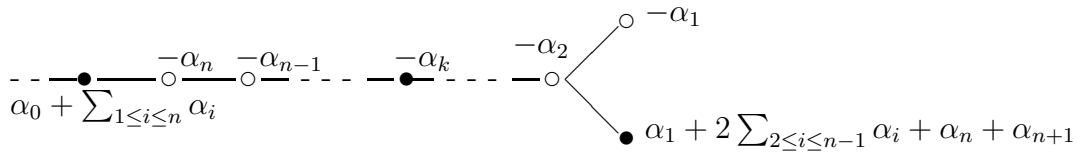
\vspace{0.5cm}

\begin{picture}(400,20) 
\put(0,-3){- - }\put(15,0){\line(1,0){10}}\put(28,0){\circle*{5}}\put(0,-13){$\alpha_0 + \sum_{1 \le i \le n}\alpha_i$}\put(33,0){\line(1,0){22}}\put(60,0){\circle{5}}
\put(55,5){$-\alpha_n$} 
\put(65,0){\line(1,0){20}}\put(90,0){\circle{5}}\put(85,5){$-\alpha_{n-1}$}
\put(95,0){\line(1,0){10}}\put(110,-3){- - -}\put(135,0){\line(1,0){10}}
\put(150,0){\circle*{5}}\put(145,5){$-\alpha_k$}\put(150,0){\line(1,0){10}}
\put(165,-3){- - -}\put(190,0){\line(1,0){10}}\put(205,0){\circle{5}}
\put(190,10){$-\alpha_2$}
\put(210,0){\line(1,1){20}}
\put(210,0){\line(1,-1){20}}
\put(233,22){\circle{5}}\put(240,22){$-\alpha_1$}
\put(233,-22){\circle*{5}}\put(240,-22){$\alpha_1 + 2\sum_{2 \le i \le n-1}\alpha_i + \alpha_n + \alpha_{n+1}$}
\end{picture}  \vspace{1.5cm}

(iii) Assume that $\mathfrak{q}$ is determined by the marked Dynkin diagram \vspace{0.5cm}

\begin{picture}(400,20) 
\put(0,-3){- - }\put(15,0){\line(1,0){10}}\put(28,0){\circle*{5}}\put(25,5){$\alpha_0$}\put(33,0){\line(1,0){22}}\put(60,0){\circle{5}}
\put(55,5){$\alpha_1$} 
\put(65,0){\line(1,0){20}}\put(90,0){\circle{5}}\put(85,5){$\alpha_2$}
\put(95,0){\line(1,0){10}}\put(110,-3){- - -}\put(135,0){\line(1,0){10}}
\put(150,0){\circle*{5}}\put(145,5){$\alpha_k$}\put(150,0){\line(1,0){10}}
\put(165,-3){- - -}\put(190,0){\line(1,0){10}}\put(205,0){\circle{5}}
\put(195,5){$\alpha_{n-1}$}\put(210,0){\line(1,0){20}}\put(235,0){\circle{5}}
\put(230,5){$\alpha_n$}\put(240,-3){$\Longrightarrow$}\put(265,0){\circle*{5}} 
\put(260,5){$\alpha_{n+1}$}    
\end{picture}  \vspace{1.0cm}

The twist $\mathfrak{q}'$ of $\mathfrak{q}$ by $\alpha_k$ is then determined by 
the marked Dynkin diagram \vspace{0.5cm}

\begin{picture}(400,20) 
\put(0,-3){- - }\put(15,0){\line(1,0){10}}\put(28,0){\circle*{5}}\put(0,-13){$\alpha_0 + \sum_{1 \le i \le n}\alpha_i$}\put(33,0){\line(1,0){22}}\put(60,0){\circle{5}}
\put(55,5){$-\alpha_n$} 
\put(65,0){\line(1,0){20}}\put(90,0){\circle{5}}\put(85,5){$-\alpha_{n-1}$}
\put(95,0){\line(1,0){10}}\put(110,-3){- - -}\put(135,0){\line(1,0){10}}
\put(150,0){\circle*{5}}\put(145,5){$-\alpha_k$}\put(150,0){\line(1,0){10}}
\put(165,-3){- - -}\put(190,0){\line(1,0){10}}\put(205,0){\circle{5}}
\put(195,5){$-\alpha_2$}\put(210,0){\line(1,0){20}}\put(235,0){\circle{5}}
\put(230,5){$-\alpha_1$}\put(240,-3){$\Longrightarrow$}\put(265,0){\circle*{5}} 
\put(260,-13){$\sum_{1 \le i \le n}\alpha_i + \alpha_{n+1}$}    
\end{picture}  \vspace{1.5cm}

(iv) Assume that $\mathfrak{q}$ is determined by the marked Dynkin diagram \vspace{0.5cm}

\begin{picture}(400,20) 
\put(0,-3){- - }\put(15,0){\line(1,0){10}}\put(28,0){\circle*{5}}\put(25,5){$\alpha_0$}\put(33,0){\line(1,0){22}}\put(60,0){\circle{5}}
\put(55,5){$\alpha_1$} 
\put(65,0){\line(1,0){20}}\put(90,0){\circle{5}}\put(85,5){$\alpha_2$}
\put(95,0){\line(1,0){10}}\put(110,-3){- - -}\put(135,0){\line(1,0){10}}
\put(150,0){\circle*{5}}\put(145,5){$\alpha_k$}\put(150,0){\line(1,0){10}}
\put(165,-3){- - -}\put(190,0){\line(1,0){10}}\put(205,0){\circle{5}}
\put(195,5){$\alpha_{n-1}$}\put(210,0){\line(1,0){20}}\put(235,0){\circle{5}}
\put(230,5){$\alpha_n$}\put(240,-3){$\Longleftarrow$}\put(265,0){\circle*{5}} 
\put(260,5){$\alpha_{n+1}$}    
\end{picture}  \vspace{1.0cm}

The twist $\mathfrak{q}'$ of $\mathfrak{q}$ by $\alpha_k$ is then determined by 
the marked Dynkin diagram \vspace{0.5cm}

\begin{picture}(400,20) 
\put(0,-3){- - }\put(15,0){\line(1,0){10}}\put(28,0){\circle*{5}}\put(0,-13){$\alpha_0 + \sum_{1 \le i \le n}\alpha_i$}\put(33,0){\line(1,0){22}}\put(60,0){\circle{5}}
\put(55,5){$-\alpha_n$} 
\put(65,0){\line(1,0){20}}\put(90,0){\circle{5}}\put(85,5){$-\alpha_{n-1}$}
\put(95,0){\line(1,0){10}}\put(110,-3){- - -}\put(135,0){\line(1,0){10}}
\put(150,0){\circle*{5}}\put(145,5){$-\alpha_k$}\put(150,0){\line(1,0){10}}
\put(165,-3){- - -}\put(190,0){\line(1,0){10}}\put(205,0){\circle{5}}
\put(195,5){$-\alpha_2$}\put(210,0){\line(1,0){20}}\put(235,0){\circle{5}}
\put(230,5){$-\alpha_1$}\put(240,-3){$\Longleftarrow$}\put(265,0){\circle*{5}} 
\put(260,-13){$2\sum_{1 \le i \le n}\alpha_i + \alpha_{n+1}$}    
\end{picture}  \vspace{1.5cm}

(v) Assume that $\mathfrak{q}$ is determined by the marked Dynkin diagram \vspace{0.5cm}

\begin{picture}(400,20) 
\put(0,-3){- - }\put(15,0){\line(1,0){10}}\put(28,0){\circle*{5}}\put(25,5){$\alpha_0$}\put(33,0){\line(1,0){22}}\put(60,0){\circle{5}}
\put(55,5){$\alpha_1$} 
\put(65,0){\line(1,0){20}}\put(90,0){\circle{5}}\put(85,5){$\alpha_2$}
\put(95,0){\line(1,0){10}}\put(110,-3){- - -}\put(135,0){\line(1,0){10}}
\put(150,0){\circle*{5}}\put(145,5){$\alpha_k$}\put(150,0){\line(1,0){10}}
\put(165,-3){- - -}\put(190,0){\line(1,0){10}}\put(205,0){\circle{5}}
\put(195,5){$\alpha_{n-1}$}\put(210,-3){$\Longrightarrow$}\put(235,0){\circle{5}}
\put(230,5){$\alpha_n$}    
\end{picture}  \vspace{1.0cm}

The twist $\mathfrak{q}'$ of $\mathfrak{q}$ by $\alpha_k$ is then determined by 
the marked Dynkin diagram \vspace{0.5cm}

\begin{picture}(400,20) 
\put(0,-3){- - }\put(15,0){\line(1,0){10}}\put(28,0){\circle*{5}}\put(0,-13){$\alpha_0 + 2 \sum_{1 \le i \le n}\alpha_i$}\put(33,0){\line(1,0){22}}\put(60,0){\circle{5}}
\put(55,5){$-\alpha_1$} 
\put(65,0){\line(1,0){20}}\put(90,0){\circle{5}}\put(85,5){$-\alpha_2$}
\put(95,0){\line(1,0){10}}\put(110,-3){- - -}\put(135,0){\line(1,0){10}}
\put(150,0){\circle*{5}}\put(145,5){$-\alpha_k$}\put(150,0){\line(1,0){10}}
\put(165,-3){- - -}\put(190,0){\line(1,0){10}}\put(205,0){\circle{5}}
\put(195,5){$-\alpha_{n-1}$}\put(210,-3){$\Longrightarrow$}\put(235,0){\circle{5}}
\put(230,5){$-\alpha_n$}    
\end{picture}  \vspace{1.5cm}

(vi) Assume that $\mathfrak{q}$ is determined by the marked Dynkin diagram \vspace{0.5cm}

\begin{picture}(400,20) 
\put(0,-3){- - }\put(15,0){\line(1,0){10}}\put(28,0){\circle*{5}}\put(25,5){$\alpha_0$}\put(33,0){\line(1,0){22}}\put(60,0){\circle{5}}
\put(55,5){$\alpha_1$} 
\put(65,0){\line(1,0){20}}\put(90,0){\circle{5}}\put(85,5){$\alpha_2$}
\put(95,0){\line(1,0){10}}\put(110,-3){- - -}\put(135,0){\line(1,0){10}}
\put(150,0){\circle*{5}}\put(145,5){$\alpha_k$}\put(150,0){\line(1,0){10}}
\put(165,-3){- - -}\put(190,0){\line(1,0){10}}\put(205,0){\circle{5}}
\put(195,5){$\alpha_{n-1}$}\put(210,-3){$\Longleftarrow$}\put(235,0){\circle{5}}
\put(230,5){$\alpha_n$}    
\end{picture}  \vspace{1.0cm}

The twist $\mathfrak{q}'$ of $\mathfrak{q}$ by $\alpha_k$ is then determined by 
the marked Dynkin diagram \vspace{0.5cm}

\begin{picture}(400,20) 
\put(0,-3){- - }\put(15,0){\line(1,0){10}}\put(28,0){\circle*{5}}\put(0,-13){$\alpha_0 + 2 \sum_{1 \le i \le n -1}\alpha_i + \alpha_n$}\put(33,0){\line(1,0){22}}\put(60,0){\circle{5}}
\put(55,5){$-\alpha_1$} 
\put(65,0){\line(1,0){20}}\put(90,0){\circle{5}}\put(85,5){$-\alpha_2$}
\put(95,0){\line(1,0){10}}\put(110,-3){- - -}\put(135,0){\line(1,0){10}}
\put(150,0){\circle*{5}}\put(145,5){$-\alpha_k$}\put(150,0){\line(1,0){10}}
\put(165,-3){- - -}\put(190,0){\line(1,0){10}}\put(205,0){\circle{5}}
\put(195,5){$-\alpha_{n-1}$}\put(210,-3){$\Longleftarrow$}\put(235,0){\circle{5}}
\put(230,5){$-\alpha_n$}    
\end{picture}  \vspace{1.5cm}

(vii) Let $n$ be an odd number and let $k$ be an integer such that $1 \le k \le n-2$. 
Assume that $\mathfrak{q}$ is determined by the marked Dynkin diagram 
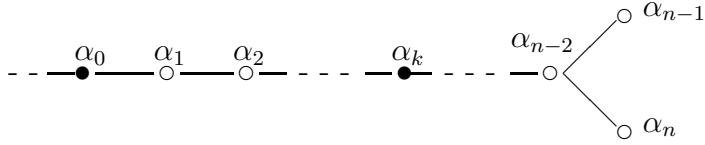
\vspace{0.5cm}

\begin{picture}(400,20) 
\put(0,-3){- - }\put(15,0){\line(1,0){10}}\put(28,0){\circle*{5}}\put(25,5){$\alpha_0$}\put(33,0){\line(1,0){22}}\put(60,0){\circle{5}}
\put(55,5){$\alpha_1$} 
\put(65,0){\line(1,0){20}}\put(90,0){\circle{5}}\put(85,5){$\alpha_2$}
\put(95,0){\line(1,0){10}}\put(110,-3){- - -}\put(135,0){\line(1,0){10}}
\put(150,0){\circle*{5}}\put(145,5){$\alpha_k$}\put(150,0){\line(1,0){10}}
\put(165,-3){- - -}\put(190,0){\line(1,0){10}}\put(205,0){\circle{5}}
\put(190,10){$\alpha_{n-2}$}
\put(210,0){\line(1,1){20}}
\put(210,0){\line(1,-1){20}}
\put(233,22){\circle{5}}\put(240,22){$\alpha_{n-1}$}
\put(233,-22){\circle{5}}\put(240,-22){$\alpha_{n}$}
\end{picture}  \vspace{1.5cm}
  
Then the twist $\mathfrak{q}'$ of $\mathfrak{q}$ by $\alpha_k$ is determined by 
the marked Dynkin diagram 
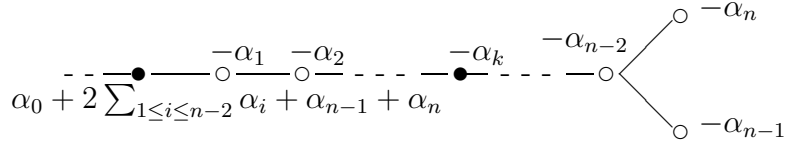
\vspace{0.5cm}

\begin{picture}(400,20) 
\put(0,-3){- - }\put(15,0){\line(1,0){10}}\put(28,0){\circle*{5}}\put(-20,-13){$\alpha_0 + 2\sum_{1 \le i \le n-2}\alpha_i + \alpha_{n-1} + \alpha_n$}\put(33,0){\line(1,0){22}}\put(60,0){\circle{5}}
\put(55,5){$-\alpha_1$} 
\put(65,0){\line(1,0){20}}\put(90,0){\circle{5}}\put(85,5){$-\alpha_2$}
\put(95,0){\line(1,0){10}}\put(110,-3){- - -}\put(135,0){\line(1,0){10}}
\put(150,0){\circle*{5}}\put(145,5){$-\alpha_k$}\put(150,0){\line(1,0){10}}
\put(165,-3){- - -}\put(190,0){\line(1,0){10}}\put(205,0){\circle{5}}
\put(180,10){$-\alpha_{n-2}$}
\put(210,0){\line(1,1){20}}
\put(210,0){\line(1,-1){20}}
\put(233,22){\circle{5}}\put(240,22){$-\alpha_n$}
\put(233,-22){\circle{5}}\put(240,-22){$-\alpha_{n-1}$}
\end{picture}  \vspace{1.5cm}

On the other hand, when $k = n-1$, $\mathfrak{q}$ is given by the marked Dynkin diagram 
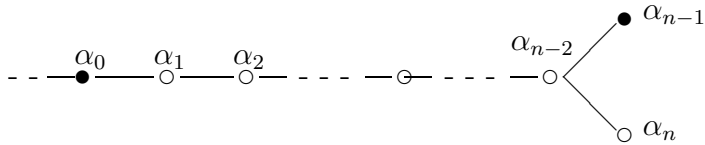
\vspace{0.5cm} 

\begin{picture}(400,20) 
\put(0,-3){- - }\put(15,0){\line(1,0){10}}\put(28,0){\circle*{5}}\put(25,5){$\alpha_0$}\put(33,0){\line(1,0){22}}\put(60,0){\circle{5}}
\put(55,5){$\alpha_1$} 
\put(65,0){\line(1,0){20}}\put(90,0){\circle{5}}\put(85,5){$\alpha_2$}
\put(95,0){\line(1,0){10}}\put(110,-3){- - -}\put(135,0){\line(1,0){10}}
\put(150,0){\circle{5}}\put(150,0){\line(1,0){10}}
\put(165,-3){- - -}\put(190,0){\line(1,0){10}}\put(205,0){\circle{5}}
\put(190,10){$\alpha_{n-2}$}
\put(210,0){\line(1,1){20}}
\put(210,0){\line(1,-1){20}}
\put(233,22){\circle*{5}}\put(240,22){$\alpha_{n-1}$}
\put(233,-22){\circle{5}}\put(240,-22){$\alpha_{n}$}
\end{picture}  \vspace{1.5cm}

Then the twist $\mathfrak{q}'$ of $\mathfrak{q}$ by $\alpha_{n-1}$ is determined by 
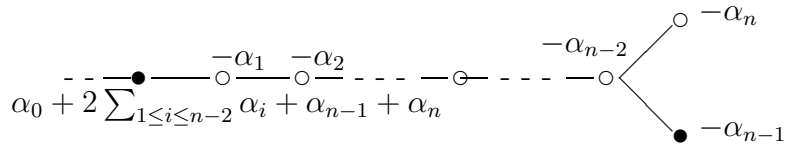
\vspace{0.5cm}

\begin{picture}(400,20) 
\put(0,-3){- - }\put(15,0){\line(1,0){10}}\put(28,0){\circle*{5}}\put(-20,-13){$\alpha_0 + 2\sum_{1 \le i \le n-2}\alpha_i + \alpha_{n-1} + \alpha_n$}\put(33,0){\line(1,0){22}}\put(60,0){\circle{5}}
\put(55,5){$-\alpha_1$} 
\put(65,0){\line(1,0){20}}\put(90,0){\circle{5}}\put(85,5){$-\alpha_2$}
\put(95,0){\line(1,0){10}}\put(110,-3){- - -}\put(135,0){\line(1,0){10}}
\put(150,0){\circle{5}}\put(150,0){\line(1,0){10}}
\put(165,-3){- - -}\put(190,0){\line(1,0){10}}\put(205,0){\circle{5}}
\put(180,10){$-\alpha_{n-2}$}
\put(210,0){\line(1,1){20}}
\put(210,0){\line(1,-1){20}}
\put(233,22){\circle{5}}\put(240,22){$-\alpha_n$}
\put(233,-22){\circle*{5}}\put(240,-22){$-\alpha_{n-1}$}
\end{picture}  \vspace{1.5cm}

(viii) Let $n$ be an even integer, and let $k$ be an integer such that  $1 \le k \le n-1$. 
Assume $\mathfrak{q}$ is determined by the marked Dynkin diagram 
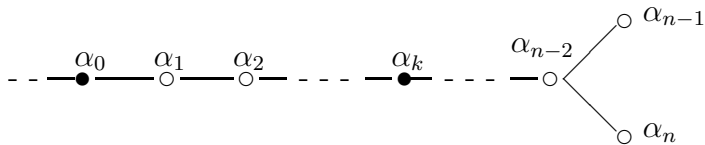
\vspace{0.5cm}

\begin{picture}(400,20) 
\put(0,-3){- - }\put(15,0){\line(1,0){10}}\put(28,0){\circle*{5}}\put(25,5){$\alpha_0$}\put(33,0){\line(1,0){22}}\put(60,0){\circle{5}}
\put(55,5){$\alpha_1$} 
\put(65,0){\line(1,0){20}}\put(90,0){\circle{5}}\put(85,5){$\alpha_2$}
\put(95,0){\line(1,0){10}}\put(110,-3){- - -}\put(135,0){\line(1,0){10}}
\put(150,0){\circle*{5}}\put(145,5){$\alpha_k$}\put(150,0){\line(1,0){10}}
\put(165,-3){- - -}\put(190,0){\line(1,0){10}}\put(205,0){\circle{5}}
\put(190,10){$\alpha_{n-2}$}
\put(210,0){\line(1,1){20}}
\put(210,0){\line(1,-1){20}}
\put(233,22){\circle{5}}\put(240,22){$\alpha_{n-1}$}
\put(233,-22){\circle{5}}\put(240,-22){$\alpha_{n}$}
\end{picture}  \vspace{1.5cm}

The twist $\mathfrak{q}'$ of $\mathfrak{q}$ by $\alpha_k$ is then determined by 
the marked Dynkin diagram 
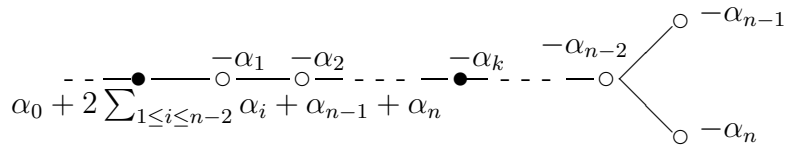
\vspace{0.5cm}

\begin{picture}(400,20) 
\put(0,-3){- - }\put(15,0){\line(1,0){10}}\put(28,0){\circle*{5}}\put(-20,-13){$\alpha_0 + 2\sum_{1 \le i \le n-2}\alpha_i + \alpha_{n-1} + \alpha_n$}\put(33,0){\line(1,0){22}}\put(60,0){\circle{5}}
\put(55,5){$-\alpha_1$} 
\put(65,0){\line(1,0){20}}\put(90,0){\circle{5}}\put(85,5){$-\alpha_2$}
\put(95,0){\line(1,0){10}}\put(110,-3){- - -}\put(135,0){\line(1,0){10}}
\put(150,0){\circle*{5}}\put(145,5){$-\alpha_k$}\put(150,0){\line(1,0){10}}
\put(165,-3){- - -}\put(190,0){\line(1,0){10}}\put(205,0){\circle{5}}
\put(180,10){$-\alpha_{n-2}$}
\put(210,0){\line(1,1){20}}
\put(210,0){\line(1,-1){20}}
\put(233,22){\circle{5}}\put(240,22){$-\alpha_{n-1}$}
\put(233,-22){\circle{5}}\put(240,-22){$-\alpha_n$}
\end{picture}  \vspace{1.5cm}

\vspace{0.2cm}

\begin{Thm}\label{Theorem 7} (1) Every element $Q \in \mathcal{S}(L)$ is obtained from 
a fixed element $Q_0 \in \mathcal{S}(L)$ by a finite succession of the twists.   
Every {\bf Q}-factorial terminalization of $Z_{\mathfrak{l}, X'}$ has the form $G \times^Q (\mathfrak{r} \times X')$ with $Q \in \mathcal{S}(L)$. Conversely, for any $Q \in \mathcal{S}(L)$,  $G \times^Q (\mathfrak{r} \times X')$ gives a {\bf Q}-factorial terminalization 
of $Z_{\mathfrak{l}, X'}$. Moreover, there is  a one-one correspondence 
between the set of all {\bf Q}-factorial terminalizations of $Z_{\mathfrak{l}, X'}$ and 
$\mathcal{S}(L)$. 

(2) Every {\bf Q}-factorial terminalization of $X$ has the form 
$G \times^Q (\mathfrak{n} \times X')$ with $Q \in \mathcal{S}(L)$. Conversely, for any $Q \in \mathcal{S}(L)$,  $G \times^Q (\mathfrak{n} \times X')$ gives a {\bf Q}-factorial terminalization 
of $X$.  \end{Thm}

\begin{Rem}\label{Remark 8} \end{Rem}

(i) In (2) different two parabolic subgroups $Q, Q' \in \mathcal{S}(L)$ 
may possibly give the same {\bf Q}-factorial terminalization of $X$. 

(ii) The pair $(X', O')$ is uniquely determined by $X$. 
We call $(X', O')$ the {\em core} of $X$.  
\vspace{0.2cm}

{\em Proof}. (1) The first statement is proved in \cite{Na 1}, Theorem 1.3. For each $Q \in \mathcal{S}(L)$, we have isomorphisms 
$$\mathrm{Pic}(G \times^Q (\mathfrak{r} \times X') \cong \mathrm{Pic}(G/Q) \cong \mathrm{Hom}_{alg.gp}(L, \mathbf{C}^*).$$ The nef cone $\overline{\mathrm{Amp}}(G \times^Q (\mathfrak{r} \times X')) \subset \mathrm{Pic}(G \times^Q(\mathfrak{r} \times X')) \otimes_{\mathbf Z}
\mathbf{R}$ is identified with a cone in $\mathrm{Hom}_{alg.gp}(L, \mathbf{C}^*)\otimes_{\mathbf Z}\mathbf{R}$. By \cite{Na 1}, Remark 1.6  we have 
$$\mathrm{Hom}_{alg.gp}(L, \mathbf{C}^*)\otimes_{\mathbf Z}\mathbf{R} = \bigcup_{Q \in 
\mathcal{S}(L)} \overline{\mathrm{Amp}}(G \times^Q
(\mathfrak{r} \times X')) $$ Fix $Q_0 \in \mathcal{S}(L)$ and identify $\mathrm{Pic}(G \times^{Q_0}(\mathfrak{r}_0 \times X')$ with $\mathrm{Hom}_{alg.gp}(L, \mathbf{C}^*)$. Then,   
by Lemma \ref{Lemma 3}, we have $$\mathrm{Hom}_{alg.gp}(L, \mathbf{C}^*)\otimes_{\mathbf Z}\mathbf{R} = \mathrm{Mov}(G \times^{Q_0} (\mathfrak{r}_0 \times X'))$$  
For each $Q \in \mathcal{S}(L)$, we can identify the nef cone  
$\overline{\mathrm{Amp}}(G \times^Q (\mathfrak{r} \times X'))$ as a cone inside 
$\mathrm{Mov}(G \times^{Q_0} (\mathfrak{r}_0 \times X'))$ by the birational map 
$$G \times^Q (\mathfrak{r} \times X') - - \to G \times^{Q_0} (\mathfrak{r}_0 \times X').$$ 
As in \cite{Na 1}, Observation 2 we can check that this identification is the same one explained 
above.  

(2)  Once (1) is proved, (2) is proved in the same way as in \cite{Na 1}, Corollary 3.4. 
$\square$ 
\vspace{0.2cm}

\begin{Thm}\label{Theorem 9} The following commutative diagram 
gives the universal Poisson deformations of $G \times^Q (\mathfrak{n} \times X')$ and $X$:
\begin{equation} 
\begin{CD} 
G \times^Q(\mathfrak{r} \times X') @>>> Z_{\mathfrak{l}, X'}/W_X \\ 
@VVV @VVV \\ 
\mathfrak{k} @>>> \mathfrak{k}/W_X     
\end{CD} 
\end{equation}  \end{Thm}
\vspace{0.2cm}


{\em Proof}. Put $0_X = \pi^{-1}(0)$ for $\pi: X \to \bar{O}$. Then $G \times^Q (\mathfrak{r} + 0_X) \subset G \times^Q(\mathfrak{r} \times X')$ is a closed Poisson subscheme. 
Moreover, $G \times^Q(\mathfrak{r} + 0_X) \to \mathfrak{k}$ is a Poisson deformation 
of $T^*(G/Q) = G \times^Q(\mathfrak{n} + 0_X)$. Let $$p\kappa_{T^*(G/Q)} : \mathfrak{k} \to H^2(T^*(G/Q), \mathbf{C})$$ 
be the Poisson Kodaira Spencer map for this Poisson deformation at $0 \in \mathfrak{k}$. Let $$p\kappa_{G \times^Q(\mathfrak{n} \times X')}: \mathfrak{k} \to H^2(G \times^Q(\mathfrak{n} \times X'), \mathbf{C})$$ be the Poisson Kodaira Spencer map for $G \times^Q(\mathfrak{r} \times X') 
\to \mathfrak{k}$ at $0 \in \mathfrak{k}$. Then $p\kappa_{T^*(G/Q)}$ factorizes as 
$$p\kappa_{T^*(G/Q)}: \mathfrak{k} \stackrel{p\kappa_{G \times^Q(\mathfrak{n} \times X')}}\to 
H^2(G \times^Q(\mathfrak{n} \times X'), \mathbf{C}) \stackrel{res}\to 
H^2(T^*(G/Q)).$$ By \cite{Na 2}, Proposition 2.7, $p\kappa_{T^*(G/Q)}$ is an isomorphism. 
Then $p\kappa_{G \times^Q(\mathfrak{n} \times X')}$ is an injection. Since $\dim \mathfrak{k} 
= \dim H^2(G \times^Q(\mathfrak{n} \times X'), \mathbf{C})$, we see that   
$p\kappa_{G \times^Q(\mathfrak{n} \times X')}$ is an isomorphism, which implies the versality of 
the Poisson deformation.  But, since   
the Poisson deformation functor $\mathrm{PD}_{G \times^Q(\mathfrak{n} \times X')}$ is 
pro-representable by \cite{Na 6}, Corollary 2.5, we see that $G \times^Q (\mathfrak{r} \times X') \to \mathfrak{k}$ is the universal Poisson deformation of $G \times^Q (\mathfrak{n} \times X')$. Then the Weyl group $W(X)$ of $X$ acts linearly on $\mathfrak{k}$ and there is a universal Poisson deformation $$\mathcal{X} \to \mathfrak{k}/W(X)$$ of $X$ (cf. \cite{Na 4}, Main Theorem). The Poisson deformation 
$Z_{\mathfrak{l}, X'} \to \mathfrak{k}$ is nothing but the pull back of the universal Poisson 
deformation by the quotient map ([ibid])
\begin{equation} 
\begin{CD} 
Z_{\mathfrak{l}, X'} @>>> \mathcal{X} \\ 
@VVV @VVV \\ 
\mathfrak{k} @>>> \mathfrak{k}/W(X)     
\end{CD} 
\end{equation}
\vspace{0.2cm}
Then $W(X)$ acts on $Z_{\mathfrak{l}, X'}$ as $G$-equivariant Poisson automorphisms, and 
$\mathcal{X} = Z_{\mathfrak{l}, X'}/W(X)$.

The next proposition completes the proof of Theorem \ref{Theorem 9}. 

\begin{Prop}\label{Proposition 10} $W(X) = W_X$. \end{Prop}

{\em Proof }. We first show that $W(X) 
\subset W_X$. 
The Poisson structure on $G \times^Q (\mathfrak{r} \times X')$ determines a symplectic 2-form 
$\omega_t$ on the regular part of $G \times^Q ((t + \mathfrak{n}) \times X')$ for each 
$t \in \mathfrak{k}$. Each symplectic 2-form $\omega_t$ determines a 2-nd cohomology class 
$[\omega_t] \in H^2(G \times^Q ((t + \mathfrak{n}) \times X'), \mathbf{C})$ (not only a 2-nd 
cohomology class on the regular part). The restriction maps of 2-nd cohomology are both isomorphisms: 
$$H^2(G \times^Q ((t + \mathfrak{n}) \times X'), \mathbf{C}) \stackrel{\mathrm{res}_t \:\: \cong}\leftarrow 
H^2(G \times^Q (\mathfrak{r} \times X'), \mathbf{C}) \stackrel{\mathrm{res}_0 \:\: \cong}\to 
H^2(G \times^Q (\mathfrak{n}) \times X'), \mathbf{C}).$$ 
We define the period map $p: \mathfrak{k} \to H^2(G \times^Q (\mathfrak{n} \times X'), \mathbf{C})$ 
by $p(t) := \mathrm{res}_0 \circ \mathrm{res}_t^{-1} ([\omega_t])$. Since 
$G \times^Q (\mathfrak{r} \times X') \to \mathfrak{k}$ is the universal Poisson deformation of 
$G \times^Q (\mathfrak{n} \times X')$, the period map $p$ is an isomorphism. We identify $\mathfrak{k}$ with 
$H^2(G \times^Q (\mathfrak{r} + X'), \mathbf{C})$ by $\mathrm{res}_0^{-1}\circ p$.      
The $W(X)$-action on $H^2(G \times^Q (\mathfrak{r} \times X'), \mathbf{C})$ is 
defined over $\mathbf{Q}$. In particular, $W(X)$ acts on $H^2(G \times^Q (\mathfrak{r} \times X'), \mathbf{R})$. By \cite{Na 4}, Main Theorem, $H^2(G \times^Q (\mathfrak{r} \times X'), \mathbf{R})$ is divided into finite 
number of chambers by finite linear spaces $\{H_i\}_{i \in I}$ of codimension 1, and each chamber corresponds to an ample cone of $G \times^{Q'} (\mathfrak{r}' \times X')$ for  $Q' \in \mathcal{S}(L)$. There is a finite subset $J$ of $I$ such that $W(X)$ is generated by the reflections $\sigma_i$ with respect to $H_i$ with $i \in J$. Take $i \in J$ and put 
$H := H_i$ and $\sigma := \sigma_i$. Take a chamber which has (an non-empty open set of ) $H$ as a codimension 1 face.  This chamber is $\mathrm{Amp}(G \times^P(\mathfrak{r}({\mathfrak p}) \times X'))$ for some $P \in \mathcal{S}(L)$. The element $\sigma$ induces 
an automorphism $\sigma_{\mathfrak{k}}: \mathfrak{k} \to \mathfrak{k}$. The fiber product  
$(G \times^P (\mathfrak{r}({\mathfrak p}) \times X')) \times_{\mathfrak k} \mathfrak{k}$ 
fits in the commutative diagram 
\begin{equation} 
\begin{CD} 
(G \times^P (\mathfrak{r}({\mathfrak p}) \times X')) \times_{\mathfrak k} \mathfrak{k} @>{\phi_{\sigma}}>> G \times^P (\mathfrak{r}({\mathfrak p}) \times X') \\ 
@VVV @VVV \\ 
Z_{\mathfrak{l}, X'} @>>> Z_{\mathfrak{l}, X'} \\   
@VVV @VVV \\ 
\mathfrak{k} @>{\sigma_{\mathfrak k}}>> \mathfrak{k}  
\end{CD} 
\end{equation}
\vspace{0.2cm}
Here each square is the Cartesian square. Since 
$(G \times^P (\mathfrak{r}({\mathfrak p}) \times X')) \times_{\mathfrak k} \mathfrak{k}$ is a 
{\bf Q}-factorial terminalization of $Z_{\mathfrak{l}, X'}$, one can write it as 
$G \times^{P'} (\mathfrak{r}(\mathfrak{p}') \times X')$ for some $P' \in \mathcal{S}(L)$. 
\vspace{0.2cm}

\begin{Claim}\label{Claim 1} There is an element $w \in W'$ such that 
$P = Ad_w(P')$. \end{Claim}

{\em Proof}. In the commutative diagram above, let us consider the central fibers. Then we 
have 
\begin{equation} 
\begin{CD} 
G \times^{P'} (\mathfrak{n}(\mathfrak{p}') \times X') @>{\cong}>> G \times^P (\mathfrak{n}({\mathfrak p}) \times X') \\ 
@V{\mathfrak{m}'_{P',0}}VV @V{\mathfrak{m}'_{P,0}}VV \\ 
X @>{id}>> X \\   
@VVV @VVV \\ 
\{0\} @>{\sigma}>> \{0\}  
\end{CD} 
\end{equation}
The isomorphism 
$G \times^{P'} (\mathfrak{n}(\mathfrak{p}') \times X') \cong G \times^P (\mathfrak{n}({\mathfrak p}) \times X')$ is $G$-equivariant. Recall that $\pi^{-1}(0)$ is a point for $\pi: X \to \bar{O}$. 
We put $0_X := \pi^{-1}(0)$. Then the $G$-equivariant isomorphism induces a $G$-equivariant 
isomorphism between ${\mathfrak{m}'_{P',0}}^{-1}(0_X) = G/P'$ and ${\mathfrak{m}'_{P,0}}^{-1}(0_X) = G/P$. If the isomorphism sends $\bar{1} \in G/P'$ to $g_0 \in G/P$, then by 
the $G$-equivariace, it sends $\bar{g}$ to $\overline{gg_0}$ for $g \in G$. 
Similarly it sends $\overline{gp'}$ to $\overline{gp'g_0}$ for $p' \in P'$. 
Since $\bar{g} = \overline{gp'}$ in $G/P'$, we have $\overline{gg_0} = \overline{gp'g_0}$ 
in $G/P$. This means that one can write $gp'g_0 = gg_0p$ with a suitable $p \in P$. 
Then $p' = g_0pg_0^{-1}$. Therefore $P'$ is conjugate to $P$. Since $P$ and $P'$ are 
contained in $\mathcal{S}(L)$,  
one can write $P = Ad_w(P')$ with a suitable element $w$ of $W'$ by \cite{Na 1}, Proposition 2.2.      
Moreover, by the same proposition, this $w$ is uniquely determined. 
$\square$ \vspace{0.2cm}
 
By the definition of $W'$, $w$ acts on $\mathfrak{k}$. We denote by $w_{\mathfrak k}$ the automorphism of $\mathfrak{k}$ induced 
by $w$.  We take a lift $\tilde{w} \in \widetilde{W'}$ of $w$ and define an isomorphism 
$$\phi_{\tilde{w}}:  
G \times^{P'} (\mathfrak{r}(\mathfrak{p}') \times X') \to G \times^P(\mathfrak{r}(\mathfrak{p}) \times X')$$ by $[g, (x, \tilde{y})] 
\mapsto [gw^{-1}, (Ad_{w}(x), \tilde{w}(\tilde{y}))]$ for $g \in G$, $x \in 
\mathfrak{r}(\mathfrak{p}')$ and $\tilde{y} \in X'$. There is a commutative diagram 
\begin{equation} 
\begin{CD} 
G \times^{P'} (\mathfrak{r}(\mathfrak{p}') \times X') @>{\phi_{\tilde w}}>> G \times^P(\mathfrak{r}(\mathfrak{p}) \times X')\\ 
@VVV @VVV \\ 
\mathfrak{k} @>{w_{\mathfrak k}}>>  \mathfrak{k}     
\end{CD} 
\end{equation}

\begin{Claim}\label{Claim 2} $w_{\mathfrak k} = \sigma_{\mathfrak k}$. \end{Claim} 
 
{\em Proof}. By using the period map $p$ for $G \times^P (\mathfrak{r} \times X') \to \mathfrak{k}$, we identify $\mathfrak{k}$ 
with $H^2(G \times^P (\mathfrak{n}(\mathfrak{p}) \times X'), \mathbf{C})$.  
The isomorphism $\phi_{\tilde w}$ induces an isomorphism of the 2-nd cohomology  
$$(\phi_{\tilde{w}})^*: H^2(G \times^P(\mathfrak{r}(\mathfrak{p}) \times X'), \mathbf{C}) \to 
H^2(G \times^{P'} (\mathfrak{r}(\mathfrak{p}') \times X'), \mathbf{C}). $$  
On the other hand, since both $G \times^P (\mathfrak{r}(\mathfrak{p}) \times X')$ and 
$G \times^{P'}(\mathfrak{r}(\mathfrak{p}') \times X')$ give {\bf Q}-factorial 
terminalizations of $Z_{\mathfrak{l}, X'}$, there is a natural identification 
$$(\mathfrak{m}_{P}^{-1} \circ \mathfrak{m}_{P'})_* :  H^2(G \times^{P'} (\mathfrak{r}(\mathfrak{p}') \times X'), \mathbf{C}) \to H^2(G \times^P(\mathfrak{r}(\mathfrak{p}) \times X'), \mathbf{C}).$$
Take $t \in \mathfrak{k}$ general. Then, as seen in the proof of Lemma \ref{Lemma 2}, the birational map $\mathfrak{m}_{P}^{-1} \circ \mathfrak{m}_{P'}$ 
identifies $G \times^P ((w_{\mathfrak k}(t) + \mathfrak{n}(\mathfrak{p})) \times X')$ with  
$G \times^{P'}(w_{\mathfrak{k}}(t) + \mathfrak{n}(\mathfrak{p}')) \times X')$.  
We then have the commutative diagram 
\begin{equation} 
\begin{CD} 
H^2(G \times^P (\mathfrak{n}(\mathfrak{p}) \times X'), \mathbf{C}) @<<< H^2(G \times^P (\mathfrak{n}(\mathfrak{p}) \times X'), \mathbf{C}) \\ 
@A{\cong \:\: \mathrm{res}_0}AA @A{\cong \:\: \mathrm{res}_0}AA \\ 
H^2(G \times^P (\mathfrak{r}(\mathfrak{p}) \times X'), \mathbf{C}) @<{(\mathfrak{m}_{P}^{-1} \circ \mathfrak{m}_{P'})_*  \circ (\phi_{\tilde{w}})^*}<< 
H^2(G \times^P (\mathfrak{r}(\mathfrak{p}) \times X'), \mathbf{C}) \\   
@V{\cong \:\: \mathrm{res}_{w_{\mathfrak{k}}(t)}}VV @V{\cong \:\: \mathrm{res}_t}VV \\ 
H^2(G \times^P ((w_{\mathfrak k}(t) + \mathfrak{n}(\mathfrak{p})) \times X'), \mathbf{C}) @<{(\phi_{\tilde{w}, t})^*}<< 
H^2(G \times^P ((t + \mathfrak{n}(\mathfrak{p})) \times X'), \mathbf{C})  
\end{CD} 
\end{equation}

By the diagram, we have 
\begin{equation} 
\begin{CD} 
w_{\mathfrak{k}}(t) @>>> t \\ 
@A{\mathrm{res}_0\circ \mathrm{res}_{w_{\mathfrak k}(t)}^{-1}}AA   @A{\mathrm{res}_0\circ \mathrm{res}_t^{-1}}AA \\ 
[\omega_{w_{\mathfrak{k}}(t)}] @>>> [\omega_t] 
\end{CD} 
\end{equation}
Hence we have $$w_{\mathfrak k}^{-1} = \mathrm{res}_0 \circ (\mathfrak{m}_{P}^{-1} \circ \mathfrak{m}_{P'})_*  \circ (\phi_{\tilde{w}})^* 
\circ (\mathrm{res}_0)^{-1}.$$ 
Similarly, for $\sigma$, we have already constructed an isomorphism $$\phi_{\sigma}: 
G \times^{P'}(\mathfrak{r}(\mathfrak{p}') \times X') \to 
G \times^P (\mathfrak{r}(\mathfrak{p}) \times X')$$ in the commutative diagram just before 
Claim \ref{Claim 1}.  

This induces an isomorphism of 2-nd cohomology 
$$(\phi_{\sigma})^*: H^2(G \times^P(\mathfrak{r}(\mathfrak{p}) \times X'), \mathbf{C}) 
\to H^2(G \times^P (\mathfrak{r}(\mathfrak{p}') \times X'), \mathbf{C}).$$
Then we have $$\sigma_{\mathfrak k}^{-1} = \mathrm{res}_0 \circ (\mathfrak{m}_{P}^{-1} \circ \mathfrak{m}_{P'})_*  \circ 
(\phi_{\sigma})^* \circ (\mathrm{res}_0)^{-1}.$$
  
It is sufficient to prove that $(\phi_{\tilde w})^*$ and $(\phi_{\sigma})^*$ coincides. 
If we restrict $\phi_{\tilde w}$ and $\phi_{\sigma}$ to the fibers over $0 \in \mathfrak{l}$, 
then we respectively have isomorphisms $\phi_{{\tilde w},0} : G \times^{P'} (\mathfrak{n}(\mathfrak{p}') \times X') \to G \times^P (\mathfrak{n}(\mathfrak{p}) \times X')$ and $\phi_{\sigma, 0} : 
G \times^{P'} (\mathfrak{n}(\mathfrak{p}') \times X') \to G \times^P (\mathfrak{n}(\mathfrak{p}) \times X')$. 
Since there is a deformation retract of $G \times^P (\mathfrak{r}(\mathfrak{p}) \times X')$ (resp. 
$G \times^{P'}(\mathfrak{r}(\mathfrak{p}') \times X')$) to $G \times^P (\mathfrak{n}(\mathfrak{p}) \times X'), \mathbf{C})$ (resp. $G \times^{P'}(\mathfrak{n}(\mathfrak{p}') \times X')$), 
it is enough to prove that the two maps $$(\phi_{{\tilde w},0})^* : H^2(G \times^{P}(\mathfrak{n}(\mathfrak{p})) \times X'), \mathbf{C}) 
\to H^2(G \times^{P'} (\mathfrak{n}(\mathfrak{p}') \times X'), \mathbf{C})$$ and $$ 
(\phi_{\sigma, 0})^*: H^2(G \times^{P}(\mathfrak{n}(\mathfrak{p})) \times X'), \mathbf{C}) 
\to H^2(G \times^{P'} (\mathfrak{n}(\mathfrak{p}') \times X'), \mathbf{C})$$
coincide. Recall that $\phi_{{\tilde w},0}$ induces an automorphism $\rho_X(\tilde{w}) \in 
\mathrm{Aut}(X/\bar{O})$:     
\begin{equation} 
\begin{CD} 
G \times^{P'} (\mathfrak{n}(\mathfrak{p}') \times X') @>{\phi_{{\tilde w},0}}>> G \times^P (\mathfrak{n}({\mathfrak p}) \times X') \\ 
@V{\mathfrak{m}'_{P',0}}VV @V{\mathfrak{m}'_{P,0}}VV \\ 
X @>{\rho_X(\tilde w)}>> X \\   
\end{CD} 
\end{equation}
By definition $X$ has a $G$-action and $\rho_X(\tilde w)$ is $G$-equivariant. 
The $G$-action on $X$ lifts to a $G$-action on $G \times^{P'} (\mathfrak{n}(\mathfrak{p}') \times X')$ (resp. $G \times^P (\mathfrak{n}({\mathfrak p}) \times X')$) by $\mathfrak{m}'_{P',0}$ (resp. 
$\mathfrak{m}'_{P,0}$). Hence $\phi_{\tilde w}$ is also $G$-equivariant. Recall that, for the 
covering map $\pi: X \to \bar{O}$, the inverse image $\pi^{-1}(0)$ consists of one point $0_X$. 
Since $\rho_X(\tilde{w})(0_X) = 0_X$, $\phi_{{\tilde w},0}$ induces a $G$-equivariant isomorphism between 
${\mathfrak{m}'}_{P',0}^{-1}(0_X) = G/P'$ and ${\mathfrak{m}'}_{P,0}^{-1}(0_X) = G/P$. 
Let $N_G(P)$ be the normalizer of $P$ in $G$.  
Fix an element $g_0 \in G$ such that $P' = g_0Pg_0^{-1}$. Then  
any $G$-equivariant isomorphism $G/P' \to G/P$ is written as $\bar{g} \mapsto \bar{g}\bar{g'_0}$, 
$\bar{g} \in G/P'$ for an element $g'_0 \in g_0N_{G}(P)$. Since $P$ is a parabolic subgroup of 
$G$, we have $N_G(P) = P$. This means that such a $G$-equivariant isomorphism is uniquely 
determined and we denote it by $f_{P', P}$. The inclusion map ${\mathfrak{m}'}_{P',0}^{-1}(0_X)  \to G \times^{P'} (\mathfrak{n}(\mathfrak{p}') \times X')$ induces an isomorphism 
$H^2(G \times^{P'} (\mathfrak{n}(\mathfrak{p}') \times X'), \mathbf{C}) \cong 
H^2(G/P', \mathbf{C})$. Similarly the inclusion map   
${\mathfrak{m}'}_{P,0}^{-1}(0_X)  \to G \times^{P} (\mathfrak{n}(\mathfrak{p}) \times X')$ induces an 
isomorphism $H^2(G \times^{P} (\mathfrak{n}(\mathfrak{p}) \times X'), \mathbf{C}) \cong 
H^2(G/P, \mathbf{C})$. Then we have a commutative diagram 
\begin{equation} 
\begin{CD} 
H^2(G \times^{P} (\mathfrak{n}(\mathfrak{p}) \times X'), \mathbf{C}) @>{(\phi_{{\tilde w},0})^*}>> H^2(G \times^{P'} (\mathfrak{n}({\mathfrak p}') \times X'), \mathbf{C}) \\ 
@V{\cong}VV @V{\cong}VV \\ 
H^2(G/P, \mathbf{C}) @>{(f_{P',P})^*}>> H^2(G/P', \mathbf{C})    
\end{CD} 
\end{equation}
On the other hand, $\phi_{\sigma, 0}$ induces the identity map of $X$ by definition. 
Hence we get a commutative diagram 
\begin{equation} 
\begin{CD} 
G \times^{P'} (\mathfrak{n}(\mathfrak{p}') \times X') @>{\phi_{\sigma, 0}}>> G \times^P (\mathfrak{n}({\mathfrak p}) \times X') \\ 
@V{\mathfrak{m}'_{P',0}}VV @V{\mathfrak{m}'_{P,0}}VV \\ 
X @>{id}>> X 
\end{CD} 
\end{equation}
The diagram is $G$-equivariant and it induces a $G$-equivariant isomorphism between 
${\mathfrak{m}'}_{P',0}^{-1}(0_X) = G/P'$ and ${\mathfrak{m}'}_{P,0}^{-1}(0_X) = G/P$.
As explained above, this $G$-equivariant isomorphism is $f_{P',P}$. We then have a commutative diagram 
\begin{equation} 
\begin{CD} 
H^2(G \times^{P} (\mathfrak{n}(\mathfrak{p}) \times X'), \mathbf{C}) @>{(\phi_{\sigma,0})^*}>> H^2(G \times^{P'} (\mathfrak{n}({\mathfrak p}') \times X'), \mathbf{C}) \\ 
@V{\cong}VV @V{\cong}VV \\ 
H^2(G/P, \mathbf{C}) @>{(f_{P',P})^*}>> H^2(G/P', \mathbf{C})    
\end{CD} 
\end{equation}
By the two diagrams of 2-nd cohomology, we conclude that 
$(\phi_{{\tilde w},0})^* = (\phi_{\sigma, 0})^*$; hence $(\phi_{\tilde w})^* = 
(\phi_{\sigma})^*$. $\square$ \vspace{0.2cm}  


\begin{Claim}\label{Claim 3} One can lift the element $w \in W'$ to an element $\tilde{w} \in  
\widetilde{W'}$ so that $\phi_{\sigma} = \phi_{\tilde{w}}$. In particular, $\tilde{w} \in W_X$. \end{Claim} 

{\em Proof}. 
Take an arbitrary lift $\tilde{w} \in \widetilde{W'}$ of $w \in W$. By Claim \ref{Claim 2} 
we have a commutative diagram  
\begin{equation} 
\begin{CD}  
G \times^P (\mathfrak{r}(\mathfrak{p}) \times X') @>{\phi_{\tilde w} \circ \phi^{-1}_{\sigma}}>>  
G \times^P (\mathfrak{r}(\mathfrak{p}) \times X') \\  
@VVV @VVV \\
Z_{\mathfrak{l}, X'} @>{\tilde{w}_Z \circ \sigma^{-1}_Z}>> Z_{\mathfrak{l}, X'}\\ 
@VVV @VVV \\   
\mathfrak{k} @>{id}>> \mathfrak{k} 
\end{CD}
\end{equation}
Take $t \in \mathfrak{k}^{reg}$. By the proof of Lemma \ref{Lemma 2}, the fiber 
of $G \times^P (\mathfrak{r}(\mathfrak{p}) \times X') \to \mathfrak{k}$ over 
$t$ is isomorphic to $G \times^L (t \times X')$.   
Therefore the map $\phi_{\tilde w} \circ \phi^{-1}_{\sigma}$ induces 
a $G$-equivariant Poisson automorphism 
$$\phi_{\tilde w,t} \circ \phi^{-1}_{\sigma,t}: 
G \times^L (t \times X') \cong G \times^L (t \times X').$$
Let us consider the map $\mu_t: G \times^L (t \times X') \to G\cdot (t + \bar{O}') \subset 
\mathfrak{g}$. Take $x \in G \times^L (t \times X')$ general. Then the $G$-orbit $G\cdot x$ is 
open dense in $G \times^L (t \times X')$. Moreover, let $O_{\mu_t(x)} \subset \mathfrak{g}$ 
be the adjoint orbit containing $\mu_t(x)$. Then $\mu_t\vert_{G\cdot x}: G\cdot x \to 
O_{\mu_t(x)} \subset \mathfrak{g}$ is the moment map for $G\cdot x$. 
Note that $O_{\mu_t(x)}$ is also contained in $G \times^L(t + \bar{O}')$ as the open dense 
$G$-orbit. $\phi_{\tilde w,t} \circ \phi^{-1}_{\sigma,t}$ induces a symplectic automorphism 
of $G\cdot x$. By \cite{B-K}, Theorem 6, it must be a deck transformation of $G\cdot x \to 
O_{\mu_t(x)}$. Let $\tilde{O}'$ be the normalization of $\bar{O}'$. Then this implies that $$\phi_{\tilde w,t} \circ \phi^{-1}_{\sigma,t} \in 
\mathrm{Aut}(G \times^L (t \times X')/G \times^L(t \times \tilde{O}')). $$

Therefore,  
we have
$$\phi_{\tilde w} \circ \phi^{-1}_{\sigma} \in \mathrm{Aut}(Z_{\mathfrak{l}, X'}/Y_{\mathfrak{l}, O'}).$$ 
By taking a suitable lift $\tilde{w} \in \tilde{W}'$ of $w \in W'$, we have 
$$\phi_{\tilde w} \circ \phi^{-1}_{\sigma} = id,$$ which means that $\phi_{\sigma} = 
\phi_{\tilde w}$. By definition, $\phi_{\sigma,0}$ induces the identity map of $Z_{\mathfrak{l}, X', 0} (= X)$. Hence we have $\tilde{w} \in W_X$. $\square$ \vspace{0.2cm}

Claim \ref{Claim 2} and Claim \ref{Claim 3} imply that $W(X) \subset W_X$. 

We next prove that $W_X \subset W(X)$. We assume that $W(X) \subsetneq W_X$ and derive 
a contradiction. Take an element $w \in W_X - W(X)$.  For simplicity we write $C_Q$ for the ample cone $\mathrm{Amp}(G \times^Q (\mathfrak{n} \times X'))$. Since $\mathrm{Mov}(G \times^Q (\mathfrak{n} \times X'))$ is a fundamental domain for the $W(X)$-action, we can take an element $w_0 \in W(X)$ so that $w_0(w\cdot C_Q)  = C_{w_0w(Q)} \subset \mathrm{Mov}(G \times^Q (\mathfrak{n} \times X')).$ By changing $w$ by $w_0w$, we may assume, from the first, that $$C_{w(Q)} \subset \mathrm{Mov}(G \times^Q (\mathfrak{n} \times X')).$$ 
Since $w \in W_X$, we can take a lift $\tilde{w} \in \tilde{W}'$ so that the the following 
diagram commutes 
\begin{equation}
\begin{CD} 
G \times^{Q} (\mathfrak{n} \times X') @>{\phi_{\tilde{w},0}}>> 
G \times^{Ad_w(Q)} (Ad_w(\mathfrak{n}) \times X') \\ 
@V{\mathfrak{m}'_{Q,0}}VV @V{\mathfrak{m}'_{w(Q),0}}VV  \\
X @>{id}>> X
\end{CD}
\end{equation}  
By the diagram $\mathfrak{m}'_{Q,0}$ and $\mathfrak{m}'_{w(Q),0}$ give the same $\mathbf{Q}$-factorization 
of $X$. This means that $C_Q$ and $C_{w(Q)}$ are the same chambers inside 
$\mathrm{Mov}(G \times^{Q} (\mathfrak{n} \times X'))$. On the other hand, 
since $w \in W'$ and $w \ne id$, we see that $C_{Q}$ and 
$C_{w(Q)}$ are different chambers. This is a contradiction.
$\square$ \vspace{0.2cm}

As in Proposition \ref{Proposition 10}, we identify $\mathfrak{k}$ with $H^2(G \times^Q (\mathfrak{r} \times X'), \mathbf{C})$ and 
put $\mathfrak{k}_{\mathbf R} = H^2(G \times^Q (\mathfrak{r} \times X'), \mathbf{R})$ 
by the identification. By Theorem 7, (1), $\mathfrak{k}_{\mathbf R}$ is decomposed into 
finite number of chambers.  \vspace{0.2cm} 

\begin{Cor}\label{Corollary 11} The Weyl group $W_X$   
acts on $\mathfrak{k}_{\mathbf R}$ so that each chamber is sent to another chamber. 
The movable cone $\mathrm{Mov}(G \times^Q (\mathfrak{n} \times X'))$ is a fundamental 
domain of the $W_X$-action. $\mathrm{Mov}(G \times^Q (\mathfrak{n} \times X'))$ contains 
exactly $\sharp \mathcal{S}(L)/ \vert W_X \vert$ chambers. In particular, $X$ has exactly 
$\sharp \mathcal{S}(L)/ \vert W_X \vert$ different {\bf Q}-factorial terminalizations. \end{Cor} 

{\em Proof}. These follow from \cite{Na 4}, Main Theorem. We remark that the second statement 
was first pointed out in \cite{BPW}, Proposition 2.19. $\square$      
\vspace{0.2cm}

\begin{Rem}\label{Remark 12} \end{Rem} We can compute the number $N$ of the conjugacy classes in 
$\mathcal{S}(L)$ by using the equivalence relation \cite{Na 5}, Definition 1, (ii) p. 95 in the set of marked Dynkin diagrams. Then, by \cite{Na 1}, Proposition 2.2, we have $$\sharp \mathcal{S}(L) = N \cdot \vert W' \vert,$$ 
where $W'$ is the finite group defined by the Equation (\ref{W'}). The group $W'$ was studied in \cite{Ho}. 
\vspace{0.2cm}

\begin{center}
\S \; {\bf 2}. {\em The map $\rho_X$}
\end{center}  

\begin{Prop}\label{Proposition 13} The map $\rho_X: \tilde{W}' \to \mathrm{Aut}(X/\bar{O})$ is 
surjective. \end{Prop}\vspace{0.2cm}

{\em Proof}. The group $\mathrm{Aut}(X/\bar{O})$ acts on the universal Poisson deformation 
$Z_{\mathfrak{l}, X'}/W_X \to \mathfrak{k}/W_X$ equivariantly. 
$\mathrm{Aut}(X/\bar{O})$ does not necessarily act effectively on the base space 
$\mathfrak{k}/W_X$. Let $K \subset \mathrm{Aut}(X/\bar{O})$ be the subgroup 
consisting of the elements which act trivially on $\mathfrak{k}/W_X$. Let us consider the surjection $\pi: \tilde{W}' \to W'$. 
By $\pi$, $W_X \subset \tilde{W}'$ is mapped isomorphically onto $\pi(W_X) \subset W'$. 
Then $W'/\pi(W_X)$ acts on $\mathfrak{k}/W_X$. 
The injection $\tilde{W}'/W_X \to \mathrm{Aut}(X/\bar{O})$ induces a map 
$W'/\pi(W_X) \to \mathrm{Aut}(X/\bar{O})/K$. We shall prove that this is an isomorphism. 
This is analogous to \cite{Lo}, Proposition 4.7. If we identify $\mathfrak{g}$ with $\mathfrak{g}^*$ by the Killing form, $G\cdot(\mathfrak{r} + \bar{O}')$ is a $\mathbf{C}^*$-invariant subvariety of $\mathfrak{g}^*$. Look at the moment map $$Z_{\mathfrak{l}, X'}/W_X \to 
G\cdot(\mathfrak{r} + \bar{O}') \subset \mathfrak{g}^*$$ for the Poisson $G$-variety 
$Z_{\mathfrak{l}, X'}/W_X$. The moment map is unique up to an element $(\mathfrak{g}^* \otimes_{\mathbf C} \mathbf{C}[\mathfrak{k}/W_X])^G$, where $G$ acts trivially on 
$\mathbf{C}[\mathfrak{k}/W_X]$. Since $(\mathfrak{g}^*)^G = \{0\}$, we see that 
the moment map is unique.  
Since the $\mathrm{Aut}(X/\bar{O})$ action on $Z_{\mathfrak{l}, X'}/W_X$ commutes with the $G$-action, the moment map induces an injection 
$$\mathfrak{g} \to \mathbf{C}[Z_{\mathfrak{l}, X'}/W_X]^{\mathrm{Aut}(X/\bar{O})}.$$

 
Let $\bar{t}, \bar{t}' \in \mathfrak{k}^{reg}/W_X$ be $\mathrm{Aut}(X/\bar{O})/K$-conjugate, and let $t, t' \in \mathfrak{k}^{reg}$ be their lifts. Look at the map $Z_{\mathfrak{l}, X'}/W_X 
\to \mathfrak{k}/W_X$. Then the fibers of the map over $\bar{t}$, $\bar{t}'$ are respectively 
identified with $Z_{\mathfrak{l}, X', t}$ and $Z_{\mathfrak{l}, X', t'}$. The $\mathrm{Aut}(X/\bar{O})$-action on $Z_{\mathfrak{l}, X'}/W_X$ induces a $G$-equivariant Poisson automorphism $Z_{\mathfrak{l}, X', t} \cong Z_{\mathfrak{l}, X', t'}$. By composing the inclusion 
$\mathfrak{g} \subset \mathbf{C}[Z_{\mathfrak{l}, X'}/W_X]$ with the surjection to 
$\mathbf{C}[Z_{\mathfrak{l}, X', t}]$, we get a map 
$\mathfrak{g} \to \mathbf{C}[Z_{\mathfrak{l}, X', t}]$. Similarly we have a map 
$\mathfrak{g} \to \mathbf{C}[Z_{\mathfrak{l}, X', t'}]$.  
Since $\mathrm{Aut}(X/\bar{O})$ acts on $\mathfrak{g}$ trivially, we have a commutative 
diagram 
\begin{equation}
\begin{CD} 
\mathbf{C}[Z_{\mathfrak{l}, X', t}] @>>> \mathbf{C}[Z_{\mathfrak{l}, X', t'}]\\ 
@AAA @AAA \\ 
\mathfrak{g} @>{id}>> \mathfrak{g}
\end{CD}
\end{equation} 
This induces a commutative diagram of the moment maps 
\begin{equation}
\begin{CD} 
Z_{\mathfrak{l}, X', t'} @>>> Z_{\mathfrak{l}, X', t}\\ 
@V{\mu_{t'}}VV @V{\mu_t}VV \\ 
\mathfrak{g}^* @>{id}>> \mathfrak{g}^*
\end{CD}
\end{equation} 

In particular, $\mu_t(Z_{\mathfrak{l}, X', t})$ and $\mu_{t'}(Z_{\mathfrak{l}, X', t'})$ are 
the closure of the same coadjoint orbit of $\mathfrak{g}^*$. This means that 
$t$ and $t'$ are $W'$-conjugate.   

We next prove that $K = \mathrm{Aut}((X')^0/O')$. Since $K$ acts on $\mathfrak{k}/W_X$ 
trivially, $K$ acts on $Z_{\mathfrak{l}, X'} = \mathfrak{k} \times_{\mathfrak{k}/W_X}
Z_{\mathfrak{l}, X'}/W_X$ in such a way that it acts on $\mathfrak{k}$ trivially. Then an 
element of $K$ determines a $G$-equivariant Poisson automorphism of  every fiber of the map $Z_{\mathfrak{l}, X'} \to \mathfrak{k}$. By the same argument as in Claim 3 in the 
proof of Proposition 10, we see that $$K = \mathrm{Aut}(Z_{\mathfrak{l}, X'}/Y_{\mathfrak{l}, O'}) 
= \mathrm{Aut}((X')^0/O').$$ 

Finally we look at the commutative diagram 
\begin{equation}
\begin{CD}
1 @. 1 \\
@VVV @VVV \\
\mathrm{Aut}((X')^0/O') @>>> K\\ 
@VVV @VVV \\ 
\tilde{W}'/W_X  @>>> \mathrm{Aut}(X/\bar{O}) \\ 
@VVV @VVV \\
W'/\pi(W_X) @>>> \mathrm{Aut}(X/\bar{O})/K \\ 
@VVV @VVV \\
1 @. 1 
\end{CD}
\end{equation} 
The 1-st horizontal map and the 3-rd horizontal map are both isomorphisms 
by the observations above. Then the middle horizontal map is an isomorphism. 
$\square$  \vspace{0.2cm}

By Corollary \ref{Corollary 11}, Remark \ref{Remark 12} and Proposition \ref{Proposition 13}, we have: 

\begin{Thm}\label{Theorem 14} Let $N$ be the number of conjugacy classes in $\mathcal{S}(L)$. 
Then the number of different {\bf Q}-factorial terminalizations of $X$ is $$\frac{N \cdot \vert \mathrm{Aut}(X/\bar{O})\vert}{\vert \mathrm{Aut}((X')^0/O')\vert}.$$ \end{Thm} 
\vspace{0.2cm}

We put $\mathfrak{k}_{\mathbf R} := H^2(G \times^Q (\mathfrak{r} \times X'), \mathbf{R})$.  
Then $\mathfrak{k}_{\mathbf R} \otimes \mathbf{C} = \mathfrak{k}$. 
$\mathfrak{k}_{\mathbf R}$ is divided into finite number of chambers by finite linear 
spaces $\{H_i\}_{i \in I}$ of codimension $1$ (= hyperplanes), and each chamber corresponds to an ample cone of $G \times^{Q'} (\mathfrak{r}' \times X')$ for $Q' \in \mathcal{S}(L)$. 
An element $w$ of $W'$ is a reflection if $w$ has order $2$ and the fixed points set of $w$ on 
$\mathfrak{k}_{\mathbf R}$ is a hyperplane. Such a hyperplane is contained in the set 
$\{H_i\}_{i \in I}$. 
Recall that there is an exact sequence 
$$ 1 \to \mathrm{Aut}((X')^0/O') \to \tilde{W'} \to W' \to 1$$
The composite $\mathrm{Aut}((X')^0/O') \to \tilde{W'} \stackrel{\rho_X}\to 
\mathrm{Aut}(X/\bar{O})$ is an injection. Put $Y := Y_{\mathfrak{l}, O',0}$. 
Then $\pi: X \to \bar{O}$ factorizes as $X \to Y \to \bar{O}$ and 
$\mathrm{Aut}((X')^0/O')$ is identified with $\mathrm{Aut}(X/Y) \subset 
\mathrm{Aut}(X/\bar{O})$ by the injection. 
Take a reflection $w \in W'$ and its lift $\tilde{w} \in \tilde{W}'$. Consider 
the subgroup $G(w) \subset \tilde{W}'$ generated by $\tilde{w}$ and 
$\mathrm{Aut}((X')^0/O')$. Note that $G(w)$ is independent of the choice of 
the lift $\tilde{w}$. There are two cases. 

(a) $G(w) \cap W_X \ne \{1\}$. 

(b) $G(w) \cap W_X = \{1\}$. \vspace{0.2cm}

In the case (a), we have $\rho_X(G(w)) = \mathrm{Aut}(X/Y)$.

\begin{Lem}\label{Lemma 15} Assume that $\mathfrak{g}$ is classical and assume that 
$X \to \bar{O}$ is the finite covering associated with the universal covering of $O$.  In the case (b), $\rho_X(G(w)) = \mathrm{Aut}(X/Y_w)$ for 
a finite covering $Y_w \to \bar{O}$ with the following properties. 

\em $\pi$ is factorized into intermediate coverings $$X \to Y \stackrel{f_w}\to Y_w \to \bar{O}$$ with $\mathrm{deg}(f_w) = 2$. \end{Lem} \vspace{0.2cm}



{\em Proof}. Let $H \subset \mathfrak{k}_{\mathbf R}$ be the fixed points set 
of $w$. Take two chambers which share a codimension 1 face $F$ lying on $H$. 
These chambers correspond to 
$\mathrm{Amp}(G \times^Q (\mathfrak{n} \times X'))$ and $\mathrm{Amp}(G \times^{Q'} (\mathfrak{n}' \times X'))$ for some $Q, Q' \in \mathcal{S}(L)$. The two parabolic subgroups are related by $Q = Ad_w(Q')$.  
The face $F$ determines two birational morphisms respectively from $G \times^Q (\mathfrak{n} \times X')$ and $G \times^{Q'} (\mathfrak{n}' \times X')$ to the same target space.  
The target space is described as follows.  
Take a parabolic subgroup $\bar{Q}$ of $G$ in such a way that 
$Q \subset \bar{Q}$, $Q' \subset \bar{Q}$ and $b_2(G/\bar{Q}) = b_2(G/Q) - 1 
= b_2(G/Q') - 1$. Let $\bar{L}$ be the Levi part of $\bar{Q}$ such that 
$T \subset \bar{L}$. We put $\bar{l} := Lie(\bar{L})$ and $\bar{\mathfrak n} := 
\mathfrak{n}(\bar{\mathfrak q})$.  Define $O_{\bar{\mathfrak l}} := \mathrm{Ind}^{\bar{\mathfrak l}}_{\mathfrak l}(O')$. Let $X_{\bar{\mathfrak l}}$ be the Stein factorization of the composite  
$$\bar{L} \times^{\bar{L} \cap Q}(\mathfrak{n}(\bar{\mathfrak l} \cap \mathfrak{q}) 
\times X') \to \bar{L} \times^{\bar{L} \cap Q}(\mathfrak{n}(\bar{\mathfrak l} \cap \mathfrak{q}) 
+ \bar{O'}) \to \bar{O}_{\bar{\mathfrak l}}.$$  Then the target space is $G \times^{\bar Q}(\bar{\mathfrak{n}} \times X_{\bar{\mathfrak l}})$.    

There is a  commutative diagram 
\begin{equation} 
\begin{CD}   
G \times^{Q'} (\mathfrak{n}' \times X') @>{\mathfrak{m}'_{Q', \bar{Q},0}}>> G \times^{\bar Q}(\bar{\mathfrak{n}} \times X_{\bar{\mathfrak l}}) @<{\mathfrak{m}'_{Q, \bar{Q},0}}<< G \times^{Q}(\mathfrak{n} \times X') \\
@VVV @VVV @VVV \\ 
G \times^{Q'} (\mathfrak{n}' + \bar{O}') @>{\mathfrak{m}_{Q', \bar{Q},0}}>> G \times^{\bar Q}(\mathfrak{n} + \bar{O}_{\bar{\mathfrak l}}) @<{\mathfrak{m}_{Q, \bar{Q},0}}<< G \times^{Q}(\mathfrak{n} + \bar{O}')
\end{CD}
\end{equation} such that $\mathfrak{m}'_{Q', \bar{Q},0}$ and $\mathfrak{m}'_{Q, \bar{Q},0}$ 
are the birational morphisms determined by $F$.  

Note that $w \in W(\bar{L})$ for the Levi part $\bar{L}$ of $\bar{Q}$. Then, since $w \in \bar{Q}$, we have a commutative diagram 
\begin{equation}
\begin{CD} 
G \times^{Q'} (\mathfrak{n}' \times X') @>{\phi_{\tilde{w},0}}>> 
G \times^Q (\mathfrak{n} \times X') \\ 
@VVV @VVV  \\
G \times^{\bar Q} (\bar{\mathfrak{n}} + \bar{O}_{\bar{\mathfrak l}})^n @>{id}>> 
G \times^{\bar Q} (\bar{\mathfrak{n}} + \bar{O}_{\bar{\mathfrak l}})^n
\end{CD}
\end{equation}

We put here $$Y_w := \mathrm{Spec}\: \Gamma (G \times^{\bar Q} (\bar{\mathfrak{n}} + \bar{O}_{\bar{\mathfrak l}})^n, \mathcal{O})$$ 
Since $Y = \mathrm{Spec} \Gamma (G \times^{Q}(\mathfrak{n} + \bar{O}')^n, 
\mathcal{O})$, the map $\mathfrak{m}_{Q, \bar{Q},0}$ induces a map $f_w: Y \to Y_w$ and 
the finite cover $Y \to \bar{O}$ factors through $Y_w$.  
The map $\mathfrak{m}_{Q, \bar{Q},0}$ is a generalized Springer map 
determined by an {\em induction} of type I or of type II (see [Na, Part I, \S 2, \S 3]).  
Then $\mathrm{deg}(\mathfrak{m}_{Q, \bar{Q},0}) = 1$ or $\mathrm{deg}(\mathfrak{m}_{Q, \bar{Q},0}) = 2$ by [ibid, (2.3.1), (2.3.2), (3.6.1) and 
(3.6.2)]. We used here the assumption that $\mathfrak{g}$ is classical. In particular, $\mathrm{deg}(f_w) = 1, \; 2$.       
The map $\phi_{\tilde{w},0}: G \times^{Q'} (\mathfrak{n}' \times X') 
\to G \times^P (\mathfrak{n} \times X')$ induces automorphisms of 
$X$, $Y$ and $Y_w$. Hence we get a commutative 
diagram 
\begin{equation}
\begin{CD} 
G \times^{Q'} (\mathfrak{n}' \times X') @>{\phi_{\tilde{w},0}}>> 
G \times^Q (\mathfrak{n} \times X') \\ 
@VVV @VVV  \\
X @>{\tilde{w}_{Z,0}}>> X \\ 
@VVV @VVV \\ 
Y @>{\tilde{w}_{Y,0}}>> Y \\
@VVV @VVV \\  
Y_w @>{id}>> Y_w
\end{CD}
\end{equation}   
If $\tilde{w}_{Y,0} = id$, then  
$\tilde{w}_{Z,0}$ is an element of $\mathrm{Aut}(X/Y)$. 
This means that $G(w) \cap W_X \ne \{1\}$, which contradicts the assumption. 
Therefore $\tilde{w}_{Y,0} \ne id$. Then $\rho_X(\tilde{w}) \in \mathrm{Aut}(X/Y_w)$, but 
$\rho_X(\tilde{w}) \notin \mathrm{Aut}(X/Y)$. In particular, $\mathrm{Aut}(X/Y_w) \ne 
\mathrm{Aut}(X/Y)$; hence, $f_w: Y \to Y_w$ is a double cover. This means that 
$\rho_X(G(w)) = \mathrm{Aut}(X/Y_w)$. $\square$   \vspace{0.2cm} 

\begin{Rem}\label{Remark 16} \end{Rem} (1) The proof of Lemma \ref{Lemma 15} shows that if Case (b) occurs, then  
$\mathrm{deg}\: \mathfrak{m}_{Q, \bar{Q}, 0} = 2$. 

(2) 
If $\mathfrak{m}'_{Q, \bar{Q},0}$ is divisorial, then Case (a) occurs. If $\mathfrak{m}'_{Q, \bar{Q},0}$ is a small birational map (i.e. $\mathrm{Codim}\:\mathrm{Exc}(\mathfrak{m}'_{Q, \bar{Q},0}) \geq 2$), 
then Case (b) occurs.  

In fact, if $\mathfrak{m}'_{Q, \bar{Q},0}$ is divisorial, then $\mathfrak{m}_{Q, \bar{Q},0}$ is birational by the next Proposition \ref{Proposition 17}. Then $Y_w = Y$; hence 
$\tilde{w}_{Y,0} = id$. 
If $\mathfrak{m}'_{Q, \bar{Q},0}$ is small, then
$\mathfrak{m}'_{Q', \bar{Q},0}$ is also small. Therefore $\mathrm{Amp}(G \times^{Q'}(\mathfrak{n}' \times X'))$ is contained in $\mathrm{Mov}(G \times^Q(\mathfrak{n} \times X')$. Since $Q = Ad_w(Q')$ and $w \ne 1$, the two cones 
$\mathrm{Amp}(G \times^Q(\mathfrak{n} \times X'))$ and $\mathrm{Amp}(G \times^{Q'}(\mathfrak{n}' \times X'))$ determine different chambers inside  
$\mathrm{Mov}(G \times^Q(\mathfrak{n} \times X'))$. We prove that $\tilde{w}_{Y,0} \ne id$ for any lift $\tilde{w}$ of $w$. In fact, if $\tilde{w}_{Y,0} = id$, then one can assume that $\tilde{w}_{Z,0} = id$ by replacing $\tilde{w}$ by a suitable lift 
of $w$. Then we have a commutative diagram 
\begin{equation}
\begin{CD} 
G \times^{Q'} (\mathfrak{n}' \times X') @>{\phi_{\tilde{w},0}}>> 
G \times^Q (\mathfrak{n} \times X') \\ 
@VVV @VVV  \\
X @>{id}>> X
\end{CD}
\end{equation}
By the diagram $G \times^{Q'} (\mathfrak{n}' \times X')$ and 
$G \times^Q (\mathfrak{n} \times X')$ give the same 
{\bf Q}-factorial terminalization of $X$. This means that 
$$\mathrm{Amp}(G \times^Q(\mathfrak{n} \times X')) = \mathrm{Amp}(G \times^{Q'}(\mathfrak{n}' \times X'))$$
This is a contradiction. $\square$ 
\vspace{0.2cm}

Take a reflection $w \in W'$. Let $H \subset \mathfrak{t}_{\mathbf R}$ be the fixed hyperplane of $w$.  
As in the proof of Lemma \ref{Lemma 15}, we find $Q \in \mathcal{S}(L)$ in such a way 
that, if we put $Q' := Ad_w(Q)$, then the ample cones $\mathrm{Amp}(G \times^Q (\mathfrak{n} \times X'))$ 
and $\mathrm{Amp}(G \times^{Q'}(\mathfrak{n}' \times X')$ 
share a codimension 1 face $F$ lying on $H$. 
The codimension 1 face $F$ determines a birational contraction map of $G \times^Q (\mathfrak{n} \times X')$. 
The birational map is described as follows.  
There is a parabolic subgroup $\bar{Q}$ of $G$ so that $Q \subset \bar{Q}$ and    
$b_2(G/\bar{Q}) = b_2(G/Q) - 1$. Let $\bar{\mathfrak l}$ be the 
Levi part of $\bar{\mathfrak{q}}$.  
There are a nilpotent orbit $O_{\bar{\mathfrak l}} \subset \mathfrak{l}$, a finite covering 
$X_{\bar{\mathfrak l}} \to \bar{O}_{\bar{\mathfrak l}}$ and a commutative diagram 
\begin{equation} 
\begin{CD}
G \times^Q (\mathfrak{n} \times X') @>{\mathfrak{m}'_{Q, \bar{Q},0}}>> G \times^{\bar Q}(\bar{\mathfrak n} \times X_{\bar{\mathfrak l}})  \\ 
@V{\pi'_Q}VV  @V{\bar{\pi}}VV \\
G \times^Q (\mathfrak{n} + \bar{O}') @>{\mathfrak{m}_{Q, \bar{Q},0}}>> G \times^{\bar Q}(\bar{\mathfrak n} + \bar{O}_{\bar{\mathfrak l}})
\end{CD} 
\end{equation}
such that $\mathfrak{m}'_{Q, \bar{Q},0}$ is a birational map. The map $\mathfrak{m}'_{Q, \bar{Q},0}$ is the desired birational contraction map. 
Note that $\mathfrak{m}_{Q, \bar{Q},0}$ 
is not necessarily a birational map. 
\vspace{0.2cm}

\begin{Prop}\label{Proposition 17} Assume that $\mathfrak{g}$ is a classical simple Lie algebra and assume that $X \to \bar{O}$ is the finite covering associated with the universal covering of $O$.  

(1) If $\mathfrak{m}'_{Q, \bar{Q},0}$ is a divisorial 
birational map, then $\mathfrak{m}_{Q, \bar{Q},0}$  
is a divisorial birational map. 

(2) If $\mathfrak{m}'_{Q, \bar{Q},0}$ is a small birational map, then $\mathrm{deg} \: \mathfrak{m}_{Q, \bar{Q},0} = 2$. \end{Prop}
\vspace{0.2cm} 

{\em Proof}.  (2) follows from Remark \ref{Remark 16}. In fact, if $\mathfrak{m}'_{Q, \bar{Q},0}$ is a small birational map, then 
Case (b) occurs by Remark \ref{Remark 16}, (2). By Remark \ref{Remark 16}, (1) we then have $\mathrm{deg} \: \mathfrak{m}_{Q, \bar{Q},0} = 2$. 
In the remainder we prove (1). To prove (1) we only have to show that $\mathfrak{m}_{Q, \bar{Q},0}$ is 
birational if $\mathfrak{m}'_{Q, \bar{Q},0}$ is a divisorial 
birational map.  In fact, if $\mathfrak{m}'_{Q, \bar{Q},0}$ and $\mathfrak{m}_{Q, \bar{Q},0}$ are both birational, then 
one is divisorial if and only if the other is divisorial.   

When $\mathfrak{g} = sl(n)$, the map $\mathfrak{m}_{Q, \bar{Q},0}$ is always birational. 
When $\mathfrak{g} = sp(2n)$ or $so(n)$, we prove that $\mathrm{Codim}\:\: \mathrm{Exc}(\mathfrak{m}'_{Q, \bar{Q},0}) \geq 2$ if $\mathrm{deg}(\mathfrak{m}_{Q, \bar{Q},0}) \geq 2$.   

Case (a): $\mathfrak{g} = sp(2n)$

In \cite{Part 1}, \S 2, we constructed a {\bf Q}-factorial terminalization of $X$. 
For such a {\bf Q}-factorial terminalization, $O'$ is contained in a simple factor $\mathfrak{l}_0$ of the Levi part $\mathfrak{l}$, and $\mathfrak{l}_0$ is of type $C_k$.   
The Jordan type $\mathbf{p} = [d^{r_d}, 
(d-1)^{r_{d-1}}, ..., 2^{r_2}, 1^{r_1}]$ of $O'$ is then a partition of $2k$ which satisfies the following conditions 

(i) $r_i$ is even for each odd $i$, 

(ii) $r_i \ne 0$ for each even $i$, and 

(iii) $r_i \ne 2$ for each even $i$.  

Moreover, $X' \to \bar{O}'$ is the finite cover associated with the universal covering 
of $O'$. Since the core $(X', O')$ is independent of the choice of a {\bf Q}-factorial terminalization of $X$, we can start with this $(X', O')$. 

There is a unique simple factor $\bar{\mathfrak l}_0$ of the Levi part $\bar{\mathfrak l}$ such that $\mathfrak{l}_0 \subset \bar{\mathfrak l}_0$. Then $\bar{\mathfrak l}_0$ is again of type 
$C$ and $$O_{\bar{\mathfrak l}} = \mathrm{Ind}^{\bar{\mathfrak l}_0}_{\mathfrak{l}_0}(O').$$ 
If $\mathrm{deg}(\mathfrak{m}_{Q, \bar{Q},0}) \geq 2$, then it is an induction of type II
(cf. \cite{Part 1}, \S 2] and $\mathrm{deg}(\mathfrak{m}_{Q, \bar{Q},0}) = 2$. Let $\bar{\mathfrak p} := [{\bar d}^{{\bar r}_{\bar d}}, (\bar{d}-1)^{\bar{d}_{\bar{r}-1}}, ..., 
2^{\bar{r}_2}, 1^{\bar{r}_1}]$ be the Jordan type of $O_{\bar{\mathfrak l}}$. 
Two partitions $\mathbf{p}$ and $\bar{\mathbf p}$ are related as follows. 
There is an even member $i$ of $\bar{\mathbf p}$ with $\bar{r}_i = 2$, and $\bar{p}$ can be written as 
$$\mathbf{p} = [(\bar{d} - 2)^{\bar{r}_{\bar d}}, ..., i^{\bar{r}_{i+2}}, (i-1)^{\bar{r}_{i+1} + 2 
+ \bar{r}_{i-1}}, (i-2)^{\bar{r}_{i-2}}, ..., 1^{\bar{r}_i}]. $$  
It is checked that $\bar{\mathbf p}$ satisfies the conditions (i) and (ii) by using the fact that 
$\mathbf{p}$ satisfies (i), (ii) and (iii). 
Let $b$ be the number of distinct even members of $\mathbf{p}$.  
Recall that $\pi_1(O') = (\mathbf{Z}/2\mathbf{Z})^{\oplus b}$. This means that 
$\mathrm{deg} \: \pi'_Q = 2^b$. Since $\mathrm{deg}(\mathfrak{m}_{Q, \bar{Q},0}) = 2$ 
and $\mathfrak{m}'_{Q, \bar{Q},0}$ is birational, we see that 
$\mathrm{deg} \: \bar{\pi} = 2^{b + 1}$. 
Let $\bar{\mathbf p}$ 
be the number of distinct even members of $\bar{\mathbf p}$. Then 
$\pi_1(O_{\bar{\mathfrak l}}) 
= (\mathbf{Z}/2\mathbf{Z})^{\oplus \bar{b}}$. 
By the description of $\mathbf{p}$ above, we have $\bar{b} = b + 1$. These facts imply 
that $X_{\bar{\mathfrak l}} \to \bar{O}_{\bar{\mathfrak l}}$ is the finite covering associated 
with the universal covering of $O_{\bar{\mathfrak l}}$. Since the Jordan type $\bar{\mathbf p}$ of $O_{\bar{\mathfrak l}}$ satisfies the conditions (i) and (ii), we see that $\mathrm{Codim}_{X_{\bar{\mathfrak l}}} \mathrm{Sing}(X_{\bar{\mathfrak l}}) \geq 4$ 
by \cite{Part 1}, Proposition (2,3). This means that $\mathrm{Codim}\:\: \mathrm{Exc}(\mathfrak{m}'_{Q, \bar{Q},0}) \geq 2$. \vspace{0.2cm}

Case (b): $\mathfrak{g} = so(n)$

In \cite{Part 1}, \S 3, we constructed a {\bf Q}-factorial terminalization of $X$. 
Assume that $\mathfrak{l}$ has simple factors $so(k)$, $sl(k_1)$, ..., and $sl(k_s)$.
There are three cases for $(X', O')$. The first case is when $O$ is rather odd, and the remaining cases are when $O$ is not rather odd.  \vspace{0.2cm} 
 
Case (b-1): $k_1$, ..., $k_s$ are all even and $$O' = O_{\mathbf p} \times O_{[2^{\frac{k_1}{2}}]} 
\times ... \times O_{[2^{\frac{k_s}{2}}]} \subset so(k) \oplus sl(k_1) \oplus ... \oplus sl(k_s)$$ 
Here $\mathbf{p} := [d^{r_d}, (d-1)^{r_{d-1}}, ..., 2^{r_2}, 1^{r_1}]$ is a {\em rather odd} partition of $k$ with the following conditions. 

(i) $r_i$ is even for each even $i$

(ii) For any couple $(i,j)$ of adjacent members, $\vert i - j \vert \leq 4$. Moreover 
$\vert i - j \vert = 4$ occurs only when $i$ and $j$ are both odd. The smallest member 
of $\mathbf{p}$ is smaller than $4$. 

Moreover, if we put $m_i := r_d + r_{d-1} + ... + r_i$, then $\frac{k_1}{2}$, ..., $\frac{k_s}{2}$ respectively coincide with some $m_i$'s. 

We shall prove that $\mathrm{deg}(\mathfrak{m}_{Q, \bar{Q},0}) \geq 2$ never happens.   
In fact, there are two possibilities for $\bar{\mathfrak l}$. 
$$\bar{\mathfrak l} = so(k + 2k_t) \oplus sl(k_1) \oplus ... \oplus sl(k_{t-1}) \oplus sl(k_{t+1}) 
\oplus ... sl(k_s) \oplus \mathfrak{z}(\bar{\mathfrak l})$$
$$\bar{\mathfrak l} = so(k) \oplus sl(k_i + k_j) \oplus sl(k_1) \oplus ... \oplus 
\check{sl(k_{i})} \oplus ... \oplus \check{sl(k_{j})} \oplus ... \oplus  
sl(k_s) \oplus \mathfrak{z}(\bar{\mathfrak l}).$$
In the first case, one can write $k_t= 2m_i$ for some $i$. If we put   
$$\bar{\mathbf p} := [(d+4)^{r_d}, (d+3))^{r_{d-1}}, ..., (i+4)^{r_i}, (i-1)^{r_{i-1}}, ..., 2^{r_2}, 1^{r_1}],$$ 
then $$O_{\bar{\mathfrak l}} = O_{\bar{\mathbf p}} \times O_{[2^{\frac{k_1}{2}}]} 
\times ...\times O_{[2^{\frac{k_{t-1}}{2}}]} \times  O_{[2^{\frac{k_{t+1}}{2}}]} \times ... \times   O_{[2^{\frac{k_s}{2}}]}.$$
In the second case, $$O_{\bar{\mathfrak l}} = O_{\mathbf p} \times O_{[2^{\frac{k_i + k_j}{2}}]} 
\times O_{[2^{\frac{k_1}{2}}]} \times ... \times \check{O}_{[2^{\frac{k_i}{2}}]}
\times ... \times \check{O}_{[2^{\frac{k_j}{2}}]} \times ... \times O_{[2^{\frac{k_s}{2}}]}.$$ 
One can check that, in both cases, $\mathfrak{m}_{Q, \bar{Q},0}$ is birational. \vspace{0.2cm}

Case (b-2) : $O' = O_{\mathbf p} \subset so(k)$. Here 
$\mathbf{p} := [d^{r_d}, (d-1)^{r_{d-1}}, ..., 2^{r_2}, 1^{r_1}]$ is a partition of $k$ with the following conditions. 
 
(i) $r_i$ is even for each even $i$,  

(ii) $r_i \ne 0$ for each odd $i$,  

(iii) $r_i \ne 2$ for each odd $i$,  and 

(iv) $r_i \geq 3$ for some odd $i$. 

Moreover, $X' \to \bar{O}'$ is the finite cover associated with the universal covering of $O$. 

For some $t$, we have 
$$\bar{\mathfrak l} = so(k + 2k_t) \oplus sl(k_1) \oplus ... \oplus sl(k_{t-1}) \oplus sl(k_{t+1}) \oplus ... sl(k_s) \oplus \mathfrak{z}(\bar{\mathfrak l})$$ and $$O_{\bar{\mathfrak l}} = \mathrm{Ind}^{so(k + 2k_t)}_{so(k)}(O') \subset so(k + 2k_t).$$
If $\mathrm{deg}(\mathfrak{m}_{Q,\bar{Q},0}) \geq 2$, then it is an induction of type II and 
$\mathrm{deg}(\mathfrak{m}_{Q,\bar{Q},0}) = 2$. Let 
$\bar{\mathbf p} := [{\bar d}^{{\bar r}_{\bar d}}, (\bar{d}-1)^{\bar{d}_{\bar{r}-1}}, ..., 
2^{\bar{r}_2}, 1^{\bar{r}_1}]$ be the Jordan type of $O_{\bar{\mathfrak l}}$. 
Two partitions $\mathbf{p}$ and $\bar{\mathbf p}$ are related as follows. 
There is an odd member $i$ of $\bar{\mathbf p}$ with $\bar{r}_i = 2$, and $\mathbf{p}$ can be written as 
$$\mathbf{p} = [(\bar{d} - 2)^{\bar{r}_{\bar d}}, ..., i^{\bar{r}_{i+2}}, (i-1)^{\bar{r}_{i+1} + 2 
+ \bar{r}_{i-1}}, (i-2)^{\bar{r}_{i-2}}, ..., 1^{\bar{r}_i}]. $$
Here $k_t = \bar{r}_{\bar d} + \bar{r}_{\bar{d}-1} + ... + \bar{r}_{i+1} + 1$. 
It is checked that $\bar{\mathbf p}$ satisfies the conditions (i), (ii) and (iv) by using the fact 
that $\mathbf{p}$ satisfies (i), (ii), (iii) and (iv). 
Let $a$ be the number of distinct odd members of $\mathbf{p}$. 
Since $\mathbf{p}$ is not rather odd and $a >0$ by the condition (iv), we have $\pi_1(O') = (\mathbf{Z}/2\mathbf{Z})^{\oplus a-1}$. This means that $\mathrm{deg}\: \pi'_Q = 2^{a-1}$. 
Since $\mathrm{deg}(\mathfrak{m}_{Q,\bar{Q},0}) = 2$ and $\mathfrak{m}'_{Q, \bar{Q},0}$ is 
birational, we see that $\mathrm{deg}\: \bar{\pi} = 2^a$. Let $\bar{a}$ be the number of distinct odd members of $\bar{\mathbf p}$. By the description of $\bar{\mathbf p}$ above, $\bar{a} = a + 1$. Then $\pi_1(O_{\bar{\mathfrak l}}) = (\mathbf{Z}/2\mathbf{Z})^{\oplus \bar{a}-1}$ because $\bar{\mathbf p}$ is not rather odd. These facts imply that $X_{\mathfrak l} 
\to \bar{O}_{\mathfrak l}$ is the finite covering associated with the universal covering of 
$O_{\mathfrak l}$. Since $\bar{\mathbf p}$ is not rather odd and satisfies the conditions (i) and (ii), we have $\mathrm{Codim}_{X_{\mathfrak l}}\mathrm{Sing}(X_{\mathfrak l}) \geq 4$ by 
\cite{Part 1}, Proposition (3.4). This means that $\mathrm{Codim}\:\: \mathrm{Exc}(\mathfrak{m}'_{Q, \bar{Q},0}) \geq 2$. \vspace{0.2cm} 

Case (b-3): $O' = O_{\mathbf p} \subset so(k)$. Here 
$\mathbf{p} := [d^{r_d}, (d-1)^{r_{d-1}}, ..., 2^{r_2}, 1^{r_1}]$ is a partition of $k$ with the following conditions. 
 
(i) $r_i$ is even for each even $i$, and  

(ii) $r_i = 1$ for each odd $i$
 
In particular, $\mathbf{p}$ is rather odd. Let $a$ be the number of distinct odd members of 
$\mathbf{p}$. Then we have an exact sequence 
$$1 \to \mathbf{Z}/2\mathbf{Z} \to \pi_1(O') \to (\mathbf{Z}/2\mathbf{Z})^{\oplus \mathrm{min}(a-1, 0)} \to 1.$$ The finite cover $X' \to \bar{O}'$ is induced from the 
etale covering of $O'$ determined by the surjection $\pi_1(O') \to (\mathbf{Z}/2\mathbf{Z})^{\oplus \mathrm{min}(a-1, 0)}$. In particular, we have 
$\mathrm{deg} (\pi'_Q) = 2^{\mathrm{min}(a-1, 0)}$.  
For some $t$, we have 
$$\bar{\mathfrak l} = so(k + 2k_t) \oplus sl(k_1) \oplus ... \oplus sl(k_{t-1}) \oplus sl(k_{t+1}) \oplus ... sl(k_s) \oplus \mathfrak{z}(\bar{\mathfrak l})$$ and $$O_{\bar{\mathfrak l}} = \mathrm{Ind}^{so(k + 2k_t)}_{so(k)}(O') \subset so(k + 2k_t).$$
If $\mathrm{deg}(\mathfrak{m}_{Q,\bar{Q},0}) \geq 2$, then it is an induction of type II and 
$\mathrm{deg}(\mathfrak{m}_{Q,\bar{Q},0}) = 2$. Let 
$\bar{\mathfrak p} := [{\bar d}^{{\bar r}_{\bar d}}, (\bar{d}-1)^{\bar{d}_{\bar{r}-1}}, ..., 
2^{\bar{r}_2}, 1^{\bar{r}_1}]$ be the Jordan type of $O_{\bar{\mathfrak l}}$. 
Two partitions $\mathbf{p}$ and $\bar{\mathbf p}$ are related as follows. 
There is an odd member $i$ of $\bar{\mathbf p}$ with $\bar{r}_i = 2$, and $\mathbf{p}$ can be written as 
$$\mathbf{p} = [(\bar{d} - 2)^{\bar{r}_{\bar d}}, ..., i^{\bar{r}_{i+2}}, (i-1)^{\bar{r}_{i+1} + 2 
+ \bar{r}_{i-1}}, (i-2)^{\bar{r}_{i-2}}, ..., 1^{\bar{r}_i}]. $$
Here $k_t = \bar{r}_{\bar d} + \bar{r}_{\bar{d}-1} + ... + \bar{r}_{i+1} + 1$.
It is checked that $\bar{\mathbf p}$ satisfies the following conditions  

(i) $\bar{r}_i$ is even for each even $i$, and 

(iii) $\bar{r}_i \ne 0$ for each odd $i$. 

We may assume that $a > 0$. In fact, if $a = 0$, then $\mathbf{p} = [0]$ because 
$\mathbf{p}$ satisfies the condition (ii). In such a case $\mathfrak{m}_{Q, \bar{Q},0}$ 
is birational by \cite{Part 1}, Claim (3.6.2). Therefore, $\mathrm{deg} (\pi'_Q) = 2^{\mathrm{min}(a-1, 0)} = 2^{a-1}$.
Since $\mathrm{deg}(\mathfrak{m}_{Q,\bar{Q},0}) = 2$ and $\mathfrak{m}'_{Q, \bar{Q},0}$ is 
birational, we see that $\mathrm{deg}\: \bar{\pi} = 2^a$. Let 
$\bar{a}$ be the number of distinct odd members of $\bar{\mathbf p}$. By the description of $\bar{\mathbf p}$ above, $\bar{a} = a + 1$. Then $\pi_1(O_{\bar{\mathfrak l}}) = (\mathbf{Z}/2\mathbf{Z})^{\oplus \bar{a}-1}$ because $\bar{\mathbf p}$ is not rather odd (bacause 
$\bar{r}_i = 2$ for some odd $i$). These facts imply that $X_{\mathfrak l} 
\to \bar{O}_{\mathfrak l}$ is the finite covering associated with the universal covering of 
$O_{\mathfrak l}$. Since $\bar{\mathbf p}$ is not rather odd and satisfies the conditions (i) and (iii), we have $\mathrm{Codim}_{X_{\mathfrak l}}\mathrm{Sing}(X_{\mathfrak l}) \geq 4$ by \cite{Part 1}, 
Proposition (3.4). This means that $\mathrm{Codim}\:\: \mathrm{Exc}(\mathfrak{m}'_{Q, \bar{Q},0}) \geq 2$.  $\square$ \vspace{0.2cm}


The homomorphism  $\rho_X : \tilde{W'} \to \mathrm{Aut}(X/\bar{O})$ induces a homomorphism 
$$\overline{\rho_X}: W' \to \mathrm{Aut}(X/\bar{O})/\mathrm{Aut}((X')^0/O').$$

\begin{Cor}\label{Corollary 18} Assume that $\mathfrak{g}$ is classical and assume that $X \to \bar{O}$ is the finite covering associated with the universal covering of $O$. 
For a reflection $w \in W'$, the map $\mathfrak{m}_{Q, \bar{Q}, 0}$ 
is a divisorial birational map or a generically finite surjective map of degree 2. In the first case, 
$\overline{\rho_X}(w) = 1$, and in the latter case $\overline{\rho_X}(w) \ne 1$. \end{Cor}  

When $\mathfrak{g}$ is classical, $\mathrm{deg}(\mathfrak{m}_{Q, \bar{Q}, 0}) = 1$ or $\mathrm{deg}(\mathfrak{m}_{Q, \bar{Q}, 0}) = 2$. 
Moreover, one can check which case occurs by using the criterion \cite{Part 1}, (2.3.1), (2.3.2), (3.6.1), (3.6.2). 
Then the Weyl group $W_X$ is the subgroup of $W'$ generated by all reflections $w$ for which  
$\mathrm{deg}(\mathfrak{m}_{Q, \bar{Q}, 0}) = 1$. \vspace{0.2cm}

\begin{center}
\S \; {\bf 3}. Examples
\end{center}

\begin{Example}\label{Example 19} (\cite{Part 1}, Example (2.5), (1)) : \end{Example}   
Let us consider the nilpotent orbit $O_{[6^2,4^2]} \subset 
sp(20)$ and let $\pi: X \to \bar{O}_{[6^2,4^2]}$ be the finite covering associated with 
the universal covering of $O_{[6^2,4^2]}$. We have $\mathrm{deg}(\pi) = 4$ because 
the partition $[6^2,4^2]$ has exactly two distinct even members.  
We can construct a {\bf Q}-factorial terminalization of $X$ 
by the induction step 
$$([1^4], sp(4)) \stackrel{\mathrm{Type II}}\to (O_{[3^2,2^2]}, sp(10)) \stackrel{\mathrm{Type II}}\to 
(O_{[4^2, 2^2]}, sp(12)) \stackrel{\mathrm{Type I}}\to (O_{[6^2,4^2]}, sp(20)).$$
We consider the three generalized Springer maps 
$$Sp(20) \times^{Q_{4,12,4}}(\mathfrak{n}(\mathfrak{q}_{4,12,4}) + \bar{O}_{[4^2, 2^2]}) \to 
\bar{O}_{[6^2,4^2]},$$
$$Sp(20) \times^{Q_{4,1,10,1,4}}(\mathfrak{n}(\mathfrak{q}_{4,1,10,1,4}) + 
\bar{O}_{[3^2,2^2]}) \to Sp(20) \times^{Q_{4,12,4}}(\mathfrak{n}(\mathfrak{q}_{4,12,4}) + \bar{O}_{[4^2, 2^2]}),$$ and  
$$Sp(20) \times^{Q_{4,1,3,4,3,1,4}} \mathfrak{n}(\mathfrak{q}_{4,1,3,4,3,1,4}) \to 
Sp(20) \times^{Q_{4,1,10,1,4}}(\mathfrak{n}(\mathfrak{q}_{4,1,10,1,4}) + 
\bar{O}_{[3^2,2^2]}).$$
Composing these 3 maps together, we have a generically 
finite map of degree 4 
$$Sp(20) \times^{Q_{4,1,3,4,3,1,4}} \mathfrak{n}(\mathfrak{q}_{4,1,3,4,3,1,4}) \to \bar{O}_{[6^2,4^2]}.$$ 
This map factors through $X$ and $Sp(20) \times^{Q_{4,1,3,4,3,1,4}} \mathfrak{n}(\mathfrak{q}_{4,1,3,4,3,1,4})$ 
gives a crepant resolution of $X$. We shall calculate how many different crepant resolutions 
$X$ has. 

The parabolic subgroup $Q_{4,1,3,4,3,1,4}$ corresponds to 
the marked Dynkin diagram 

\begin{picture}(300,20)
\put(30,0){\circle{5}}\put(35,0){\line(1,0){20}}
\put(60,0){\circle{5}}\put(65,0){\line(1,0){20}}
\put(90,0){\circle{5}}\put(95,0){\line(1,0){20}}
\put(120,0){\circle*{5}}\put(125,0){\line(1,0){20}}
\put(150,0){\circle*{5}}\put(155,0){\line(1,0){20}}
\put(180,0){\circle{5}}\put(185,0){\line(1,0){20}}
\put(210,0){\circle{5}}\put(215,0){\line(1,0){20}}
\put(240,0){\circle*{5}}\put(245,0){\line(1,0){20}}
\put(270,0){\circle{5}}\put(275,-3){$\Longleftarrow$}
\put(300,0){\circle{5}}
\end{picture} 
\vspace{0.5cm}

Let $L$ be the Levi part of $Q_{4,1,3,4,3,1,4}$. One can read from the marked 
Dynkin diagram that the semi simple part of $\mathfrak{l}$ is $A_3$ + $A_2$ + $C_2$. 
By \cite{Ho}, p.71, we see that $$W' = (\mathbf{Z}/2\mathbf{Z})^{\oplus 3}.$$  
$\mathcal{S}(L)$ contains  
exactly 6 conjugacy classes. Other 5 conjugacy classes correspond to the following 
marked Dynkin diagrams (cf. Remark 11)
 
\begin{picture}(300,20)
\put(30,0){\circle*{5}}\put(35,0){\line(1,0){20}}
\put(60,0){\circle{5}}\put(65,0){\line(1,0){20}}
\put(90,0){\circle{5}}\put(95,0){\line(1,0){20}}
\put(120,0){\circle{5}}\put(125,0){\line(1,0){20}}
\put(150,0){\circle*{5}}\put(155,0){\line(1,0){20}}
\put(180,0){\circle{5}}\put(185,0){\line(1,0){20}}
\put(210,0){\circle{5}}\put(215,0){\line(1,0){20}}
\put(240,0){\circle*{5}}\put(245,0){\line(1,0){20}}
\put(270,0){\circle{5}}\put(275,-3){$\Longleftarrow$}
\put(300,0){\circle{5}}
\end{picture}
\vspace{0.5cm}
  
\begin{picture}(300,20)
\put(30,0){\circle*{5}}\put(35,0){\line(1,0){20}}
\put(60,0){\circle{5}}\put(65,0){\line(1,0){20}}
\put(90,0){\circle{5}}\put(95,0){\line(1,0){20}}
\put(120,0){\circle*{5}}\put(125,0){\line(1,0){20}}
\put(150,0){\circle{5}}\put(155,0){\line(1,0){20}}
\put(180,0){\circle{5}}\put(185,0){\line(1,0){20}}
\put(210,0){\circle{5}}\put(215,0){\line(1,0){20}}
\put(240,0){\circle*{5}}\put(245,0){\line(1,0){20}}
\put(270,0){\circle{5}}\put(275,-3){$\Longleftarrow$}
\put(300,0){\circle{5}}
\end{picture}
\vspace{0.5cm}

\begin{picture}(300,20)
\put(30,0){\circle{5}}\put(35,0){\line(1,0){20}}
\put(60,0){\circle{5}}\put(65,0){\line(1,0){20}}
\put(90,0){\circle*{5}}\put(95,0){\line(1,0){20}}
\put(120,0){\circle*{5}}\put(125,0){\line(1,0){20}}
\put(150,0){\circle{5}}\put(155,0){\line(1,0){20}}
\put(180,0){\circle{5}}\put(185,0){\line(1,0){20}}
\put(210,0){\circle{5}}\put(215,0){\line(1,0){20}}
\put(240,0){\circle*{5}}\put(245,0){\line(1,0){20}}
\put(270,0){\circle{5}}\put(275,-3){$\Longleftarrow$}
\put(300,0){\circle{5}}
\end{picture}
\vspace{0.5cm}

\begin{picture}(300,20)
\put(30,0){\circle{5}}\put(35,0){\line(1,0){20}}
\put(60,0){\circle{5}}\put(65,0){\line(1,0){20}}
\put(90,0){\circle*{5}}\put(95,0){\line(1,0){20}}
\put(120,0){\circle{5}}\put(125,0){\line(1,0){20}}
\put(150,0){\circle{5}}\put(155,0){\line(1,0){20}}
\put(180,0){\circle{5}}\put(185,0){\line(1,0){20}}
\put(210,0){\circle*{5}}\put(215,0){\line(1,0){20}}
\put(240,0){\circle*{5}}\put(245,0){\line(1,0){20}}
\put(270,0){\circle{5}}\put(275,-3){$\Longleftarrow$}
\put(300,0){\circle{5}}
\end{picture}
\vspace{0.5cm}

\begin{picture}(300,20)
\put(30,0){\circle{5}}\put(35,0){\line(1,0){20}}
\put(60,0){\circle{5}}\put(65,0){\line(1,0){20}}
\put(90,0){\circle{5}}\put(95,0){\line(1,0){20}}
\put(120,0){\circle*{5}}\put(125,0){\line(1,0){20}}
\put(150,0){\circle{5}}\put(155,0){\line(1,0){20}}
\put(180,0){\circle{5}}\put(185,0){\line(1,0){20}}
\put(210,0){\circle*{5}}\put(215,0){\line(1,0){20}}
\put(240,0){\circle*{5}}\put(245,0){\line(1,0){20}}
\put(270,0){\circle{5}}\put(275,-3){$\Longleftarrow$}
\put(300,0){\circle{5}}
\end{picture}
\vspace{0.5cm}

It follows from this fact that $$N = 6, \:\: \sharp \mathcal{S}(L) = 6 \vert W' \vert = 48.$$
In our case $X' = \bar{O'} = \{0\}$. 
On the other hand, $\mathrm{Aut}(X/\bar{O}_{[6^2,4^2]}) = (\mathbf{Z}/2\mathbf{Z})^{\oplus 2}$ (cf. \cite{Part 1}, Proposition (2.1)). 
Then, by Theorem 13,  
we have $$\sharp \{\mathrm{crepant}\: \mathrm{resolutions}\: \mathrm{of}\: X\} = N \cdot \vert \mathrm{Aut}(X/\bar{O}_{[6^2,4^2]}) 
\vert = 6\cdot 4  = 24.$$  

We next closely study the map $\rho_X: \tilde{W'} \to \mathrm{Aut}(X/\bar{O}_{[6^2,4^2]})$. 
Look at the first marked Dynkin diagram 

\begin{picture}(300,20)
\put(30,0){\circle{5}}\put(35,0){\line(1,0){20}}
\put(60,0){\circle{5}}\put(65,0){\line(1,0){20}}
\put(90,0){\circle{5}}\put(95,0){\line(1,0){20}}
\put(120,0){\circle*{5}}\put(125,0){\line(1,0){20}}
\put(150,0){\circle*{5}}\put(155,0){\line(1,0){20}}
\put(180,0){\circle{5}}\put(185,0){\line(1,0){20}}
\put(210,0){\circle{5}}\put(215,0){\line(1,0){20}}
\put(240,0){\circle*{5}}\put(245,0){\line(1,0){20}}
\put(270,0){\circle{5}}\put(275,-3){$\Longleftarrow$}
\put(300,0){\circle{5}}
\put(235, -9){$\alpha$}
\end{picture} 
\vspace{0.5cm}

Let $Q'_{4,1,3,4,3,1,4}$ be a parabolic subgroup of $Sp(20)$ obtained from 
$Q_{4,1,3,4,3,1,4}$ by twisting at $\alpha$. Then $Q_{4,1,3,4,3,1,4} = Ad_{w_1}(Q'_{4,1,3,4,3,1,4})$ for a (uniquely determined) element $w_1 \in W'$. Let $\bar{Q}_{4,1,10,1,4}$ be rhe parabolic subgroup of $Sp(20)$ with flag type $(4,1,10,1,4)$ such that $Q_{4,1,3,4,3,1,4} \subset 
\bar{Q}_{4,1,10,1,4}$ and $Q'_{4,1,3,4,3,1,4} \subset 
\bar{Q}_{4,1,10,1,4}$. Put $Q_I = Q_{4,1,3,4,3,1,4}$, $Q'_I = Q'_{4,1,3,4,3,1,4}$ and $\bar{Q}_I = \bar{Q}_{4,1,10,1,4}$. 
Let $\mathfrak{n}_I$, $\mathfrak{n}'_I$ and $\bar{\mathfrak{n}}_I$ be respectively the nilradicals of $Q_I$, $Q'_I$ and $\bar{Q}_I$. Then we have a commutative diagram 
\begin{equation} 
\begin{CD}    
Sp(20) \times^{Q'_I} \mathfrak{n}'_I  @>{\mathfrak{m}'_{Q'_I, \bar{Q}_I, 0}}>> Sp(20) \times^{\bar{Q}_I}(\bar{\mathfrak{n}}_I \times X_{[3^2,2^2]}) @<{\mathfrak{m}'_{Q_I, \bar{Q}_I,0}}<< Sp(20) \times^{Q_I} \mathfrak{n}_I \\
@VVV @VVV @VVV \\ 
Sp(20) \times^{Q'_I} \mathfrak{n}'_I  @>{\mathfrak{m}_{Q'_I, \bar{Q}_I,0}}>> Sp(20) \times^{\bar{Q}_I}(\bar{\mathfrak{n}}_I + \bar{O}_{[3^2,2^2]}) @<{\mathfrak{m}_{Q_I, \bar{Q}_I,0}}<< Sp(20) \times^{Q_I} \mathfrak{n}_I
\end{CD}
\end{equation}
The horizontal maps on the first row are birational maps, and the horizontal maps on the second row have degree 2 because they are determined by a type II induction. Then the birational maps on the first row are both small  
by Proposition \ref{Proposition 17}. We put $$Y_{w_1} := \mathrm{Spec}\: \Gamma (Sp(20) \times^{\bar{Q}_I}(\bar{\mathfrak{n}}_I + \bar{O}_{[3^2,2^2]}), \mathcal{O}).$$
$Y_{w_1}$ is a double covering of $\bar{O}_{[6^2,4^2]}$.    

Next let $Q_{4,3,1,4,1,3,4}$ be the parabolic subgroup corresponding to the 
marked Dynkin diagram 

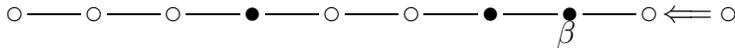
\begin{picture}(300,20)
\put(30,0){\circle{5}}\put(35,0){\line(1,0){20}}
\put(60,0){\circle{5}}\put(65,0){\line(1,0){20}}
\put(90,0){\circle{5}}\put(95,0){\line(1,0){20}}
\put(120,0){\circle*{5}}\put(125,0){\line(1,0){20}}
\put(150,0){\circle{5}}\put(155,0){\line(1,0){20}}
\put(180,0){\circle{5}}\put(185,0){\line(1,0){20}}
\put(210,0){\circle*{5}}\put(215,0){\line(1,0){20}}
\put(240,0){\circle*{5}}\put(245,0){\line(1,0){20}}
\put(270,0){\circle{5}}\put(275,-3){$\Longleftarrow$}
\put(300,0){\circle{5}}
\put(235, -11){$\beta$}
\end{picture} 
\vspace{0.5cm}

Let $Q'_{4,3,1,4,1,3,4}$ be a parabolic subgroup of $Sp(20)$ obtained from 
$Q_{4,3,1,4,1,3,4}$ by twisting at $\beta$. Then $Q_{4,3,1,4,1,3,4} = Ad_{w_2}(Q'_{4,3,1,4,1,3,4})$ for a (uniquely determined) element $w_2 \in W'$. Let $\bar{Q}_{4,3,6,3,4}$ be rhe parabolic subgroup of $Sp(20)$ with flag type $(4,3,6,3,4)$ such that $Q_{4,3,1,4,1,3,4} \subset 
\bar{Q}_{4,3,6,3,4}$ and $Q'_{4,3,1,4,1,3,4} \subset 
\bar{Q}_{4,3,6,3,4}$. Put $Q_{II}  = Q_{4,3,1,4,1,3,4}$, $Q'_{II} = Q'_{4,3,1,4,1,3,4}$ and $\bar{Q}_{II} = \bar{Q}_{4,3,6,3,4}$. 
Let $\mathfrak{n}_{II}$, $\mathfrak{n}'_{II}$ and $\bar{\mathfrak{n}}_{II}$ be respectively the nilradicals of $Q_{II}$, $Q'_{II}$ and $\bar{Q}_{II}$.
Then we have a commutative diagram 
\begin{equation} 
\begin{CD}   
Sp(20) \times^{Q'_{II}} \mathfrak{n}'_{II}  @>{\mathfrak{m}'_{Q'_{II}, \bar{Q}_{II},0}}>> Sp(20) \times^{\bar{Q}_{II}}(\bar{\mathfrak{n}}_{II} \times  X_{[2^2,1^2]}) @<{\mathfrak{m}'_{Q_{II}, \bar{Q}_{II},0}}<< Sp(20) \times^{Q_{II}} \mathfrak{n}_{II} \\
@VVV @VVV @VVV \\ 
Sp(20) \times^{Q'_{II}} \mathfrak{n}'_{II}  @>{\mathfrak{m}_{Q'_{II}, \bar{Q}_{II},0}}>> Sp(20) \times^{\bar{Q}_{II}}(\bar{\mathfrak{n}}_{II} + \bar{O}_{[2^2, 1^2]}) @<{\mathfrak{m}_{Q, \bar{Q}_{II},0}}<< Sp(20) \times^{Q_{II}} \mathfrak{n}_{II}
\end{CD}
\end{equation}
The horizontal maps on the first row are birational maps, and the horizontal maps on the second row have degree 2 because they  are determined by a type II induction. Then the birational maps on the first row are both small  
by Proposition \ref{Proposition 17}. We put $$Y_{w_2} := \mathrm{Spec}\: \Gamma (Sp(20) \times^{\bar{Q}_{II}}(\bar{\mathfrak{n}}_{II} + \bar{O}_{[2^2,1^2]}), \mathcal{O}).$$
$Y_{w_2}$ is a double covering of $\bar{O}_{[6^2,4^2]}$.    

In our case, $Y = X$ and $\tilde{W'} = W'$. By the next lemma $\mathrm{Aut}(X/Y_{\omega_1})$ and 
$\mathrm{Aut}(X/Y_{\omega_2})$ are mutually different, index 2, subgroups of 
$\mathrm{Aut}(X/\bar{O}_{[6^2,4^2]}) = (\mathbf{Z}/2\mathbf{Z})^{\oplus 2}$. 
By Lemma 14, $\mathrm{Im}(\rho_X)$ contains both $\mathrm{Aut}(X/Y_{\omega_i})$. 
This reproves that $$\mathrm{Im}(\rho_X) = \mathrm{Aut}(X/\bar{O}_{[6^2,4^2]}).$$  
In particular, $$W_X = \mathrm{Ker}(\rho_X) = \mathbf{Z}/2\mathbf{Z}.$$

\begin{Lem}\label{Lemma 20} The double coverings $\pi_1: Y_{\omega_1} \to \bar{O}_{[6^2,4^2]}$ and 
$\pi_2: Y_{\omega_2} \to \bar{O}_{[6^2,4^2]}$ are different coverings in the function 
field $\mathbf{C}(X)$ of $X$.  \end{Lem}

{\em Proof}. $\bar{O}_{[6^2,4^2]}$ has $A_3$-surface singularity along the orbit $O_{[6^2, 4, 2^2]}$. This means that the transverse slice $S \subset \bar{O}_{[6^2,4^2]}$ for $O_{[6^2, 4, 2^2]}$ is a complex analytic germ of an $A_3$-singularity. We claim that $\pi_1^{-1}(S)$ 
is isomorphic to a disjoint union of two copies of $S$, but $\pi_2^{-1}(S)$ is isomorphic to 
an $A_1$-germ. We first look at the generalized Springer map 
$$f: Sp(20) \times^{Q_{4,12,4}} (\mathfrak{n}(\mathfrak{q}_{4,12,4}) + \bar{O}_{[4^2,2^2]}) 
\to \bar{O}_{[6^2,4^2]}$$ Since this map is determined by the induction of type I, 
it is a birational map. The closed subvariety $Sp(20) \times^{Q_{4,12,4}} (\mathfrak{n}(\mathfrak{q}_{4,12,4}) + \bar{O}_{[4^2,2,1^2]}) \subset Sp(20) \times^{Q_{4,12,4}} (\mathfrak{n}(\mathfrak{q}_{4,12,4}) + \bar{O}_{[4^2,2^2]})$ is mapped onto the subvariety $\bar{O}_{[6^2,4,2^2]} \subset 
\bar{O}_{[6^2,4^2]}$ by $f$. The map $$Sp(20) \times^{Q_{4,12,4}} (\mathfrak{n}(\mathfrak{q}_{4,12,4}) + \bar{O}_{[4^2,2,1^2]}) \to \bar{O}_{[6^2,4,2^2]}$$ is a generically finite map of degree 2 because 
it is determined by the type II induction. Note that $Sp(20) \times^{Q_{4,12,4}} (\mathfrak{n}(\mathfrak{q}_{4,12,4}) + \bar{O}_{[4^2,2^2]})$ has $A_1$-surface singularity along    
$Sp(20) \times^{Q_{4,12,4}} (\mathfrak{n}(\mathfrak{q}_{4,12,4}) + O_{[4^2,2,1^2]})$. Therefore we 
conclude that $f^{-1}(S)$ is a crepant partial resolution of the $A_3$-singularity with one 
exceptional curve $C \cong \mathbf{P}^1$ and $f^{-1}(S)$ has two $A_1$-singularities 
at two points $p_{+}, p_{-} \in C$.  We next look at the generalized Springer maps 
$$g_1: Sp(20) \times^{\bar{Q}_{4,1,10,1,4}}(\mathfrak{n}(\bar{\mathfrak{q}}_{4,1,10,1,4}) + \bar{O}_{[3^2,2^2]}) \to 
Sp(20) \times^{Q_{4,12,4}} (\mathfrak{n}(\mathfrak{q}_{4,12,4}) + \bar{O}_{[4^2,2^2]})$$ 
$$g_2: Sp(20) \times^{\bar{Q}_{4,3,6,3,4}}(\mathfrak{n}(\bar{\mathfrak{q}}_{4,3,6,3,4}) + \bar{O}_{[2^2, 1^2]}) 
\to Sp(20) \times^{Q_{4,12,4}} (\mathfrak{n}(\mathfrak{q}_{4,12,4}) + \bar{O}_{[4^2,2^2]})$$ and 
observe $g_1^{-1}(p)$ and $g_2^{-1}(p)$ for a point $p \in C$.   
The map $g_1$ can be described as 
$$Sp(20) \times^{\bar{Q}_{4,1,10,1,4}}(\mathfrak{n}(\bar{\mathfrak{q}}_{4,1,10,1,4}) + \bar{O}_{[3^2,2^2]}) =  
Sp(20) \times^{Q_{4,12,4}}(Q_{4,12,4} \times^{\bar{Q}_{4,1,10,1,4}} (\mathfrak{n}(\bar{\mathfrak{q}}_{4,1,10,1,4}) + 
\bar{O}_{[3^2,2^2]}))$$ $$= Sp(20) \times^{Q_{4,12,4}}\{\mathfrak{n}(\mathfrak{q}_{4,12,4}) 
\times (L(Q_{4,12,4}) \times^{L(Q_{4,12,4}) \cap \bar{\mathfrak q}_{4,1,10,1,4}} (\mathfrak{n}(\mathfrak{l}(\mathfrak{q}_{4,12,4}) \cap \bar{\mathfrak q}_{4,1,10,1,4}) + \bar{O}_{[3^2,2^2]}))\}$$
$$\longrightarrow Sp(20) \times^{Q_{4,12,4}}\{\mathfrak{n}(\mathfrak{q}_{4,12,4}) \times L(Q_{4,12,4}) \cdot 
(\mathfrak{n}(\mathfrak{l}(\mathfrak{q}_{4,12,4}) \cap \bar{\mathfrak q}_{4,1,10,1,4}) + \bar{O}_{[3^2,2^2]})\}$$ 
$$= Sp(20) \times^{Q_{4,12,4}}(\mathfrak{n}(\mathfrak{q}_{4,12,4}) + \bar{O}_{[4^2,2^2]})$$
Therefore $g_1$ is obtained from the map 
$$L(Q_{4,12,4}) \times^{L(Q_{4,12,4}) \cap \bar{Q}_{4,1,10,1,4}} (\mathfrak{n}(\mathfrak{l}(\mathfrak{q}_{4,12,4}) \cap \bar{\mathfrak q}_{4,1,10,1,4}) + \bar{O}_{[3^2,2^2]})$$
$$\to L(Q_{4,12,4}) \cdot 
(\mathfrak{n}(\mathfrak{l}(\mathfrak{q}_{4,12,4}) \cap \bar{\mathfrak q}_{4,1,10,1,4}) + \bar{O}_{[3^2,2^2]})$$ 
by taking $Sp(20) \times^{Q_{4,12,4}}(\mathfrak{n}(\mathfrak{q}_{4,12,4}) \times \cdot)$. 
This map is nothing but the generalized Springer map 
$$\bar{g}_1: Sp(12) \times^{Q_{1,10,1}} (\mathfrak{n}(\mathfrak{q}_{1,10,1}) + \bar{O}_{[3^2,2^2]}) 
\to \bar{O}_{[4^2,2^2]}.$$
Similarly, $g_2$ is obtained from the generalized Springer map 
$$\bar{g}_2: Sp(12) \times^{Q_{3,6,3}} (\mathfrak{n}(\mathfrak{q}_{3,6,3}) + \bar{O}_{[2^2,1^2]}) 
\to \bar{O}_{[4^2,2^2]}.$$
Write $p \in Sp(20) \times^{Q_{4,12,4}} (\mathfrak{n}(\mathfrak{q}_{4,12,4}) + \bar{O}_{[4^2,2^2]})$ 
as $p = [h, z + x]$ with $h \in Sp(20)$, $z \in \mathfrak{n}(\mathfrak{q}_{4,12,4})$ and $x \in \bar{O}_{[4^2,2^2]}$. Then $$g_1^{-1}(p) = [h, z + \bar{g}_1^{-1}(x)] \cong \bar{g}_1^{-1}(x), \:\: g_2^{-1}(p) = [h, z + \bar{g}_2^{-1}(x)] \cong \bar{g}_2^{-1}(x).$$

{\bf Case 1}: $p = p_+$ or $p = p_{-}$.  

In this case we have $x \in O_{[4^2,2,1^2]}$. 
In order to see the fibers $\bar{g}_1^{-1}(x)$ and $\bar{g}_2^{-1}(x)$, we introduce a basis $\{e(i,j)\}$ of the symplectic vector space $\mathbf{C}^{12}$ so that 
 
(a) $\{e(i,j)\}$ is a Jordan basis for $x$, i.e. 
$x\cdot e(i,j) = e(i-1,j)$ for $i > 1$ and $x\cdot e(1,j) = 0$. 

(b) $\langle e(i,j), e(q,r) \rangle \ne 0$ if and only if $q = d_j -i + 1$ and 
$r = \beta(j)$. Here $\beta$ is a permutation of $\{1,2, ..., s\}$ such that 
$\beta^2 = id$, $d_{\beta(j)} = d_j$, and $\beta (j) \ne j$ if $d_j$ is odd 
(cf. \cite{S-S}, p.259, see also \cite{C-M}, 5.1).  

In our case $d_1 = d_2 = 4$, $d_3 = 2$ and $d_4 = d_5 = 1$. 
Hence $\beta (4) = 5$, $\beta (5) = 4$ and $\beta (3) = 3$. After a suitable basis 
change of the subspace $\Sigma_{1 \leq i \leq 4,\: 1 \leq j \leq 2}\mathbf{C}e(i,j)$, 
we may also assume that $\beta (1) = 2$ and $\beta (2) = 1$.   

\begin{picture}(200, 120)(0, 0)
\put(00,  00){\line(1, 0){120}}
\put(00,  20){\line(1, 0){120}}
\put(2, 5){e(1,1)} 
\put(32, 5){e(2,1)}
\put(62, 5){e(3,1)} 
\put(92, 5){e(4,1)}
\put(2, 25){e(1,2)}
\put(32, 25){e(2,2)}
\put(62, 25){e(3,2)} 
\put(92, 25){e(4,2)}
\put(00, 40){\line(1, 0){120}}
\put(00, 60){\line(1, 0){60}} 
\put(00, 80){\line(1, 0){30}} 
\put(00, 100){\line(1,0){30}}
\put(00, 100){\line(0,-1){100}}
\put(30, 100){\line(0, -1){100}} 
\put(60, 60){\line(0,-1){60}}
\put(90, 40){\line(0, -1){40}} 
\put(120, 40){\line(0, -1){40}}
\put(2, 45){e(1,3)}
\put(32, 45){e(2,3)}
\put(2, 65){e(1,4)}
\put(2, 85){e(1,5)}
\end{picture}

\vspace{1.0cm}

If we put $F_0 = \mathbf{C}e(1,1)$, then $F_0 \subset F_0^{\perp}$ is an isotropic flag 
such that $x \cdot F_0 = 0$, $x \cdot \mathbf{C}^{12} \subset F_0^{\perp}$, and 
$x$ is an endomorphism of $F_0^{\perp}/F_0$ with Jordan type $[3^2, 2, 1^2]$. 
On the other hand, $F'_0 := \mathbf{C}e(1,2)$ and $(F'_0)^{\perp}$ also satisfy the 
same property. These two flags are only those flags of type $(1,10,1)$ which have the properties  

(i) $x \cdot F= 0$, $x \cdot \mathbf{C}^{12} \subset F^{\perp}$, and 

(ii) $x$ is an endomorphism of $F^{\perp}/F = \mathbf{C}^{10}$ with 
$x \in \bar{O}_{[3^2,2^2]}$. 

This means that $\bar{g}_1^{-1}(p)$ consists of two points. Moreover, these two points 
are lying on the locus $$Sp(12) \times^{Q_{1,10,1}} (\mathfrak{n}(\mathfrak{q}_{1,10,1}) + {O}_{[3^2,2, 1^2]}) \subset Sp(12) \times^{Q_{1,10,1}} (\mathfrak{n}(\mathfrak{q}_{1,10,1}) + \bar{O}_{[3^2,2^2]})$$ 
Note that $Sp(12) \times^{Q_{1,10,1}} (\mathfrak{n}(\mathfrak{q}_{1,10,1}) + \bar{O}_{[3^2,2^2]})$ has 
$A_1$ singularities along $Sp(12) \times^{Q_{1,10,1}} (\mathfrak{n}(\mathfrak{q}_{1,10,1}) + {O}_{[3^2,2, 1^2]})$.  

We next consider the fiber $\bar{g}_2^{-1}(p)$. 
If we put $F_0 := \mathbf{C}e(1,1) \oplus \mathbf{C}e(1,2) \oplus \mathbf{C}e(1,3)$, 
then $F_0 \subset F_0^{\perp}$ is an isotropic flag 
such that $x \cdot F_0 = 0$, $x \cdot \mathbf{C}^6 \subset F_0^{\perp}$, and 
$x$ is an endomorphism of $F_0^{\perp}/F_0$ with Jordan type $[2^2, 1^2]$.
This is the only isotropic flag of type $(3, 6, 3)$ which has the properties  

(i) $x \cdot F= 0$, $x \cdot \mathbf{C}^6 \subset F^{\perp}$, and 

(ii) $x$ is an endomorphism of $F^{\perp}/F = \mathbf{C}^{6}$ with 
$x \in \bar{O}_{[2^2,1^2]}$.

Therefore $\bar{g}_2^{-1}(p)$ consists of one point. \vspace{0.2cm}

{\bf Case 2}: $p \ne p_{+}, p_{-}$. 

In this case $x \in O_{[4^2, 2^2]}$. Since $\bar{g}_1$ and $\bar{g}_2$ are both 
generalized Springer map made from type II inductions, we conclude that 
$\bar{g}_1^{-1}(p)$ and $\bar{g}_2^{-1}(p)$ respectively consist of two points. 
\vspace{0.2cm}

By these observations, we see that \vspace{0.2cm} 

1. $(f \circ g_1)^{-1}(S)$ is isomorphic to a disjoint union of two 
copies of $f^{-1}(S)$. \vspace{0.2cm} 

2. $ (f \circ g_2)^{-1}(S)$ is isomorphic to the minimal resolution 
of an $A_1$-singularity. The map $(f \circ g_2)^{-1}(S) \to f^{-1}(S)$ is a 
double covering only ramified at $p_+$ and $p_{-}$. Here $g_2^{-1}(C)$ corresponds to the 
exceptional curve of the minimal resolution of an $A_1$-singularity. 
\vspace{0.2cm} 

Let us consider the commutative diagrams
\begin{equation} 
\begin{CD}   
Sp(20) \times^{\bar{Q}_{4,1,10,1,4}}(\mathfrak{n}(\bar{\mathfrak{q}}_{4,1,10,1,4}) + \bar{O}_{[3^2,2^2]})  @>{g_1}>> 
Sp(20) \times^{Q_{4,12,4}} (\mathfrak{n}(\mathfrak{q}_{4,12,4}) + \bar{O}_{[4^2,2^2]}) \\
@V{f_1}VV @V{f}VV \\ 
Y_{\omega_1}  @>{\pi_1}>> \bar{O}_{[6^2,4^2]}
\end{CD}
\end{equation}

\begin{equation} 
\begin{CD}   
Sp(20) \times^{\bar{Q}_{4,3,6,3,4}}(\mathfrak{n}(\bar{\mathfrak{q}}_{4,3,6,3,4}) + \bar{O}_{[2^2,1^2]})  @>{g_2}>> 
Sp(20) \times^{Q_{4,12,4}} (\mathfrak{n}(\mathfrak{q}_{4,12,4}) + \bar{O}_{[4^2,2^2]}) \\
@V{f_2}VV @V{f}VV \\ 
Y_{\omega_2}  @>{\pi_2}>> \bar{O}_{[6^2,4^2]}
\end{CD}
\end{equation}


Here $(f \circ g_1)^{-1}(S)$ is mapped to a disjoint union of two copies of an 
$A_3$-singularity by $f_1$. On the other hand, $(f \circ g_2)^{-1}(S)$ is mapped 
to an $A_1$-singularity by $f_2$.  $\square$ 
\vspace{0.2cm}

\begin{Example}\label{Example 21}  (\cite{Part 1}, Example (3.10), (2)): \end{Example} 
Let $X$ be a finite covering of $\bar{O}_{[11^3, 3^2, 1]} \subset so(40)$ 
associated with the universal covering of $O_{[11^3, 3^2, 1]}$. 
By inductions of type I and of type II 
$$[11^3, 3^2, 1] \stackrel{Type I}\leftarrow [9^3, 3^2, 1] \stackrel{Type I}\leftarrow  
[7^3, 3^2, 1] \stackrel{Type I}\leftarrow [5^3, 3^2, 1] \stackrel{Type II}\leftarrow  
[3^3, 2^2, 1]$$ we finally get a partition $[3^3, 2^2, 1]$ of $14$. 
Let $\bar{Q} \subset SO(40)$ be a parabolic subgroup stabilizing an 
isotropic flag of type $(3,3,3,4,14,4,3,3,3)$ and put $Q := \rho_{40}^{-1}(\bar{Q})$ for 
$\rho_{40}: Spin (40) \to SO(40)$. Let us consider the nilpotent orbit 
$O_{[3^3, 2^2, 1]} \subset so(14)$. 
Let $X_{[3^3,2^2,1]} \to \bar{O}_{[3^3, 2^2, 1]}$ be a finite covering associated with 
the universal covering of $O_{[3^3, 2^2, 1]}$.  
Then $$Spin (40) \times^Q (\mathfrak{n} \times X_{[3^3,2^2,1]})$$ is a {\bf Q}-factorial 
terminalization of $X$. Fix a Borel subgroup $B$ of $Spin(40)$ with $B \subset Q$ and a  maximal torus $T \subset B$. We denote by $\{\alpha_i\}_{1 \le i \le 20}$ the simple 
roots determined by $B$. 
Then the parabolic subgroup $Q$ corresponds to the marked Dynkin diagram of type $D_{20}$: 
\vspace{0.5cm}
   
\begin{picture}(400,20)
\put(30,0){\circle{5}}\put(35,0){\line(1,0){10}}\put(50,0){\circle{5}} 
\put(55,0){\line(1,0){10}}\put(70,0){\circle*{5}}\put(65,5){$\alpha_3$}
\put(75,0){\line(1,0){10}}    
\put(90,0){\circle{5}}
\put(95,0){\line(1,0){10}}\put(110,0){\circle{5}}\put(115,0){\line(1,0){10}}
\put(130,0){\circle*{5}}\put(125,5){$\alpha_6$}
\put(135,0){\line(1,0){10}}\put(150,0){\circle{5}}
\put(155,0){\line(1,0){10}}\put(170,0){\circle{5}}\put(175,0){\line(1,0){10}} 
\put(190,0){\circle*{5}}\put(185,5){$\alpha_9$} 
\put(195,0){\line(1,0){10}}\put(210,0){\circle{5}}
\put(215,0){\line(1,0){10}}\put(230,0){\circle{5}}\put(235,0){\line(1,0){10}} 
\put(250,0){\circle{5}}\put(255,0){\line(1,0){10}}\put(270,0){\circle*{5}}\put(265,5){$\alpha_{13}$} 
\put(275,0){\line(1,0){10}}\put(290, -3){- - -}\put(310,0){\line(1,0){10}}
\put(325,0){\circle{5}}\put(327,0){\line(1,1){10}}\put(340,12){\circle{5}}
\put(327,0){\line(1,-1){10}}\put(340,-12){\circle{5}}
\end{picture}
\vspace{1.0cm}

Put $\mathfrak{t} := Lie(T)$ and define $\mathfrak{t}^*_{\mathbf R} \subset \mathfrak{t}^*$  to be the real vector space spanned by all simple roots. Let $I := \{\alpha_i\}_{1 \le i \le 20, i \ne 3, 6, 9, 13}$ and denote by $\Phi_I$ the root subsystem of $\Phi$ generated by $I$. 
The semisimple part of the Levi part $L$ of $Spin(40)$ is $3A_2 + 
A_3 + D_7$ and its root system is nothing but $\Phi_I$. 
By [Ho, p.71] we have $$W' = W(B_3) \times W(B_1).$$

We shall explicitly describe $W'$ in terms of twists. 
  
We have 
$$H^2(Spin(40)/Q, \mathbf{R}) \cong \mathrm{Hom}_{alg.gp}(L, \mathbf{C}^*)\otimes_{\mathbf Z}\mathbf{R} = \{x \in \mathfrak{t}^*_{\mathbf R}; (x, \alpha_i) = 0 \;\; \mathrm{for}\;\; \mathrm{all}\;\; \alpha_i \in I\},$$
where $(\;, \;)$ is the Killing form on $\mathfrak{t}^*_{\mathbf R}$. 
By the Killing form we identify $\mathfrak{t}^*_{\mathbf R}$ with 
$\mathfrak{t}_{\mathbf R} (= \mathrm{dual}\; \mathrm{space}\; \mathrm{of}\; 
\mathfrak{t}^*_{\mathbf R})$.    
Thus $\mathrm{Hom}_{alg.gp}(L, \mathbf{C}^*)\otimes_{\mathbf Z}\mathbf{R}$ can be identified with the subspace of $\mathfrak{t}_{\mathbf R}$ consisting of the elements 
which kill all $\alpha_i \in I$.     

We twist the $Q$ by $\alpha_3$: 
 
\begin{picture}(400,20) 
\put(30,0){\circle{5}}\put(25,5){$\alpha_1$}\put(35,0){\line(1,0){20}}\put(60,0){\circle{5}}
\put(55,5){$\alpha_2$} 
\put(65,0){\line(1,0){20}}\put(90,0){\circle*{5}}\put(85,5){$\alpha_3$}
\put(95,0){\line(1,0){20}}    
\put(120,0){\circle{5}}\put(115,5){$\alpha_4$} 
\put(125,0){\line(1,0){20}}\put(150,0){\circle{5}}\put(145,5){$\alpha_5$} 
\put(155,0){\line(1,0){20}}
\put(180,0){\circle*{5}}\put(175,5){$\alpha_6$}
\put(185,0){\line(1,0){10}}\put(200,-3){- - -} 
\end{picture}  \vspace{1.0cm}

\begin{picture}(400,20) 
\put(30,0){\circle{5}}\put(25,5){$-\alpha_5$}\put(35,0){\line(1,0){20}}\put(60,0){\circle{5}}
\put(55,5){$-\alpha_4$} 
\put(65,0){\line(1,0){20}}\put(90,0){\circle*{5}}\put(85,5){$-\alpha_3$}
\put(95,0){\line(1,0){20}}    
\put(120,0){\circle{5}}\put(115,5){$-\alpha_2$} 
\put(125,0){\line(1,0){20}}\put(150,0){\circle{5}}\put(145,5){$-\alpha_1$} 
\put(155,0){\line(1,0){20}}
\put(180,0){\circle*{5}}\put(175,8){$\sum_{1 \le i \le 6}\alpha_i$}
\put(185,0){\line(1,0){10}}\put(200,-3){- - -} 
\end{picture}
\vspace{0.5cm}

Then the twist by $\alpha_3$ induces an automorphism $T_3$ of 
$\mathfrak{t}^*_{\mathbf R}$. In particular, we have     
$$T_3(\alpha_3) = -\alpha_3, \;\;T_3( \alpha_6) = \sum_{1 \le i \le 6}\alpha_i, \;\; 
T_3(\alpha_9) = \alpha_9, \;\; T_3(\alpha_{13}) = \alpha_{13}$$
We define an automorphism $T_3^t$ of $\mathfrak{t}_{\mathbf R}$  
by $$\langle T_3^t(\phi), \cdot \rangle = \langle \phi, \: T_3(\cdot) \rangle, \; \phi \in 
\mathfrak{t}_{\mathbf R}$$
Since $T_3(\sum_{\alpha_i \in I}\mathbf{R}\alpha_i) = \sum_{\alpha_i \in I}\mathbf{R}\alpha_i$, 
the automorphism $T_3^t$ induces an automorphism of  $\mathrm{Hom}_{alg.gp}(L, \mathbf{C}^*)\otimes_{\mathbf Z}\mathbf{R} \subset \mathfrak{t}_{\mathbf R}$.  
If we write 
$$\mathfrak{t}^*_{\mathbf R}/\sum_{\alpha_i \in I}\alpha_i = \mathbf{R}\alpha_3 
\oplus \mathbf{R}\alpha_6 \oplus \mathbf{R}\alpha_9 \oplus \mathbf{R}\alpha_{13},$$
then $$\mathrm{Hom}_{alg.gp}(L, \mathbf{C}^*)\otimes_{\mathbf Z}\mathbf{R} = 
\mathbf{R}\alpha_3^* 
\oplus \mathbf{R}\alpha_6^* \oplus \mathbf{R}\alpha_9^* \oplus \mathbf{R}\alpha_{13}^*.$$ 
Then we have 
$$T_3^t(\alpha_3^*) = -\alpha_3^* + \alpha_6^*, \; T_3^t(\alpha_6^*) = \alpha_6^*, \; 
T_3^t(\alpha_9^*) = \alpha_9^*, \; T_3^t(\alpha_{13}^*) = \alpha_{13}^*.$$  

We next twist $Q$ by $\alpha_6$: 

\begin{picture}(400,20)
\put(0,-3){- - }\put(15,0){\line(1,0){10}}  
\put(30,0){\circle*{5}}\put(25,5){$\alpha_3$}\put(35,0){\line(1,0){20}}\put(60,0){\circle{5}}
\put(55,5){$\alpha_4$} 
\put(65,0){\line(1,0){20}}\put(90,0){\circle{5}}\put(85,5){$\alpha_5$}
\put(95,0){\line(1,0){20}}    
\put(120,0){\circle*{5}}\put(115,5){$\alpha_6$} 
\put(125,0){\line(1,0){20}}\put(150,0){\circle{5}}\put(145,5){$\alpha_7$} 
\put(155,0){\line(1,0){20}}
\put(180,0){\circle{5}}\put(175,5){$\alpha_8$}
\put(185,0){\line(1,0){20}}\put(205,0){\circle*{5}}\put(200,5){$\alpha_9$}\put(210,0){\line(1,0){10}}
\put(225,-3){- - -} 
\end{picture}  \vspace{1.0cm}

\begin{picture}(400,20) 
\put(0,-3){- - }\put(15,0){\line(1,0){10}}  
\put(30,0){\circle*{5}}\put(10,10){$\sum_{3 \le i \le 8}\alpha_i$}\put(35,0){\line(1,0){20}}\put(60,0){\circle{5}}
\put(55,-10){$-\alpha_8$} 
\put(65,0){\line(1,0){20}}\put(90,0){\circle{5}}\put(85,5){$-\alpha_7$}
\put(95,0){\line(1,0){20}}    
\put(120,0){\circle*{5}}\put(115,5){$-\alpha_6$} 
\put(125,0){\line(1,0){20}}\put(150,0){\circle{5}}\put(145,5){$-\alpha_5$} 
\put(155,0){\line(1,0){20}}
\put(180,0){\circle{5}}\put(175,5){$-\alpha_4$}
\put(185,0){\line(1,0){20}}\put(205,0){\circle*{5}}\put(200,10){$\sum_{4 \le i \le 9}\alpha_i$}\put(210,0){\line(1,0){10}}
\put(230,-3){- - -} 
\end{picture}  \vspace{1.0cm}
 
Let $T_6$ be the automorphism of $\mathfrak{t}^*_{\mathbf R}$ determined by the twist 
by $\alpha_6$. Then 
$$T_6(\alpha_3) = \sum_{3 \le i \le 8}\alpha_i, \; T_6(\alpha_6) = -\alpha_6,\; 
T_6(\alpha_9) = \sum_{4 \le i \le 9}\alpha_i,\; T_6(\alpha_{13}) = \alpha_{13}.$$ 
By using this we see that the automorphism $T^t_6$ of 
$\mathrm{Hom}_{alg.gp}(L, \mathbf{C}^*)\otimes_{\mathbf Z}\mathbf{R}$ is given by   
$$T^t_6(\alpha^*_3) = \alpha^*_3,\; T^t_6(\alpha^*_6) = \alpha^*_3 - \alpha^*_6 + 
\alpha^*_9, \; T^t_6(\alpha^*_9) = \alpha^*_9,\; T^t_6(\alpha^*_{13}) = \alpha^*_{13}.$$

We finally twist $Q$ by $\alpha_{13}$: 

\begin{picture}(400,40)
\put(0,-3){- - }\put(15,0){\line(1,0){10}}
\put(30,0){\circle*{5}}\put(25,5){$\alpha_9$}
\put(35,0){\line(1,0){20}}\put(60,0){\circle{5}}\put(55,5){$\alpha_{10}$} 
\put(65,0){\line(1,0){20}}\put(90,0){\circle{5}}\put(85,5){$\alpha_{11}$}
\put(95,0){\line(1,0){20}}\put(120,0){\circle{5}}\put(115,5){$\alpha_{12}$}
\put(125,0){\line(1,0){20}}\put(150,0){\circle*{5}}\put(145,5){$\alpha_{13}$} 
\put(155,0){\line(1,0){20}}
\put(180,0){\circle{5}}\put(175,5){$\alpha_{14}$}
\put(185,0){\line(1,0){20}}\put(210,0){\circle{5}}\put(205,5){$\alpha_{15}$} 
\put(215,0){\line(1,0){20}}\put(240,0){\circle{5}}\put(235,5){$\alpha_{16}$} 
\put(245,0){\line(1,0){20}}\put(270,0){\circle{5}}\put(265,5){$\alpha_{17}$} 
\put(275,0){\line(1,0){20}}\put(300,0){\circle{5}}\put(292,5){$\alpha_{18}$}
\put(305,0){\line(1,1){20}}\put(328,20){\circle{5}}\put(328,25){$\alpha_{19}$} 
\put(305,0){\line(1,-1){20}}\put(328,-20){\circle{5}}\put(330, -25){$\alpha_{20}$} 
\end{picture}
\vspace{1.0cm}

\begin{picture}(400,40)
\put(0,-3){- - }\put(15,0){\line(1,0){10}}
\put(30,0){\circle*{5}}\put(-10,-15){$\alpha_9 + 2\sum_{10 \le i \le 18}\alpha_i + \alpha_{19} 
+ \alpha_{20}$}
\put(35,0){\line(1,0){20}}\put(60,0){\circle{5}}\put(52,5){$-\alpha_{10}$} 
\put(65,0){\line(1,0){20}}\put(90,0){\circle{5}}\put(82,5){$-\alpha_{11}$}
\put(95,0){\line(1,0){20}}\put(120,0){\circle{5}}\put(112,5){$-\alpha_{12}$}
\put(125,0){\line(1,0){20}}\put(150,0){\circle*{5}}\put(142,5){$-\alpha_{13}$} 
\put(155,0){\line(1,0){20}}
\put(180,0){\circle{5}}\put(172,5){$-\alpha_{14}$}
\put(185,0){\line(1,0){20}}\put(210,0){\circle{5}}\put(202,5){$-\alpha_{15}$} 
\put(215,0){\line(1,0){20}}\put(240,0){\circle{5}}\put(232,5){$-\alpha_{16}$} 
\put(245,0){\line(1,0){20}}\put(270,0){\circle{5}}\put(262,5){$-\alpha_{17}$} 
\put(275,0){\line(1,0){20}}\put(300,0){\circle{5}}\put(310,0){$-\alpha_{18}$}
\put(305,0){\line(1,1){20}}\put(328,20){\circle{5}}\put(328,25){$-\alpha_{20}$} 
\put(305,0){\line(1,-1){20}}\put(328,-20){\circle{5}}\put(330, -25){$-\alpha_{19}$} 
\end{picture}
\vspace{1.5cm}

Let $T_{13}$ be the automorphism of $\mathfrak{t}^*_{\mathbf R}$ determined by the twist 
by $\alpha_{13}$. Then 
$$T_{13}(\alpha_3) = \alpha_3, \; T_{13}(\alpha_6) = \alpha_6,$$ 
$$T_{13}(\alpha_9) = \alpha_9 + 2\sum_{10 \le i \le 18}\alpha_i + \alpha_{19} 
+ \alpha_{20} \; T_{13}(\alpha_{13}) = -\alpha_{13}.$$ 
By using this we see that the automorphism $T^t_{13}$ of 
$\mathrm{Hom}_{alg.gp}(L, \mathbf{C}^*)\otimes_{\mathbf Z}\mathbf{R}$ is given by   
$$T^t_{13}(\alpha^*_3) = \alpha^*_3,\; T^t_{13}(\alpha^*_6) = \alpha^*_6, \; T^t_{13}(\alpha^*_9) = \alpha^*_9,\; T^t_{13}(\alpha^*_{13}) = 2\alpha^*_9 - \alpha^*_{13}.$$

We need one more twisting transformation to describe $W'$. In order to 
find the twisting, we must start with the parabolic subgroup $Q'$ obtained from 
$Q$ by the twist at $\alpha_9$. The parabolic subgroup $Q'$ corresponds to 
the marked Dynkin diagram 

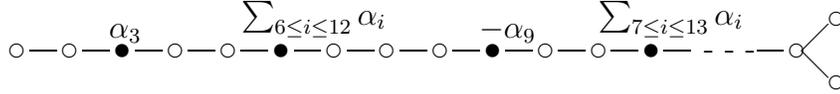
\begin{picture}(400,40)
\put(30,0){\circle{5}}\put(35,0){\line(1,0){10}}\put(50,0){\circle{5}} 
\put(55,0){\line(1,0){10}}\put(70,0){\circle*{5}}\put(65,5){$\alpha_3$}
\put(75,0){\line(1,0){10}}    
\put(90,0){\circle{5}}
\put(95,0){\line(1,0){10}}\put(110,0){\circle{5}}\put(115,0){\line(1,0){10}}
\put(130,0){\circle*{5}}\put(115,10){$\sum_{6 \le i \le 12}\alpha_i$}
\put(135,0){\line(1,0){10}}\put(150,0){\circle{5}}
\put(155,0){\line(1,0){10}}\put(170,0){\circle{5}}\put(175,0){\line(1,0){10}} 
\put(190,0){\circle{5}} 
\put(195,0){\line(1,0){10}}\put(210,0){\circle*{5}}\put(205,5){$-\alpha_9$} 
\put(215,0){\line(1,0){10}}\put(230,0){\circle{5}}\put(235,0){\line(1,0){10}} 
\put(250,0){\circle{5}}\put(255,0){\line(1,0){10}}\put(270,0){\circle*{5}}\put(250,10){$\sum_{7 \le i \le 13}\alpha_i$} 
\put(275,0){\line(1,0){10}}\put(290, -3){- - -}\put(310,0){\line(1,0){10}}
\put(325,0){\circle{5}}\put(327,0){\line(1,1){10}}\put(340,12){\circle{5}}
\put(327,0){\line(1,-1){10}}\put(340,-12){\circle{5}}
\end{picture}
\vspace{1.0cm}

We illustrate the marked Dynkin diagram in more details around the vertex 
$\sum_{7 \le i \le 13}\alpha_i$.  

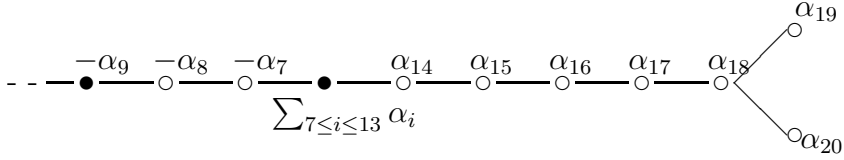
\begin{picture}(400,40)
\put(0,-3){- - }\put(15,0){\line(1,0){10}}
\put(30,0){\circle*{5}}\put(25,5){$-\alpha_9$}
\put(35,0){\line(1,0){20}}\put(60,0){\circle{5}}\put(55,5){$-\alpha_8$} 
\put(65,0){\line(1,0){20}}\put(90,0){\circle{5}}\put(85,5){$-\alpha_7$}
\put(95,0){\line(1,0){20}}\put(120,0){\circle*{5}}\put(100,-15){$\sum_{7 \le i \le 13}\alpha_i$}
\put(125,0){\line(1,0){20}}\put(150,0){\circle{5}}\put(145,5){$\alpha_{14}$} 
\put(155,0){\line(1,0){20}}
\put(180,0){\circle{5}}\put(175,5){$\alpha_{15}$}
\put(185,0){\line(1,0){20}}\put(210,0){\circle{5}}\put(205,5){$\alpha_{16}$} 
\put(215,0){\line(1,0){20}}\put(240,0){\circle{5}}\put(235,5){$\alpha_{17}$} 
\put(245,0){\line(1,0){20}}\put(270,0){\circle{5}}\put(265,5){$\alpha_{18}$} 
\put(275,0){\line(1,1){20}}\put(298,20){\circle{5}}\put(298,25){$\alpha_{19}$} 
\put(275,0){\line(1,-1){20}}\put(298,-20){\circle{5}}\put(300, -25){$\alpha_{20}$} 
\end{picture}
\vspace{1.0cm}

If we twist $Q'$ by $\sum_{7 \le i \le 13}\alpha_i$, then we get 

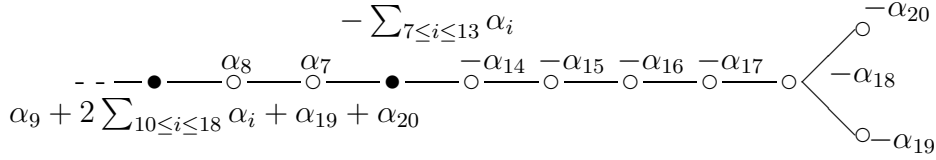
\begin{picture}(400,40)
\put(0,-3){- - }\put(15,0){\line(1,0){10}}
\put(30,0){\circle*{5}}\put(-25,-15){$\alpha_9 + 2\sum_{10 \le i \le 18}\alpha_i + \alpha_{19} 
+ \alpha_{20}$}
\put(35,0){\line(1,0){20}}\put(60,0){\circle{5}}\put(55,5){$\alpha_8$} 
\put(65,0){\line(1,0){20}}\put(90,0){\circle{5}}\put(85,5){$\alpha_7$}
\put(95,0){\line(1,0){20}}\put(120,0){\circle*{5}}\put(100,20){$-\sum_{7 \le i \le 13}\alpha_i$}
\put(125,0){\line(1,0){20}}\put(150,0){\circle{5}}\put(145,5){$-\alpha_{14}$} 
\put(155,0){\line(1,0){20}}
\put(180,0){\circle{5}}\put(175,5){$-\alpha_{15}$}
\put(185,0){\line(1,0){20}}\put(210,0){\circle{5}}\put(205,5){$-\alpha_{16}$} 
\put(215,0){\line(1,0){20}}\put(240,0){\circle{5}}\put(235,5){$-\alpha_{17}$} 
\put(245,0){\line(1,0){20}}\put(270,0){\circle{5}}\put(285,0){$-\alpha_{18}$} 
\put(275,0){\line(1,1){20}}\put(298,20){\circle{5}}\put(298,25){$-\alpha_{20}$} 
\put(275,0){\line(1,-1){20}}\put(298,-20){\circle{5}}\put(300, -25){$-\alpha_{19}$} 
\end{picture}
\vspace{1.5cm}

Let $T_{9 + 13}$ be the automorphism of $\mathfrak{t}^*_{\mathbf R}$ determined by the twist 
by $\sum_{7 \le i \le 13}\alpha_i$. Then 
$$T_{9 + 13}(\alpha_3) = \alpha_3, \; T_{9 + 13}(\sum_{6 \le i \le 12}\alpha_i) = \sum_{6 \le i \le 12}\alpha_i,$$
$$T_{9 + 13}(-\alpha_9) = \alpha_9 + 2\sum_{10 \le i \le 18}\alpha_i + \alpha_{19} 
+ \alpha_{20},\; T_{9 + 13}(\sum_{7 \le i \le 13}\alpha_i) = -\sum_{7 \le i \le 13}\alpha_i.$$ 
By using this we see that the automorphism $T^t_{9 + 13}$ of 
$\mathrm{Hom}_{alg.gp}(L, \mathbf{C}^*)\otimes_{\mathbf Z}\mathbf{R}$ is given by   
$$T^t_{9 + 13}(\alpha^*_3) = \alpha^*_3,\; T^t_{9+13}(\alpha^*_6) = \alpha^*_6, \; T^t_{9+13}(\alpha^*_9) = 2\alpha^*_6 - \alpha^*_9,\; T^t_{9+13}(\alpha^*_{13}) = 2\alpha^*_6 - 2\alpha^*_9 + \alpha^*_{13}.$$

We introduce here a new basis $\alpha^*_3$, $\alpha^*_6$, $\alpha^*_9$ and 
$\alpha^*_9 - \alpha^*_{13}$ on $\mathrm{Hom}_{alg.gp}(L, \mathbf{C}^*)\otimes_{\mathbf Z}\mathbf{R}$. With respect to this basis, we have the following matrix representation. 

$$T^t_3 = \left( \begin{array}{cccc} 
-1 & 0 & 0 & 0 \\ 
1 & 1 & 0 & 0 \\ 
0 & 0 & 1 & 0 \\
0 & 0 & 0 & 1 
\end{array}\right), \;\;\; 
T^t_6 = 
\left( \begin{array}{cccc} 
1 & 1 & 0 & 0 \\ 
0 & -1 & 0 & 0 \\ 
0 & 1 & 1 & 0 \\
0 & 0 & 0 & 1 
\end{array}\right), \;\;\; 
T^t_{9+13} = 
\left( \begin{array}{cccc} 
1 & 0 & 0 & 0 \\ 
0 & 1 & 2 & 0 \\ 
0 & 0 & -1 & 0 \\
0 & 0 & 0 & 1 
\end{array}\right)$$ and 
$$T^t_{13} = 
\left( \begin{array}{cccc} 
1 & 0 & 0 & 0 \\ 
0 & 1 & 0 & 0 \\ 
0 & 0 & 1 & 0 \\
0 & 0 & 0 & -1 
\end{array}\right).$$ 

One can easily read off that $$\langle T^t_3, \; T^t_6, \; T^t_{9+13}\rangle = W(B_3), \;\;   
\langle T^t_{13} \rangle = W(B_1), \;\; \langle T^t_3, \; T^t_6, \; T^t_{9+13}, T^t_{13} \rangle = W(B_3) \times W(B_1)$$

Let us compute how many {\bf Q}-factorial terminalizations $X$ has. 
$\mathcal{S}(L)$ contains 4 conjugacy classes. Two of them are those 
classes which respectively contain $Q$ and $Q'$.  
Other two classes are given by the marked Dynkin diagrams (cf. Remark 11):  

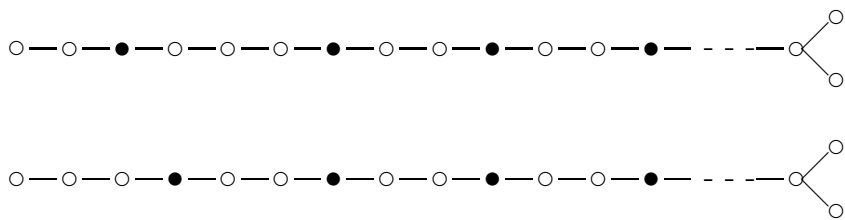
\begin{picture}(400,20)
\put(30,0){\circle{5}}\put(35,0){\line(1,0){10}}\put(50,0){\circle{5}} 
\put(55,0){\line(1,0){10}}\put(70,0){\circle*{5}}
\put(75,0){\line(1,0){10}}    
\put(90,0){\circle{5}}
\put(95,0){\line(1,0){10}}\put(110,0){\circle{5}}\put(115,0){\line(1,0){10}}
\put(130,0){\circle{5}}
\put(135,0){\line(1,0){10}}\put(150,0){\circle*{5}}
\put(155,0){\line(1,0){10}}\put(170,0){\circle{5}}\put(175,0){\line(1,0){10}} 
\put(190,0){\circle{5}} 
\put(195,0){\line(1,0){10}}\put(210,0){\circle*{5}}
\put(215,0){\line(1,0){10}}\put(230,0){\circle{5}}\put(235,0){\line(1,0){10}} 
\put(250,0){\circle{5}}\put(255,0){\line(1,0){10}}\put(270,0){\circle*{5}} 
\put(275,0){\line(1,0){10}}\put(290, -3){- - -}\put(310,0){\line(1,0){10}}
\put(325,0){\circle{5}}\put(327,0){\line(1,1){10}}\put(340,12){\circle{5}}
\put(327,0){\line(1,-1){10}}\put(340,-12){\circle{5}}
\end{picture}
\vspace{1.0cm} 

\begin{picture}(400,20)
\put(30,0){\circle{5}}\put(35,0){\line(1,0){10}}\put(50,0){\circle{5}} 
\put(55,0){\line(1,0){10}}\put(70,0){\circle{5}}
\put(75,0){\line(1,0){10}}    
\put(90,0){\circle*{5}}
\put(95,0){\line(1,0){10}}\put(110,0){\circle{5}}\put(115,0){\line(1,0){10}}
\put(130,0){\circle{5}}
\put(135,0){\line(1,0){10}}\put(150,0){\circle*{5}}
\put(155,0){\line(1,0){10}}\put(170,0){\circle{5}}\put(175,0){\line(1,0){10}} 
\put(190,0){\circle{5}} 
\put(195,0){\line(1,0){10}}\put(210,0){\circle*{5}}
\put(215,0){\line(1,0){10}}\put(230,0){\circle{5}}\put(235,0){\line(1,0){10}} 
\put(250,0){\circle{5}}\put(255,0){\line(1,0){10}}\put(270,0){\circle*{5}} 
\put(275,0){\line(1,0){10}}\put(290, -3){- - -}\put(310,0){\line(1,0){10}}
\put(325,0){\circle{5}}\put(327,0){\line(1,1){10}}\put(340,12){\circle{5}}
\put(327,0){\line(1,-1){10}}\put(340,-12){\circle{5}}
\end{picture}
\vspace{1.0cm}

Therefore $$\sharp \mathcal{S}(L)  = 4\cdot \vert W' \vert 
= 4 \cdot 96 = 384$$
Note that $X_{[3^3,2^2,1]} \to \bar{O}_{[3^3,2^2,1]}$ is a double covering.  
Recall that 
there is an exact sequence 
$$1 \to W_X \to W' \stackrel{\bar{\rho}_X}\to \mathrm{Aut}(X/\bar{O})/
\mathrm{Aut}(X_{[3^3,2^2,1]}/\bar{O}_{[3^3,2^2,1]}) (= \mathbf{Z}/2\mathbf{Z})$$
Let $\omega_3$ (resp. $\omega_6$, $\omega_{9 + 13}$, $\omega_{13}$) 
be an element of $W'$ determined by the twist $T^t_3$ (resp. $T^t_6$, 
$T^t_{9+13}$, $T^t_{13}$). Then $$\bar{\rho}_X(\omega_3) = 
\bar{\rho}_X(\omega_6) = \bar{\rho}_X(\omega_{9+13}) = 0, \;\; 
\bar{\rho}_X(\omega_{13}) = \bar{1}$$ This means that 
$$W_X = W(B_3).$$ 
As a consequence, we have $$\sharp \{\mathbf{Q}\mathrm{-factorial}\; \mathrm{terminalizations}\: \mathrm{of}\: X\} = \sharp \mathcal{S}(L)/ \vert W_X \vert 
= 8.$$  $\square$  \vspace{0.2cm}

\begin{center} 
Research Institute for Mathematical Science, Kyoto University, Japan

namikawa@kurims.kyoto-u.ac.jp 
\end{center}

\end{document}